\documentclass[11pt]{article}
\usepackage{benstyle}
 
\usepackage{enumitem}
 
\usepackage[nocompress]{cite}

\usepackage[font=small]{caption}

\newcommand\two{0.35\textwidth}

\newcommand\four{0.225\textwidth}

\title{Frame approximation with bounded coefficients}
\author{Ben Adcock \& Mohsen Seifi \\ Department of Mathematics \\ Simon Fraser University \\ Canada \\ \url{ben_adcock@sfu.ca}, \url{smseifi@sfu.ca}}

\begin{document}

\maketitle

\begin{abstract}
Due to their flexibility, frames of Hilbert spaces are attractive alternatives to bases in approximation schemes for problems where identifying a basis is not straightforward or even feasible. Computing a best approximation using frames, however, can be challenging since it requires solving an ill-conditioned linear system. One consequence of this ill-conditioning is that the coefficients of {such a frame approximation} can grow large. In this paper we resolve this issue by introducing two methods for frame approximation that possess bounded coefficients. As we show, these methods typically lead to little or no deterioration in the approximation accuracy, but successfully avoid the large coefficients inherent to previous approaches, thus making them attractive in situations where large coefficients are undesirable. We also present theoretical analysis to support these conclusions.
\end{abstract}

\noindent
\textbf{Keywords:} frames, function approximation, ill-conditioning, singular value decomposition

\vspace{1pc}
\noindent
\textbf{AMS subject classifications:} 42C15, 42C30, 41A10, 65T40

\section{Introduction}

In scientific computing, it is often convenient to use systems of functions that are non-orthogonal and near-redundant {(i.e.\ close to linearly dependent)} as part of an approximation scheme. Orthogonal or well-conditioned bases can be difficult to construct in certain problems, so relaxing this requirement can offer substantial flexibility. In a series of recent works \cite{BADHframespart,BADHFramesPart2,PEHighDim}, a framework for approximation in certain redundant systems of functions arising as so-called \textit{frames} of Hilbert spaces has been developed. 
Frames can provide viable alternatives in approximation problems where it can be challenging or possibly infeasible to devise good orthonormal bases. Examples include spectral approximations on irregular domains in one or more dimensions, approximation of structured functions using orthonormal bases augmented with a finite collection of feature functions, and concatenation of orthonormal bases to approximate function decomposable as sums of functions well approximated in different systems.

A challenge when computing an approximation in near-redundant systems is dealing with the ill-conditioning of the linear algebra problem to be solved to obtain the coefficients of the approximation. In particular, {the first $N$ elements of an infinite frame always lead to an ill-conditioned linear systems (for large enough $N$)}, and this ill-conditioning can be arbitrarily bad \cite{BADHframespart}. Fortunately, as proposed in \cite{BADHframespart,BADHFramesPart2,PEHighDim}, this ill-conditioning can be addressed through regularization, for instance, Truncated Singular Value Decompositions (TSVD). The result is a well-conditioned approximation with provable error guarantees.

However, the coefficients of such an approximation can be large. In the problem considered in \cite{BADHframespart} they can be as large as $\ord{1/\sqrt{\epsilon}}$, where $\epsilon$ is the SVD truncation parameter if a TSVD is used. If one uses a black-box solver such as \texttt{backslash} in Matlab, then the same applies with $\epsilon = \epsilon_{\mathrm{mach}}$ (see also Fig.\ \ref{fig:introfig2}). Worse still, in the method of \cite{BADHFramesPart2} (which delivers better limiting accuracy) the coefficients can grow as large as $\ord{1/\epsilon}$. Notably, this occurs in the pre-asymptotic regime in $N$ (the number of terms in the approximation), when the approximation error is still moderate. As $N \rightarrow \infty$, the coefficients become $\ord{1}$ in magnitude. 

With this issue in mind, in this paper we revisit the topic of frame approximations, and introduce two new schemes which compute frame approximations while maintaining bounded, or at worst slowly growing coefficients. As we show, these schemes lead to small-norm coefficients with typically little deterioration in the accuracy of the approximation over the TSVD approach of \cite{BADHframespart,BADHFramesPart2}, or black-box approaches such as \texttt{backslash}.

\subsection{Background}

To describe our main contribution, we first recap the setup of \cite{BADHframespart}.
A \textit{frame} of a separable Hilbert space $\rH$ is a countable system $\Phi = \{ \phi_n \}_{n \in I} \subset \rH$ satisfying the so-called \textit{frame condition}
\be{
\label{ineq : Frame condition}
A\nm{f}^{2} \leq \sum_{n \in I} \left\vert \left\langle f, \phi_{n} \right\rangle \right\vert^{2} \leq B \nm{f}^{2} \quad \forall f\in \rH,
}
for constants $0 < A \leq B < \infty$, referred to as the \textit{frame bounds}. While \R{ineq : Frame condition} implies that $\mathrm{span}(\Phi)$ is dense in $\rH$, a frame is generally not a basis. Indeed, there may exist nonzero coefficients $\bm{x} = \{ x_{n} \}_{n \in I} \in \ell^2(I)$ for which the sum $\sum_{n \in I} x_{n} \phi_{n}$ converges in $\rH$ and satisfies $\sum_{n \in I} x_{n} \phi_{n} = 0$. This is often referred to as \textit{overcompleteness} or \textit{redundancy} of the frame.


The concern of \cite{BADHframespart} and this paper is computations using $N$ frame elements $\Phi_N = \{ \phi_n \}_{n \in I_N}$, where $I_N \subset I$ is an index set of cardinality $N$. Computing the best approximation $\cP_N f = \sum_{n \in I_N} x_n \phi_n$ to $f \in \rH$ from the subspace $\rH_N = \spn(\Phi_N)$ equates to solving the linear system
\be{
\label{GNsystemintro}
\bm{G}_N \bm{x} = \bm{y},
}
for the coefficients $\bm{x} = \{ x_n \}_{n \in I_N}$, where $\bm{y} = \{ \ip{f}{\phi_n} \}_{n \in I_N}$ and $\bm{G}_N = \{ \ip{\phi_n}{\phi_m} \}_{m,n \in I_N}$ is the Gram matrix of $\Phi_N$. This system is generally ill-conditioned for large $N$.  In \cite{BADHframespart,BADHFramesPart2} the authors regularize this system using the TSVD of $\bm{G}_N$ with a truncation parameter $\epsilon > 0$. This leads to a regularized solution $\bm{x}^{\epsilon}$ of \R{GNsystemintro} and an approximation $\cP^{\epsilon}_{N} f = \sum^{N}_{n=1} x^{\epsilon}_n \phi_n$ to $f$.  It was shown that this TSVD approximation satisfies
\be{
\label{TSVDerr}
\nmu{f - \cP^{\epsilon}_{N} f } \leq E_{N,\epsilon}(f),\qquad E_{N,\epsilon}(f) = \inf \left \{ \nmu{f - \cT_N \bm{z} } + \sqrt{\epsilon} \nm{\bm{z}} : \bm{z} \in \bbC^N \right \},
}
where $\cT_N : \bbC^N \rightarrow \rH_N,\ \bm{z} = \{z_n \}_{n \in I_N} \mapsto \sum_{n \in I_N} z_n \phi_n$.
Notice that when $\epsilon = 0$, this bound simply expresses the fact that $\cP^0_N f \equiv \cP_N f$ is the best approximation to $f$ from $\rH_N$. For $\epsilon > 0$, this bound determines the effect of discarding the small singular values. It asserts that the accuracy of the approximation $\cP^{\epsilon}_N f$ {for fixed $N$} depends on how well $f$ can be approximated by an element $\sum_{n \in I_N} z_n \phi_n$ with coefficients that are not too large. It also implies that
\be{
\label{TVSDlimit}
\limsup_{N \rightarrow \infty}\nmu{f - \cP^{\epsilon}_{N} f } \leq \sqrt{\epsilon/A} \nm{f},
}
i.e.\ the limiting accuracy is with $\ord{\sqrt{\epsilon}}$ of $f$. This can be improved to $\ord{\epsilon}$ via the approach of \cite{BADHFramesPart2} -- see also Remark \ref{r:FNA2setup}.

\subsection{{New contributions}}

\begin{figure}[t]
\begin{center}
\small{
\begin{tabular}{cc}
\includegraphics[width=\two]{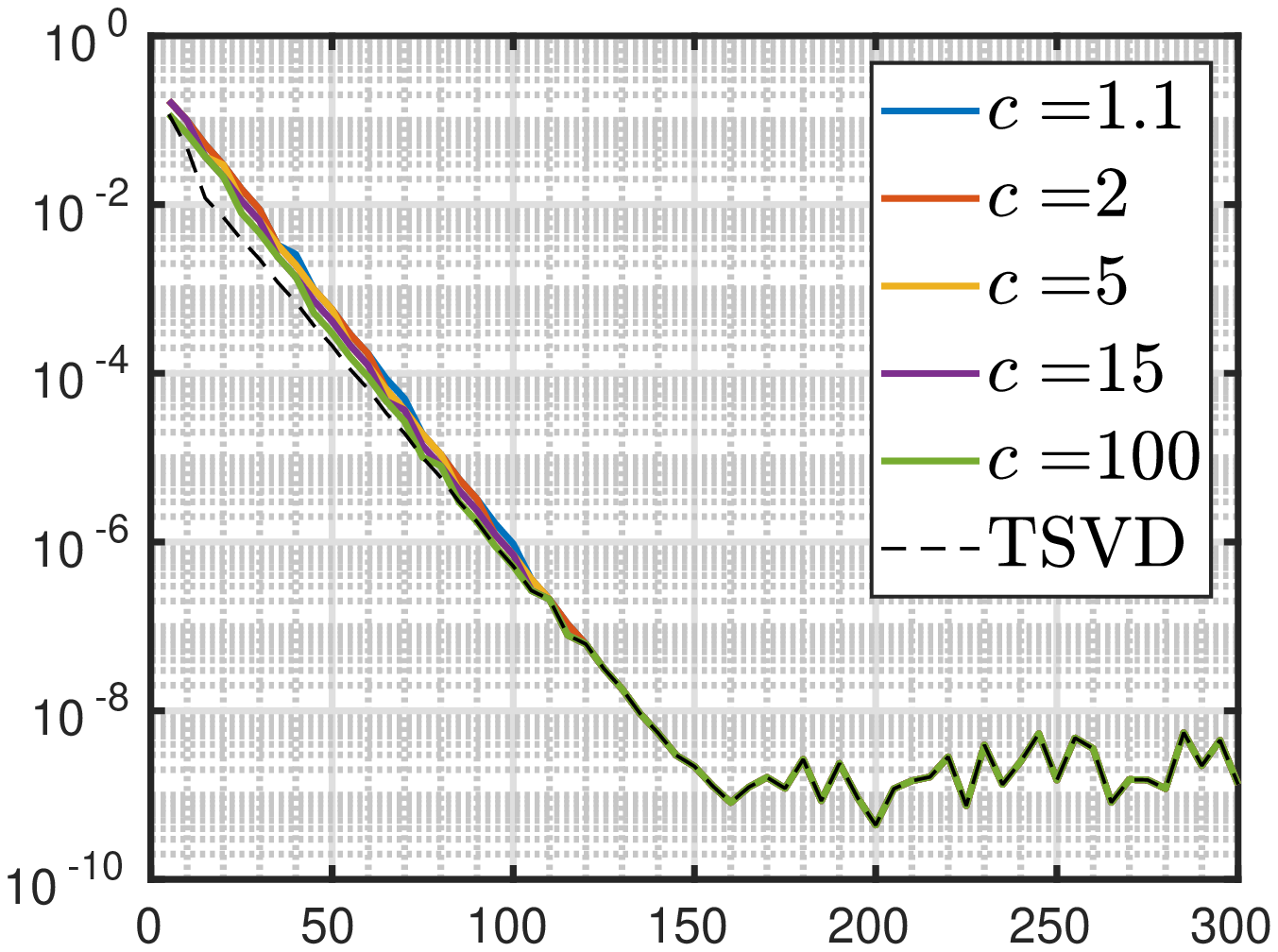}&
\includegraphics[width=\two]{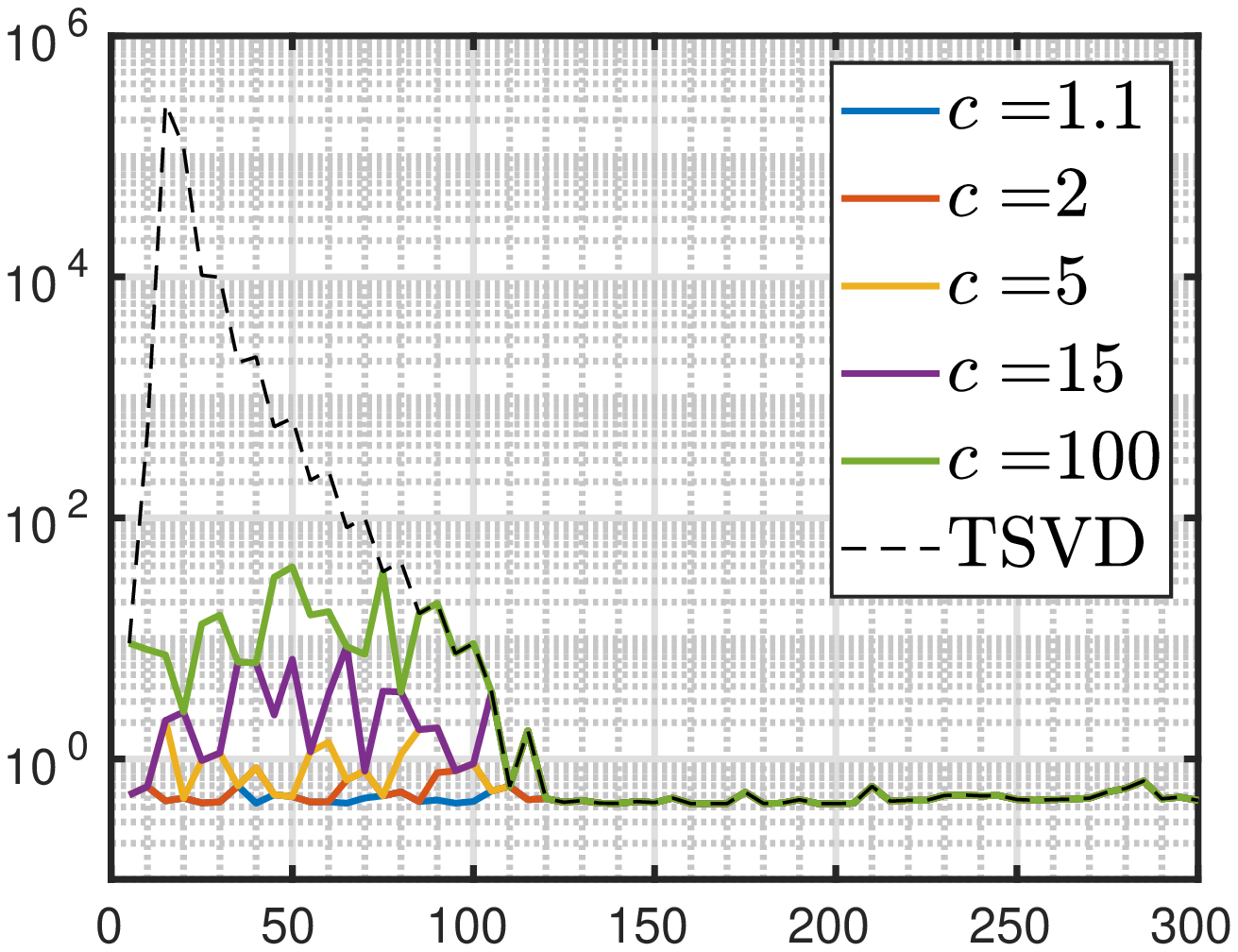}
\end{tabular}
}
\end{center}
\vspace*{-5mm}
\caption{
Error (left) and coefficient norm (right) versus $N$ for approximating $f(t) = \frac{1}{1+75 t^2}$ in the interval $[-1/2,1/2]$ using the \textit{Legendre polynomial} frame (as described in \S \ref{ss:mainexamp}, this is the frame constructed by restricting the orthonormal Legendre polynomials on $[-1,1]$ to $[-1/2,1/2]$).  `TSVD' is the TSVD approximation $\cP^{\epsilon}_N f$. The other approximations are the ASVD1 approximation $\cA^{\epsilon,c}_N f$ introduced in this paper with different values of $c$. The threshold $\epsilon$ is set to be $10^{-15}$ in all cases.
}
\label{fig:introfig}
\end{figure}

In \cite{BADHframespart} it was also shown that the coefficients $\bm{x}^{\epsilon}$ of the TSVD approximation behave roughly like the approximation error divided by $\sqrt{\epsilon}$. Specifically,
\be{
\label{TSVDcoeff}
\nm{\bm{x}^{\epsilon}} \leq E_{N,\epsilon}(f) / \sqrt{\epsilon}.
}
Therefore, while  the coefficients are $\ord{1}$ in the limit as $N \rightarrow \infty$, specifically,
\be{
\label{TSVDcoefflimit}
\limsup_{N \rightarrow \infty} \nm{\bm{x}^{\epsilon}}  \leq \nm{f}/\sqrt{A}, 
}
they may behave at worst like $1/\sqrt{\epsilon}$ before the onset of this asymptotic behaviour. 
Fig.\ \ref{fig:introfig} gives a typical example of this phenomenon. For the frame in question, the coefficient norm first increases exponentially fast, before decaying to $\ord{1}$ in the limit. Notice that the approximation errors and coefficient norm obey the rough relationship described in \R{TSVDcoeff}: namely their ratio is approximately $\sqrt{\epsilon}$ for all $N$.

\begin{figure}[t]
\begin{center}
\small{
\begin{tabular}{cc}
\includegraphics[width=\two]{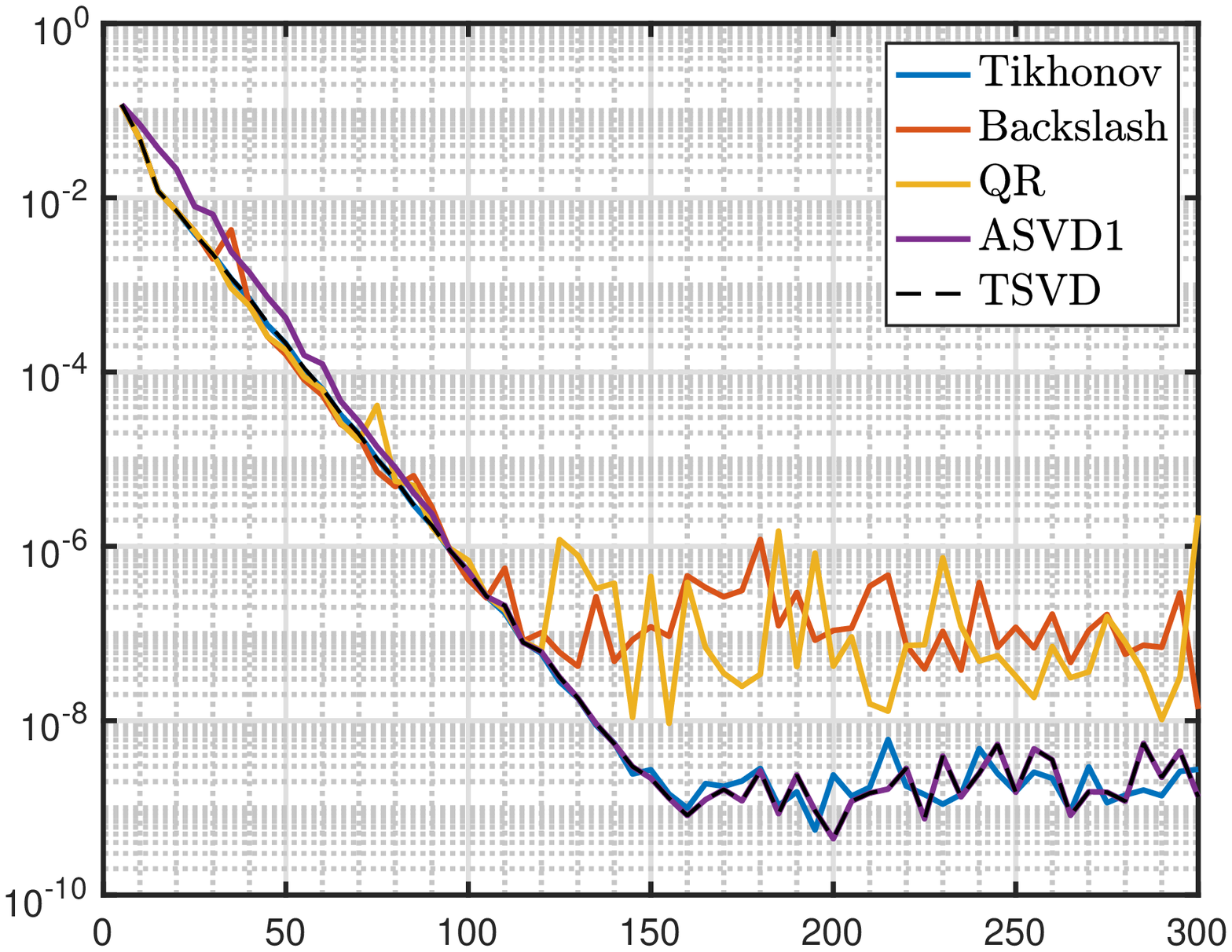}&
\includegraphics[width=\two]{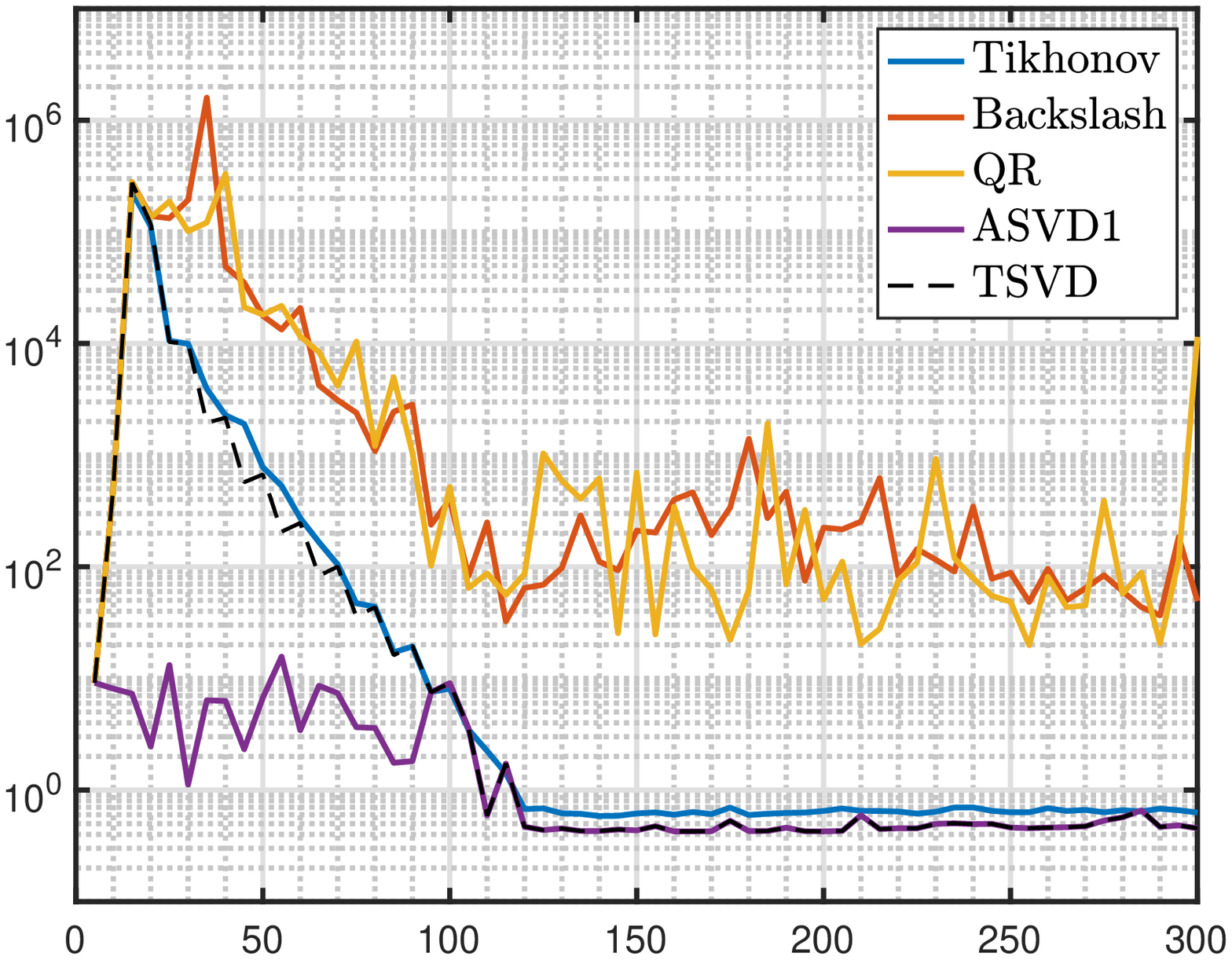} \\
\includegraphics[width=\two]{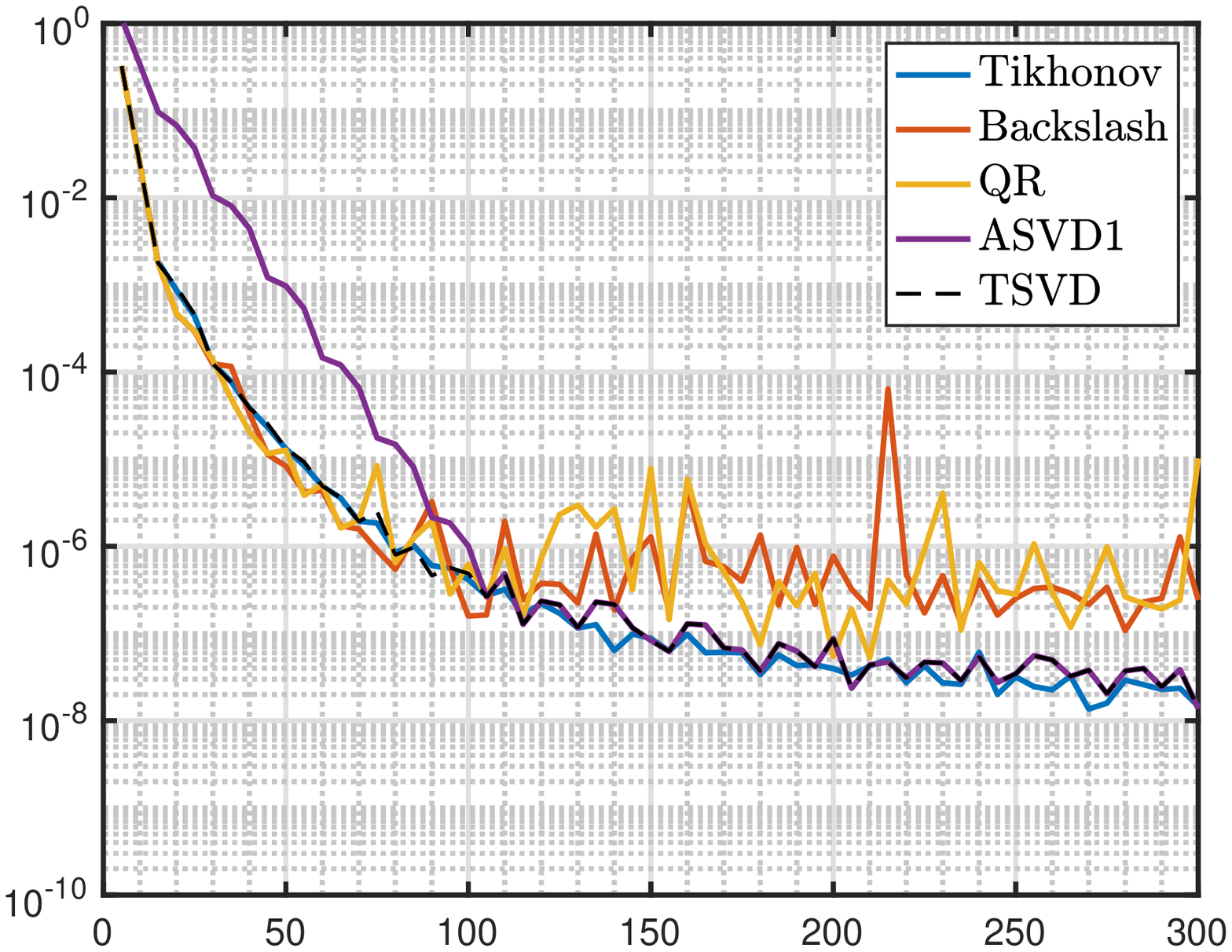}&
\includegraphics[width=\two]{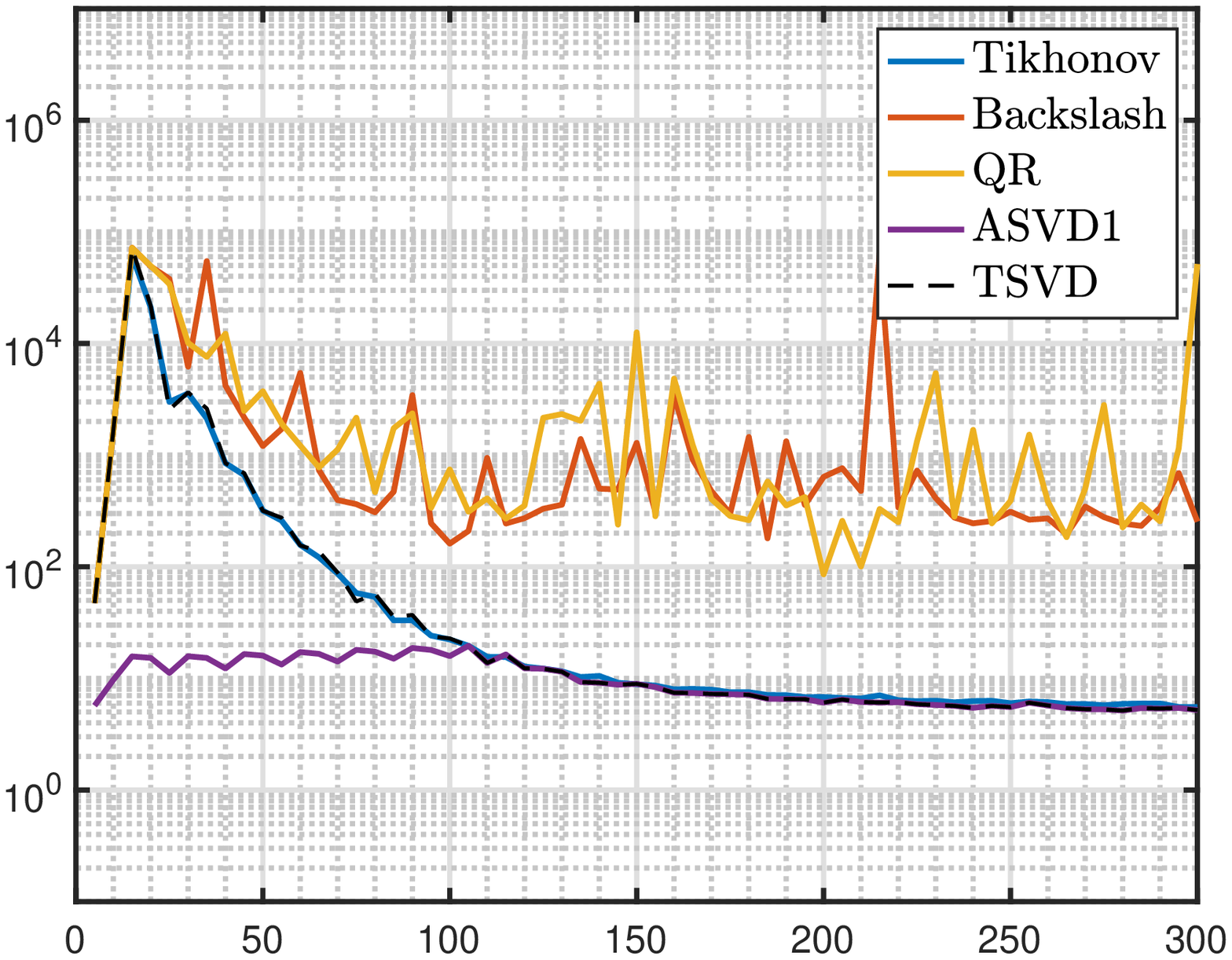}
\end{tabular}
}
\end{center}
\vspace*{-5mm}
\caption{
Error (left) and coefficient norm (right) versus $N$ for approximating $f(t) = \frac{1}{1+75 t^2}$ (top row) and $f(t) = \frac{1}{0.57 - t}$ (bottom row) in the interval $[-1/2,1/2]$ using the Legendre polynomial frame (see Fig.\ \ref{fig:introfig}). `TSVD' is the TSVD approximation $\cP^{\epsilon}_N f$ with $\epsilon = 10^{-15}$. `Backslash' solves \R{GNsystemintro} using Matlab's \texttt{backslash} command. `QR' solves \R{GNsystemintro} by computing the QR factorization of $\bm{G}_N$. `Tikhonov' applies Tikhonov regularization to \R{GNsystemintro}, solving $\min_{\bm{z} \in \bbC^N} \nm{\bm{G}_N \bm{z} - \bm{y}}^2 + \lambda^2 \nm{\bm{z}}^2$. The parameter $\lambda$ is set as $\lambda = 10^{-15}$. `ASVD1' is the ASVD1 approximation $\cA^{\epsilon,c}_N f$ of this paper with $c = 15$ and $\epsilon = 10^{-15}$. 
}
\label{fig:introfig2}
\end{figure}

This phenomenon is not specific to the TSVD approximation. The same behaviour (in fact, sometimes somewhat worse) is seen for Tikhonov regularization, as well as QR factorization and Matlab's \texttt{backslash}. See Fig.\ \ref{fig:introfig2}. It is interesting to note that in the limit $N \rightarrow \infty$ the TSVD and Tikhonov approximations (which are quite closely related) obtain both smaller errors and smaller coefficient values than QR or \texttt{backslash}.

This situation is bizarre. When $f$ is poorly approximated, i.e.\ $E_{N,\epsilon}(f) \gg \sqrt{\epsilon}$, the computed approximation has large coefficients. Yet when $f$ is well approximated, i.e.\ the error is close to $\sqrt{\epsilon}$, the approximation has small coefficients. Surely the approximation can be modified so that the coefficients remain small in the former regime? Given the redundancy of frames, there may be many ways to achieve such a relatively poor approximation without having to endure large coefficients.

In this paper, we solve this problem by modifying the TSVD approach. Rather than simply thresholding all singular values below the cutoff $\epsilon$,  the idea is to adaptively choose the singular values based on the function being approximated. Specifically, we identify those singular values that are of insignificant importance in approximating $f$, but are responsible for growth of norm of the solution $\bm{x}$. In our first method, \textit{Adaptive SVD 1} (\textit{ASVD1}), we discard a singular value $\sigma_{n}$ of the truncated Gram matrix $\bm{G}_N$ if either $\sigma_{n} \leq \epsilon$ or
\be{ 
\label{ineq : ASVD1 constraint}
\frac{\left\vert \left\langle \bm{y}, \bm{v}_{n} \right\rangle \right \vert}{\sigma_{n}} > c \nm{\bm{y}}, 
} 
where $\bm{v}_{n}$ is the corresponding $n^{\text{th}}$ singular vector of $\bm{G}_{N}$ and $c>0$ is a parameter. Imposing constraints of the type (\ref{ineq : ASVD1 constraint}) leads to regularized approximations for which coefficients are at worst of magnitude of $\mathcal{O}(\sqrt{N})$. While this bound is independent of the choice of the parameter $\epsilon$, it still depends on the degree of the approximation. Our second method, \textit{Adaptive SVD 2} (\textit{ASVD2}), tackles this issue by imposing further regularization and leads to solutions whose norms are always bounded above by $c\nm{\bm{y}}$. See \S\ref{section : algorithms} for further details of this method. 

We analyze both schemes, ASVD1 and ASVD2, from the point of view of the approximation error and coefficient norm.  For the ASVD1 approximation $\cA^{\epsilon,c}_{N}f $ we show that
\be{
\label{ASVD1_err_intro}
\Vert f-\cA^{\epsilon,c}_{N} f\Vert\leq\Vert f-\mathcal{T}_{N}\bm{z}\Vert+\max\left\{\sqrt{\epsilon}, \frac{\Vert f-\mathcal{T}_{N}\bm{z}\Vert}{c\nm{\bm{y}}-\Vert\bm{z}\Vert}\right\}\Vert\bm{z}\Vert,\qquad\forall\bm{z}\in\mathbb{C}^{N}, \ \nm{\bm{z}} < c \nm{\bm{y}},
}
and if $\bm{x}^{\epsilon} = ( x^{\epsilon}_n)^{N}_{n=1}$ are the coefficients of $\cA^{\epsilon,c}_N f = \sum^{N}_{n=1} x^{\epsilon}_n \phi_n$,
\bes{
\nm{\bm{x}^{\epsilon}} \leq \min \left \{ c \sqrt{B N} \nm{f} , E_{N,\epsilon}(f) / \sqrt{\epsilon} \right \}.
}
We also prove similar results for the ASVD2 scheme. At first glance, \R{ASVD1_err_intro} may seem somewhat restrictive in comparison to the bound \R{TSVDerr} for the TSVD approximation. Our analysis and numerical results show that this is not the case, provided $c$ is not too small. Indeed, since $\liminf_{N \rightarrow \infty} \nm{\bm{y}} \geq \sqrt{A} \nm{f}$, the bound \R{ASVD1_err_intro} effectively states that the accuracy of the ASVD1 approximation is determined by how well $f$ can be approximated using coefficients $\bm{z}$ whose norm is at most $c \sqrt{A}$ times bigger than $\nm{f}$. For the Legendre polynomial frame example (used in Figs.\ \ref{fig:introfig} and \ref{fig:introfig2}) we discuss how such coefficients arise in practice.
We also examine the asymptotic behaviour of the schemes, and show that the ASVD1 and ASVD2 approximations behave like the TSVD approximation as $N \rightarrow \infty$. Our numerical examples show that ASVD1 and ASVD2 can often provide nearly as good approximations as TSVD without undergoing a regime of large coefficients. An example of this was already shown in Fig.\ \ref{fig:introfig}.

\subsection{Remarks}

The line of work \cite{BADHFramesPart2,PEHighDim} on numerical frame approximation originated with the study of Fourier extensions \cite{boyd2005fourier}, which are motivated by embedding methods for PDEs \cite{pasquettiFourEmbed,lui2009embedding,boffi2015immersed,kasolis2015fictitious,shirokoff2015volumepenalty}. More recently, fast algorithms have been developed \cite{LyonFast,matthysen2015fastfe,matthysen2017fastfe2d}, as well as approaches suitable for higher-dimensional problems \cite{PEHighDim}.
Frames are well-known tools in modern signal and image processing, coding theory and sampling theory. While they are less well-known in numerical analysis, there are numerous methods that rely on approximation systems that are nearly redundant. See \cite{BADHframespart} and references therein. Imposing the structure of a frame bestows problem with several pleasant properties. In particular, as discussed in \S \ref{s:preliminaries}, it guarantees the existence of at least one set of coefficients with small norm (the so-called \textit{frame coefficients}). However, it is worth noting that the results in \cite{BADHframespart,BADHFramesPart2} pertaining to the error and coefficient norm of the TSVD approximation do not require a frame structure.

{The same is true for the corresponding results proved} in this paper. Yet, by designing approximations with small-norm coefficients, we are implicitly assuming that the system admits such representations. This will not be the case in general redundant systems. For instance, the monomial basis $\phi_n(t) = t^{n-1}$ {for, say, $L^2(0,1)$,} is nearly redundant for large $N$, but most functions cannot be approximated  using coefficients with small norm. {Even within the frame setting, as we will see, if the constant $c$ is chosen too small, then ASVD1 and ASVD2 perform poorly.} On the other hand, if the system is known to admit accurate, small-norm coefficient approximations {and $c$ is chosen reasonably}, then ASVD1 and ASVD2 are viable approaches. {However, since the question of whether or not this holds is, in the absence of a frame, highly dependent on the system, we shall continue to work exclusively with frames in this paper.}

Finally, we remark that one could consider alternatives to hard singular value thresholding as we do in this paper, for example, Tikhonov regularization (see Fig.\ \ref{fig:introfig2}) and the closely-related constrained least squares. However, thresholded SVDs are rather simpler to interpret and analyze, and it is not clear what benefits Tikhonov regularization might convey.


\section{Preliminaries}\label{s:preliminaries}

We commence with some background material on frames. For more in-depth overviews, see \cite{christensen2003introduction}. In what follows, we assume that $\Phi = \{ \phi_{n} \}_{n \in I}$ is an indexed family in a separable Hilbert space $\rH$ over $\mathbb{C}$, where $I$ is a countable index set. We write $\left\langle \cdot, \cdot \right\rangle$ and $\nm{\cdot}$ for the inner product and the induced norm on $\rH$ respectively.  We write $\ell^2(I)$ for the space of square-summable complex sequences indexed by $I$, and denote its norm by $\nm{\cdot}$, i.e.\ $\nm{\bm{x}} = \sqrt{\sum_{n \in I} \abs{x_{n}}^{2}}$ for $\bm{x} = \{ x_n \}_{n \in I}$.

\subsection{Frames}

As observed, $\Phi$ is a frame if it satisfies the frame condition \R{ineq : Frame condition}. If $A = B$, then $\Phi$ is a \textit{tight} frame. Associated with any frame $\Phi$ is the \textit{synthesis operator}
\bes{
\mathcal{T} : \ell^{2}(I) \rightarrow \rH, \quad \bm{y} = \{ y_{n} \}_{n \in I} \mapsto \sum_{n \in I} y_{n}\phi_{n}.
}
Its adjoint, the \textit{analysis operator}, is given by 
\bes{
\mathcal{T}^{*} : \rH \rightarrow \ell^{2}(I), \quad f \mapsto \{ \langle f, \phi_{n} \rangle \}_{n \in I},
}
and the composition $\mathcal{S} = \mathcal{T}\mathcal{T}^{*}$, known as the \textit{frame operator}, is
\bes{
\mathcal{S} : \rH \rightarrow \rH, \quad f \mapsto \sum_{n \in I} \langle f, \phi_{n} \rangle\phi_{n}. 
}
The frame operator is self-adjoint, invertible, bounded and positive. If $\Phi$ is tight then $\cS = A \cI$.
%

A \textit{dual frame} $\Psi = \{ \psi_n \}_{n \in I} \subset \rH$ of $\Phi$ is an indexed family that is frame for $\rH$ and satisfies
\be{\label{eq : Dual frame}
f = \sum_{n \in I} \left\langle f, \psi_{n} \right\rangle \phi_{n} = \sum_{n \in I} \left\langle f, \phi_{n} \right\rangle \psi_{n}, \quad \forall f\in \rH.
}
The \textit{canonical dual frame} is the frame $\Psi = \{  \cS^{-1} \phi_n \}_{n \in I} \subset \rH$. For this frame, \R{eq : Dual frame} gives
\be{ \label{eq: Canonical dual frame}
f = \sum_{n \in I} \left\langle f, \mathcal{S}^{-1}\phi_{n} \right\rangle \phi_{n} = \sum_{n \in I} \left\langle \mathcal{S}^{-1}f, \phi_{n} \right\rangle \phi_{n}.
}
The canonical dual frame has frame bounds $1/B$ and $1/A$, namely
\be{
\label{dualframebd}
1/B \nm{f}^{2} \leq \sum_{n \in I} \vert \left\langle f, \mathcal{S}^{-1}\phi_{n} \right\rangle \vert^{2} \leq 1/A \nm{f}^{2}, \quad \forall f\in \rH. 
} 
Notice that \R{eq: Canonical dual frame} gives $f = \sum_{n \in I} a_n \phi_n$, where $\bm{a} = \left\{ \left\langle f, \mathcal{S}^{-1}\phi_{n} \right\rangle \right\}_{n \in I}$. These are the \textit{frame coefficients} of $f$. They have minimal $\ell^{2}$-norm amongst all representations of $f$ in $\Phi$. In other words, if $f = \sum_{n \in I} a_{n}\phi_{n} = \sum_{n \in I} c_{n}\phi_{n}$ for some $\bm{c} = \left\{ c_{n} \right\}_{n \in I} \in \ell^{2}(I)$, then $\nm{\bm{c}} \geq \nm{\bm{a}}$. {These coefficients will play a role in our later analysis. However, we note in passing that approximating $f$ by $\sum_{n \in I_N} a_n \phi_n$, while it has the advantage of bounded coefficients, does not generally yield a good approximation for finite $N$ \cite{BADHframespart}.}


\subsection{Truncation of frames and best approximations}
For each $N \in \mathbb{N}$, let $I_N \subset I$ be index sets satisfying $|I_N| = N$ and
\bes{
I_1 \subset I_2 \subset \cdots ,\qquad \bigcup_{N = 1}^{\infty} I_{N} = I.
}
Let $\Phi_{N} = \left\{ \phi_{n} \right\}_{n \in I_{N}}$ be the finite set of frame elements with indices in $I_N$ and define the finite-dimensional subspace $\rH_N = \spn(\Phi_N)$. Notice that $\Phi_N$ is a frame for $\rH_N$. We let $\mathcal{T}_{N} : \mathbb{C}^{N} \rightarrow \rH_{N}$, $\mathcal{T}^{*}_{N} : \rH_{N} \rightarrow \mathbb{C}^{N}$, and $\mathcal{S}_{N} : \rH_{N} \rightarrow \rH_{N}$ denote the truncated synthesis, analysis and frame operators respectively.

%

We consider approximating an element $f \in \rH$ from $\rH_N$. Since $\rH$ is a Hilbert space, the best approximation is the orthogonal projection onto $\rH_N$. Write $\cP_N : \rH \rightarrow \rH_N$ for the orthogonal projection operator.  Using $\Phi_N$, we can write this projection as $\cP_N f = \sum_{n \in I_N} x_n \phi_n = \cT_N \bm{x}$, where the coefficients $\bm{x} = \{ x_n \}_{n \in I_N }$ are a solution of the linear system 
\be{ \label{eq: Best approximation}
\bm{G}_{N} \bm{x} = \bm{y}, \qquad \bm{y} = \left\{ \left\langle f, \phi_{n} \right\rangle \right\}_{n \in I_{N}},
} 
where $\bm{G}_N = \{ \ip{\phi_n}{\phi_m} \}_{m,n \in I_N }$ is the Gram matrix of $\Phi_N$. Since $\Phi_N$ may be linearly dependent, this linear system may be singular. In and of itself, this is of little consequence, since the projection $\cP_N f$ itself is unique. However, even when $\Phi_N$ is linearly independent, $\bm{G}_N$ is ill-conditioned for large $N$ unless $\Phi$ happens to be a Riesz basis \cite{BADHframespart}.


\subsection{Truncated SVD frame approximation}

To deal with this ill-conditioning, \cite{BADHframespart,BADHFramesPart2} used Truncated Singular Value Decompositions (TSVDs). Since $\bm{G}_{N}$ is nonnegative definite its singular values $\sigma_n$, $n \in I_N$, are its eigenvalues and its SVD takes the form $\bm{G}_{N}=\bm{V}\bm{\Sigma}\bm{V}^{*}$,
where $\bm{V}\in\mathbb{C}^{N\times N}$ is unitary and $\bm{\Sigma}=\diag\left(\sigma_n\right)_{n \in I_N}$. Write $\{\bm{v}_{n}:n\in I_{N}\}$ for the columns of $\bm{V}$, which are left/right singular vectors of $\bm{G}_{N}$, i.e. $\bm{G}_{N}\bm{v}_{n}=\sigma_{n}\bm{v}_{n}, n\in I_{N}$.
Given a threshold $\epsilon > 0$, define the regularized solution $\bm{x}^{\epsilon}$ of \R{eq: Best approximation} as
\bes{
\bm{x}^{\epsilon} = (\bm{G}^{\epsilon}_N)^{\dag} \bm{y},\qquad \bm{G}^{\epsilon}_{N} = \bm{V} \bm{\Sigma}^{\epsilon} \bm{V}^*,
}
where $\dag$ denotes the pseudoinverse and $\bm{\Sigma}^{\epsilon}$ has $n^{\rth}$ diagonal entry equal to $\sigma_n$ if $\sigma_n > \epsilon$ and zero otherwise. Applying $\cT_N$, this gives rise to a regularized approximation, defined as
\bes{
\cP^{\epsilon}_N f = \cT_N \bm{x}^{\epsilon} .
}

\rem{
\label{r:FNA2setup}
This approximation satisfies the bounds \R{TSVDerr} and \R{TVSDlimit}. In particular, the limiting accuracy is $\ordu{\sqrt{\epsilon}}$. To attain $\ord{\epsilon}$ accuracy, one can replace the square system \R{eq: Best approximation}, corresponding to a finite section \cite{lindner2006}, by a least-squares system $\bm{G}_{M,N} \bm{x} \approx \bm{y}$, where $\bm{G}_{M,N} = \{ \ip{\phi_n}{\phi_m} \}_{m \in I_M,n \in I_N }$ is the tall Gram matrix of size $M \times N$ (a so-called \textit{uneven} section) and $\bm{y} = \{ \ip{f}{\phi_m} \}_{m \in I_M}$. As shown in \cite{BADHFramesPart2}, if $M$ is chosen sufficiently large in relation to $N$, one achieves $\ord{\epsilon}$ accuracy in the limit $N \rightarrow \infty$.

This approach also addresses the case where it undesirable to compute $\bm{G}_{M,N}$ and $\bm{y}$ (typically their entries involve evaluating integrals). If $\bm{G}_{M,N}$ is replaced by the matrix $\bm{G}_{M,N} =  \{ \ell_{m,M}(\phi_n) \}_{m \in J_M,n \in I_N} \}$ and $\bm{y}$ by $\bm{y} =  \{ \ell_{m,M}(f) \}_{m \in J_M}$, where the $\ell_{m,M}$ are linear `measurement' functionals -- for example, pointwise samples -- then the results of \cite{BADHFramesPart2} guarantee $\ord{\epsilon}$ accuracy in the limit under suitable conditions on the functionals $\ell_{m,M}$, $M$ and $N$.

In this paper, for our analysis we consider the simpler setup outlined above. However, we expect our theoretical results can be extended to setting of \cite{BADHFramesPart2} with additional effort.
}

\section{Adaptive SVD frame approximations}
\label{section : algorithms}

As observed, the coefficients $\bm{x}^{\epsilon}$ of the TSVD approximation can initially grow large, a behaviour that can be traced to certain small singular values that are retained in the truncation. We now develop new methods that avoid this behaviour.

\subsection{A general thresholded SVD approximation}

Before introducing our two main methods, we first define general thresholded SVD approximation.  Let $\Lambda\subseteq I_{N}$ be a subset of singular values to be used in the approximation and $\bm{\Sigma}_{\Lambda}$ be the diagonal matrix with $n^{\text{th}}$ entry $\sigma_{n}$ if $n\in\Lambda$ and zero otherwise, and set
\be{
\bm{G}_{\Lambda}=\bm{V}\bm{\Sigma}_{\Lambda}\bm{V}^{*}.
}
Then we define the approximation
\bes{
\mathcal{P}_{\Lambda}f=\mathcal{T}_{N} \bm{x}_{\Lambda} ,\qquad \bm{x}^{\Lambda}=\left(\bm{G}_{\Lambda}\right)^{\dagger}\bm{y}=\bm{V}\left(\bm{\Sigma}_
{\Lambda}\right)^{\dagger}\bm{V}^{*}\bm{y}.
}
Notice that the orthogonal projection $\cP_N f$ and the TSVD approximation $\cP^{\epsilon}_N f$ are special cases of $\cP_{\Lambda} f$, corresponding to $\Lambda = I_N$ and  $\Lambda = \Lambda^{\epsilon} = \{ n \in I_N : \sigma_n > \epsilon \}$ respectively.

\subsection{Projections}

Since it will be useful later, we now show that $\cP_{\Lambda} : \rH \rightarrow \rH$ is an orthogonal projection for any fixed $\Lambda$. First, to each singular vector $\bm{v}_{n}$ we associate an element $\xi_{n} \in \rH_{N}$ defined by
\bes{
\xi_{n}=\mathcal{T}_{N}\bm{v}_{n}=\sum_{m\in I_{N}}\left(\bm{v}_{n}\right)_{m}\phi_{m}.
}
Notice that the $\xi_{n}$ are orthogonal in $\rH$ with 
\be{
\label{prolate_orthog}
\langle\xi_{n},\xi_{m}\rangle=\langle\mathcal{T}_{N}\bm{v}_{n},\mathcal{T}_{N}\bm{v}_{m}\rangle=\langle\bm{v}_{n},\mathcal{T}_{N}^{*}\mathcal{T}_{N}\bm{v}_{m}\rangle=\sigma_{m}\langle\bm{v}_{n},\bm{v}_{m}\rangle=\sigma_{m}\delta_{n,m},
}
and hence form an orthogonal basis for $\rH_{N}$. Next, notice that we may write
\be{
\label{coef_formula}
\bm{x}=\sum_{n\in I_{N}}\frac{\langle\bm{y},\bm{v}_{n}\rangle}{\sigma_{n}}\bm{v}_{n},\quad\bm{x}_{\Lambda}=\sum_{n\in \Lambda}\frac{\langle\bm{y},\bm{v}_{n}\rangle}{\sigma_{n}}\bm{v}_{n},
}
where $\langle\cdot,\cdot\rangle$ denotes the Euclidean inner product on $\mathbb{C}^{N}$, and $\langle\bm{y},\bm{v}_{n}\rangle=\sum_{m\in I_{N}}\langle f,\phi_{m}\rangle\overline{\left(\bm{v}_{n}\right)_{m}}=\langle f,\xi_{n}\rangle$.
This gives
\be{
\label{orthog_proj}
\mathcal{P}_{N}f = \cT_N \bm{x} =\sum_{n\in I_{N}}\frac{\langle f,\xi_{n}\rangle}{\sigma_{n}}\xi_{n},\qquad\mathcal{P}_{\Lambda}f= \cT_N \bm{x}_{\Lambda} = \sum_{n\in \Lambda}\frac{\langle f,\xi_{n}\rangle}{\sigma_{n}}\xi_{n}.
}
From this, we deduce the claim. The operator $\cP_{\Lambda}$ is precisely the orthogonal projection onto the subspace $\rH_{\Lambda} = \spn \{ \xi_n : n \in \Lambda \}$.


\subsection{Methods}

For any fixed $\Lambda$, the expression \R{coef_formula} gives $\nm{\bm{x}_{\Lambda}}^2 = \sum_{n \in \Lambda} \frac{|\ip{\bm{y}}{\bm{v}_n} |^2}{\sigma^2_n}$.
Since our goal is to maintain small-norm coefficients, this immediately leads to a possible strategy for choosing $\Lambda$: we simply exclude any $n$ for which the term $|\ip{\bm{y}}{\bm{v}_n} |/\sigma_n$ exceeds some tolerance.  To determine a suitable tolerance, we recall \R{TSVDcoefflimit}.
This motivates choosing the tolerance as a multiple of $\nm{f}$. Of course, $\nm{f}$ is  generally unknown.  However, notice that
\bes{
\nm{\bm{y}}^2 = \sum_{n \in I_N} | \ip{f}{\phi_n} |^2 \rightarrow \sum_{n \in I} | \ip{f}{\phi_n} |^2 \geq A \nm{f}^2,\qquad N \rightarrow \infty,
}
where $\bm{y}$ is the right-hand side of the linear system \R{eq: Best approximation}. Hence, we take the tolerance to be some multiple of $\nm{\bm{y}}$ instead. This leads to the following:

\defn{
\label{d:ASVD1}
The \textit{Adaptive SVD 1 (ASVD1)} frame approximation corresponds to the choice
\bes{
\Lambda= \Lambda^{\epsilon,c}= \Lambda^{\epsilon,c}(f) = \left \{  n \in I_N : \sigma_{n} > \epsilon ,\ \frac{\vert\langle\bm{y}, \bm{v}_{n}\rangle\vert}{\sigma_{n}}\leq c\nm{\bm{y}} \right \},
}
where $\epsilon > 0$ is some threshold and $c > 0$ is a constant.
}
Notice that this definition of $\Lambda$ still retains the truncation of all small singular values below a threshold $\epsilon > 0$. We will discuss the rationale for this in \S \ref{ss:numexamp}. 

Due to the definition of $\Lambda$, ASVD1 leads to coefficients satisfying
$
\nm{\bm{x}_{\Lambda}} \leq c \sqrt{N} \nm{\bm{y}} \leq c \sqrt{B} \sqrt{N} \nm{f},
$
where the second inequality follows from the frame condition \R{ineq : Frame condition}.
Thus, the coefficients may not be truly bounded independently of $N$, although, as we show in Theorem \ref{thm:limsup}, they will be $\ord{1}$ in the limit $N \rightarrow \infty$. In our second method we avoid this behaviour:

\defn{
\label{d:ASVD2}
The \textit{Adaptive SVD 2 (ASVD2)} frame approximation corresponds to the choice $\Lambda = \Lambda^{\epsilon,c}$, where $\Lambda$ is the largest set satisfying
 \bes{
\Lambda \subseteq \left \{ n \in I_N : \sigma_n > \epsilon \right \} ,\qquad \sqrt{\sum_{n\in \Lambda}\frac{|\langle\bm{y},\bm{v}_{n}\rangle|^2}{\sigma^2_{n}}} \leq c \nm{\bm{y}},
 }
 and where $\epsilon > 0$ is some threshold and $c > 0$ is a constant.
}

{Note that computing such a set $\Lambda$ is easily done by sorting the terms in the sum from smallest to largest.}
By definition, in this method we have $\nm{\bm{x}_{\Lambda}} \leq c \nm{\bm{y}}$. Since $\nm{\bm{y}} \leq \sqrt{B} \nm{f}$ due to \R{ineq : Frame condition}, the coefficients are uniformly bounded in $N$.

\section{Examples}\label{ss:numexamp}

We now present several examples comparing ASVD1 and ASVD2 to the TSVD approximation.

\subsection{The Legendre polynomial frame}\label{ss:mainexamp}

In our first example, we consider an example of a general frame construction where an orthonormal basis is restricted to a subdomain. Such frames are particularly useful for approximating functions irregular-shaped domains, where finding orthonormal bases can be challenging \cite{PEHighDim}.
For simplicity, we consider the one-dimensional case, where $\rH = L^2(-1/2,1/2)$ is the space of square-integrable functions on the interval $(-1/2,1/2)$.  Let
\bes{
\psi_n(t) = \sqrt{n+1/2} P_n(t),\qquad n \in I : = \bbN_0,
}
be the orthonormal Legendre basis of $L^2(-1,1)$, where $P_n$ is the classical Legendre polynomial with normalization $P_n(1) = 1$, and define $\phi_n = \psi_n |_{(-1/2,1/2)}$. The family $\{ \phi_n \}^{\infty}_{n=0}$ forms a tight (in fact, Parseval) and {overcomplete frame} for $\rH$ with frame bounds $A = B = 1$.

Before showing numerical examples with this frame, we briefly discuss its approximation properties for smooth functions. Let $I_N = \{0,\ldots,N-1\}$ and $H^k(-1/2,1/2)$ be the standard Sobolev space of order $k \in \bbN$ with norm $\nm{\cdot}_{H^k}$. Then, simplifying and slightly modifying the proof of \cite[Thm.\ 5.2]{PEHighDim}, one can show that if $f \in H^k(-1/2,1/2)$ then, for each $N$, there exists a vector of coefficients $\bm{z} \in \bbC^N$ such that
\be{
\label{LegFrameSob}
\nm{f - \cT_N \bm{z}} \leq C_k N^{-k} \nm{f}_{H^k},\qquad \nm{\bm{z}} \leq C_k \nm{f}_{H^k},
}
where $C_k > 0$ is a constant depending on $k$ only.  Hence, the error term in \R{TSVDerr} satisfies
\bes{
E_{N,\epsilon}(f) \leq \min_{0 \leq l \leq k} \left \{ C_l  ( N^{-l} + \sqrt{\epsilon}  ) \nm{f}_{H^l} \right \}.
}
Therefore, for smooth functions $f$ TSVD approximation error is guaranteed decrease rapidly down to the level $\ord{\sqrt{\epsilon}}$. An example of this behaviour was shown previously Fig.\ \ref{fig:introfig}.

\begin{figure}
\begin{center}
\small{
\begin{tabular}{cccc}
\includegraphics[width=\four]{f_1_IG_FNA1_ASVD1_err.eps} &
\includegraphics[width=\four]{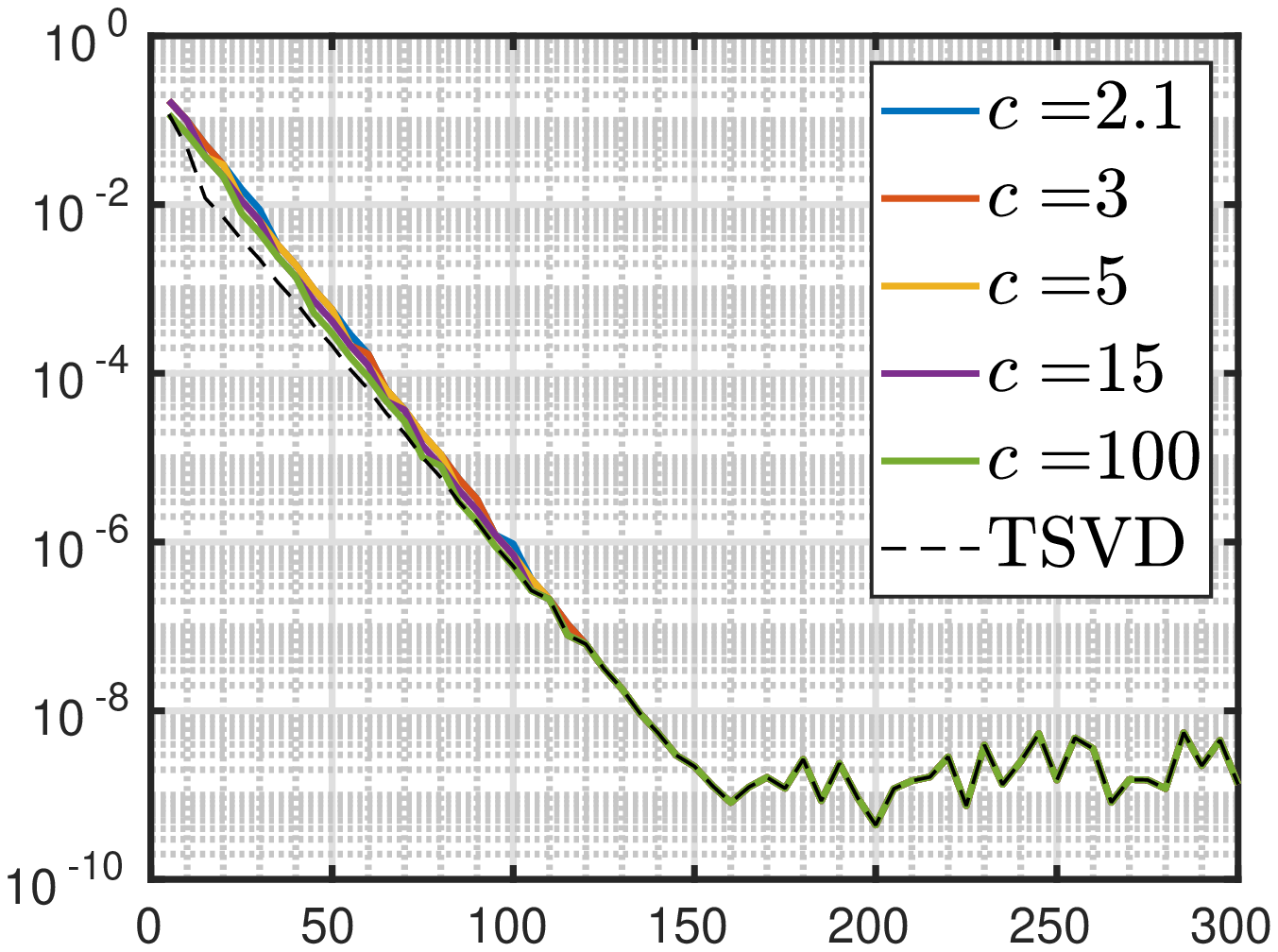}&
\includegraphics[width=\four]{f_1_IG_FNA1_ASVD1_coef.eps}&
\includegraphics[width=\four]{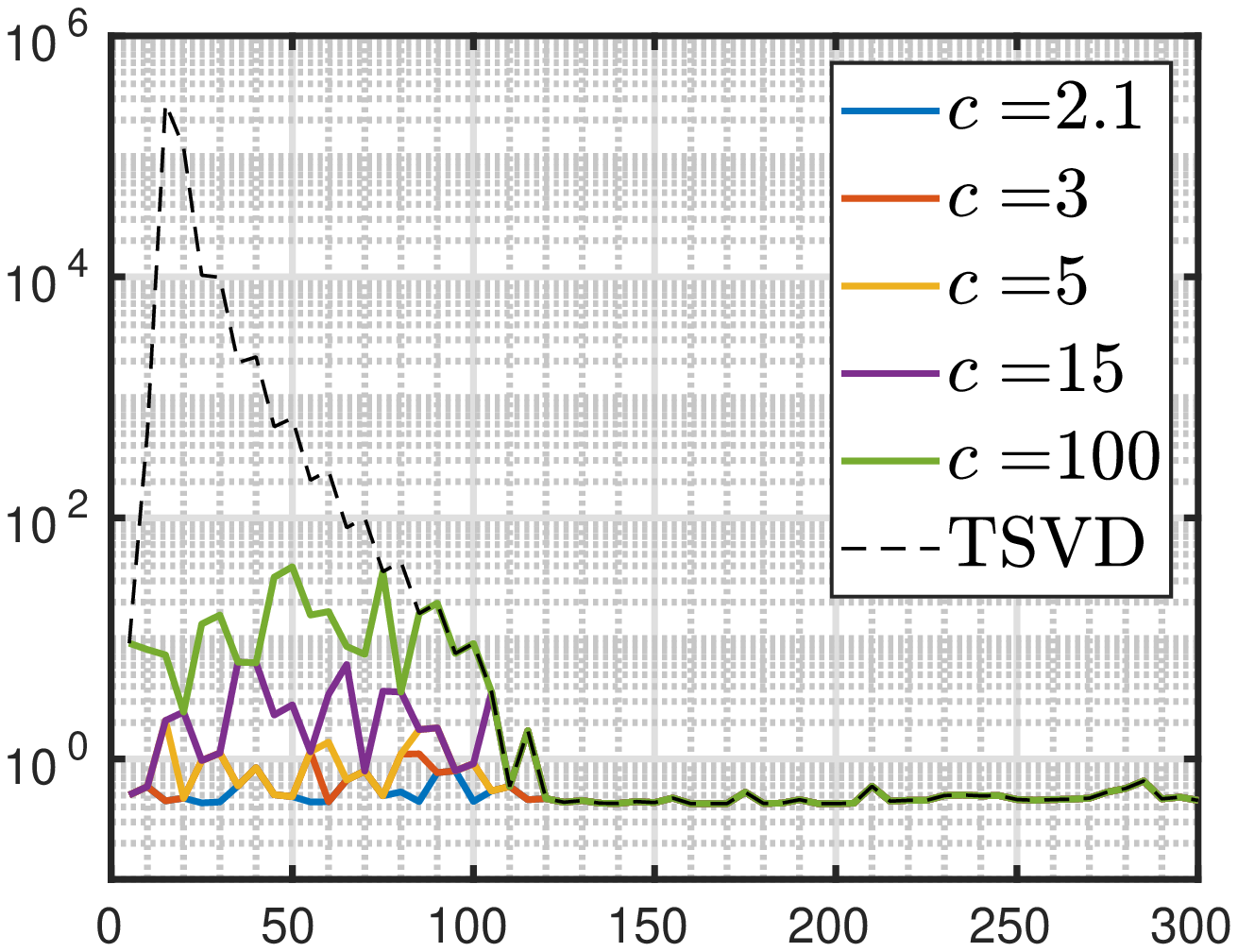}
\\
\includegraphics[width=\four]{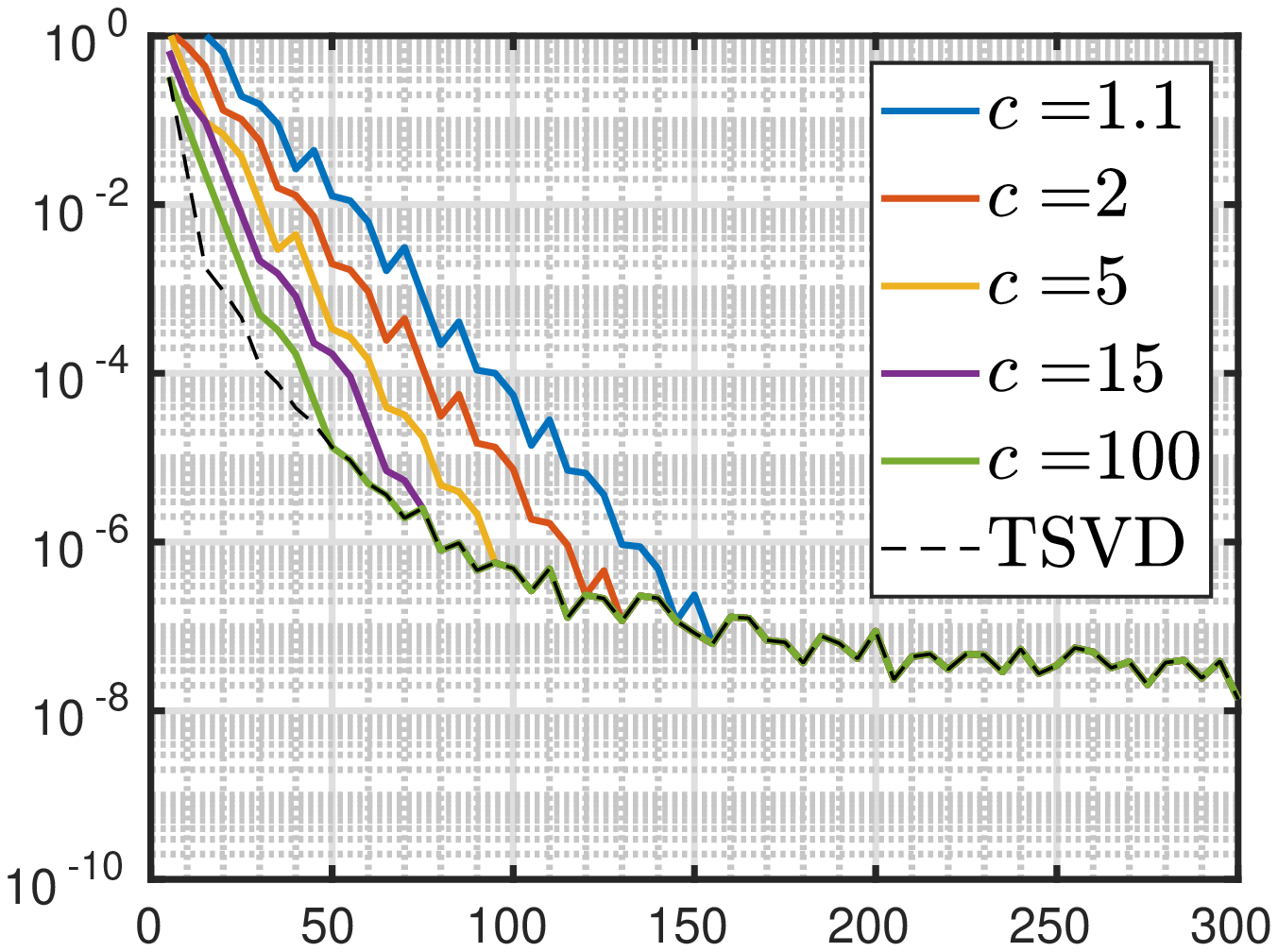} &
\includegraphics[width=\four]{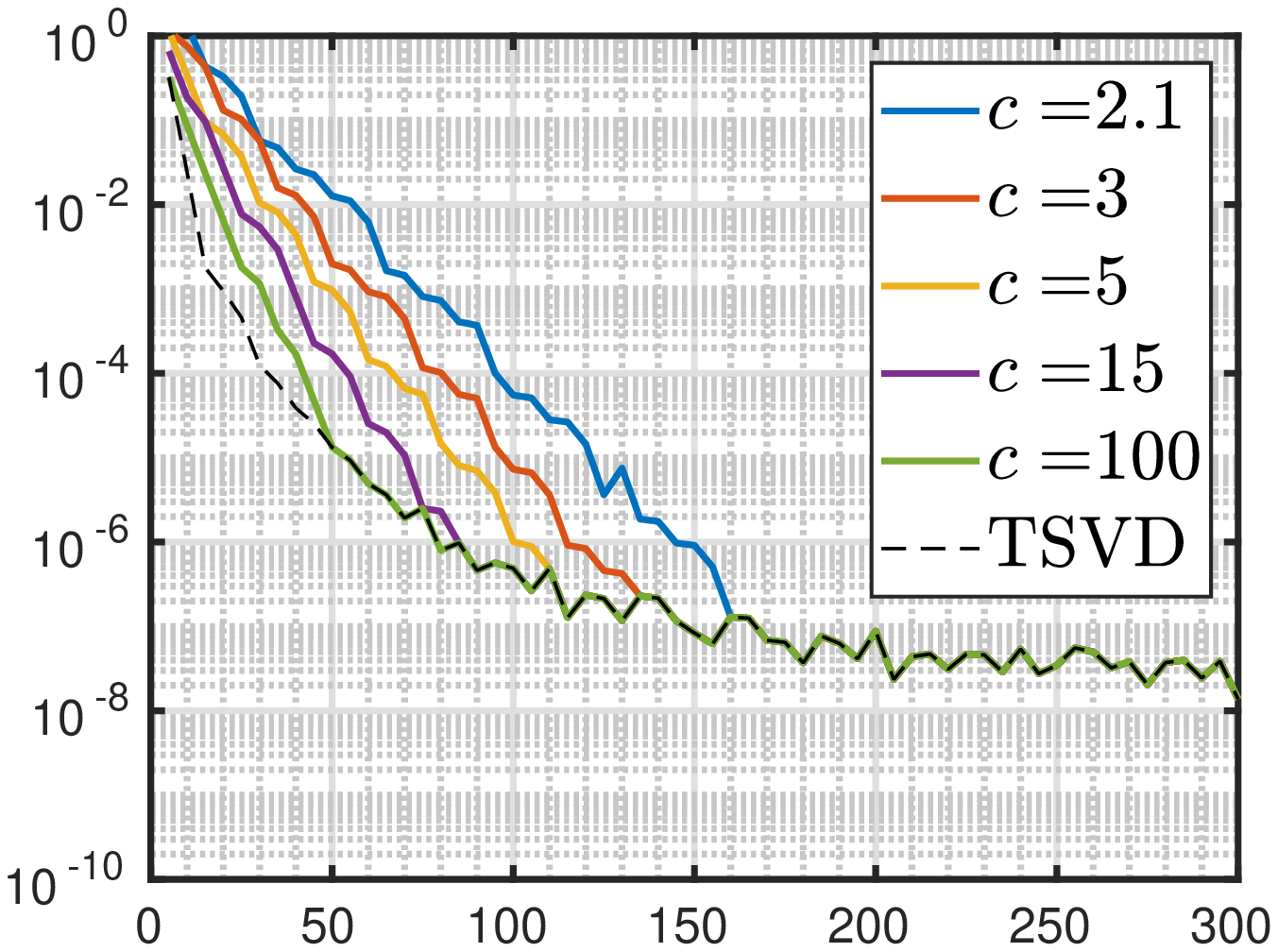}&
\includegraphics[width=\four]{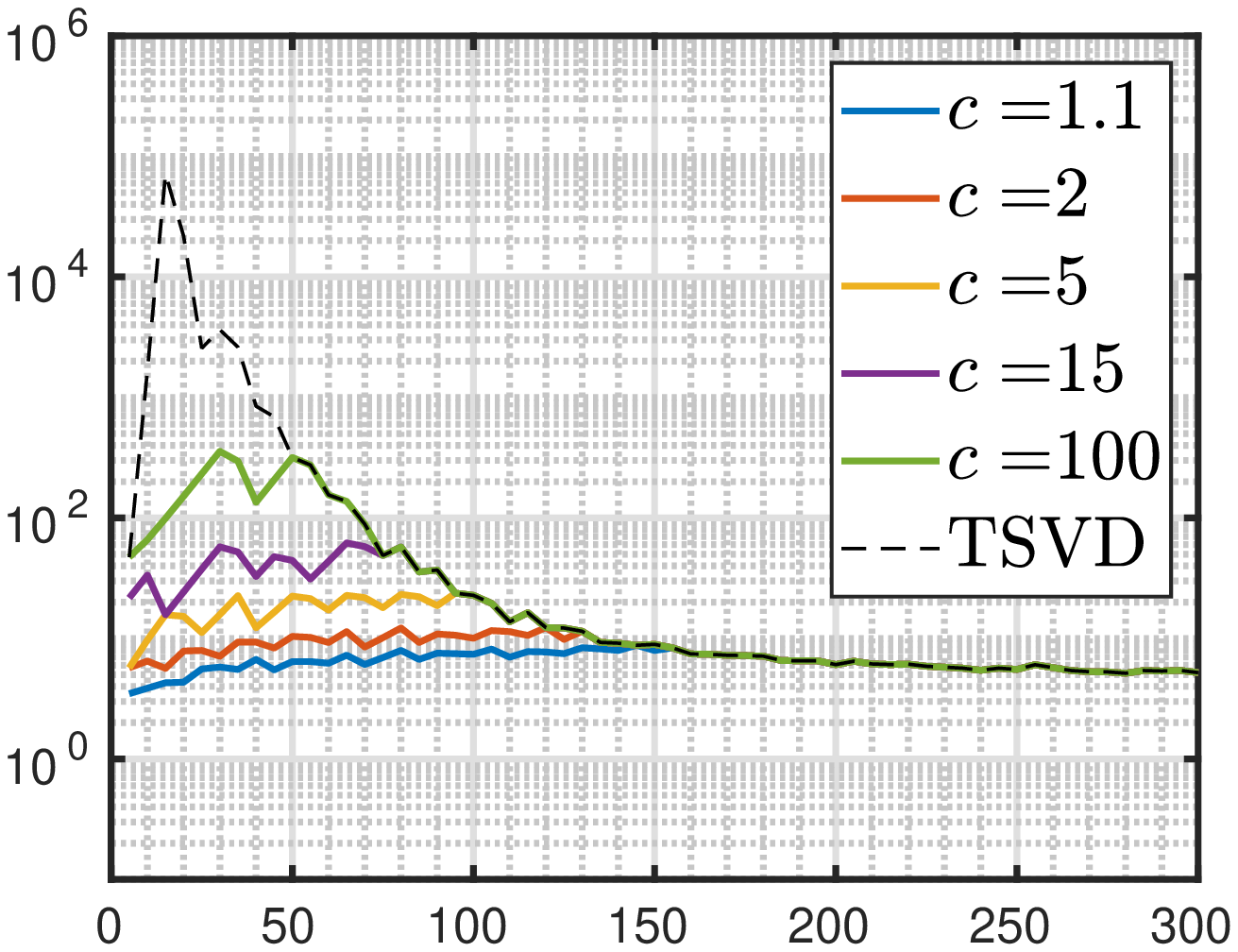}&
\includegraphics[width=\four]{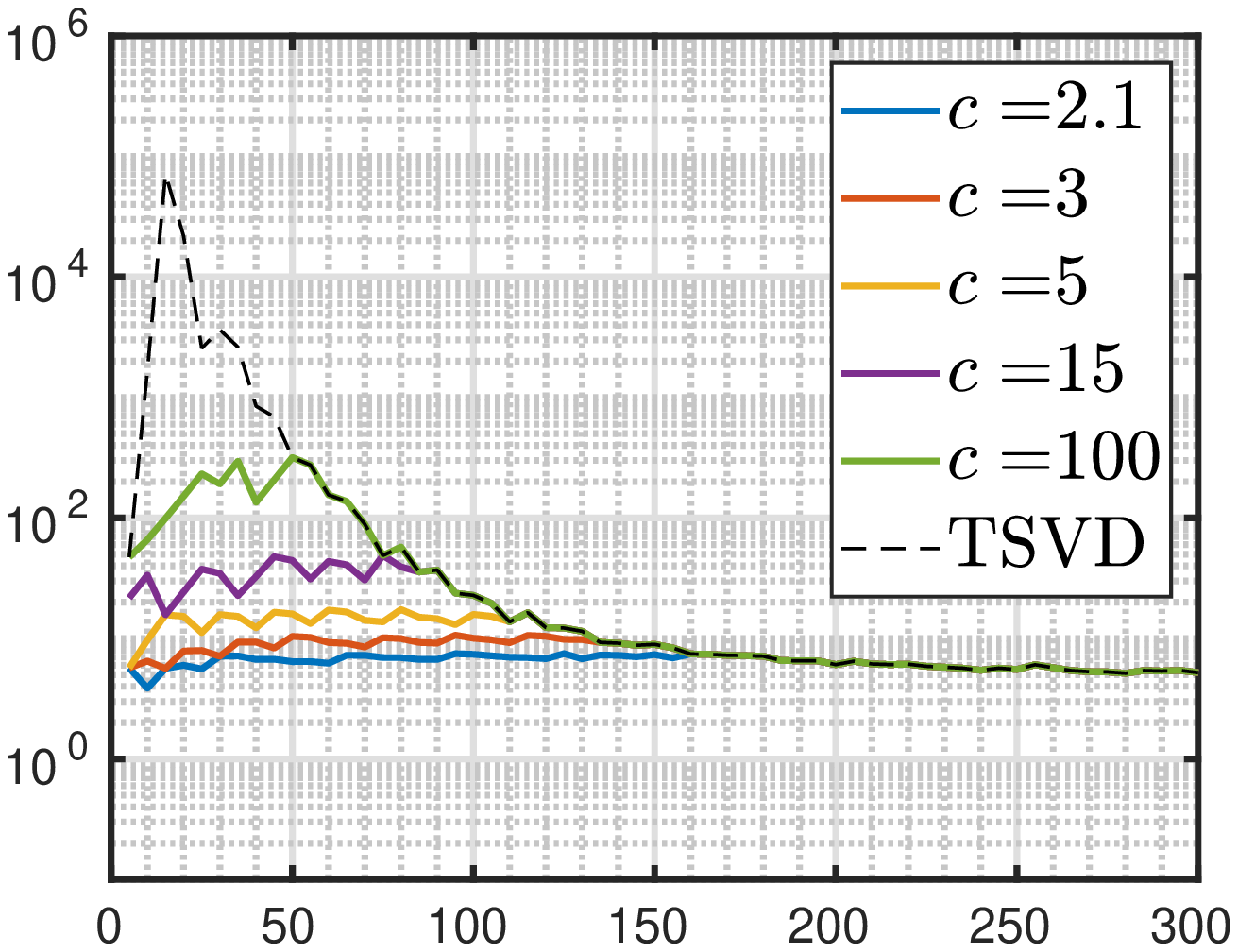}
\\
\includegraphics[width=\four]{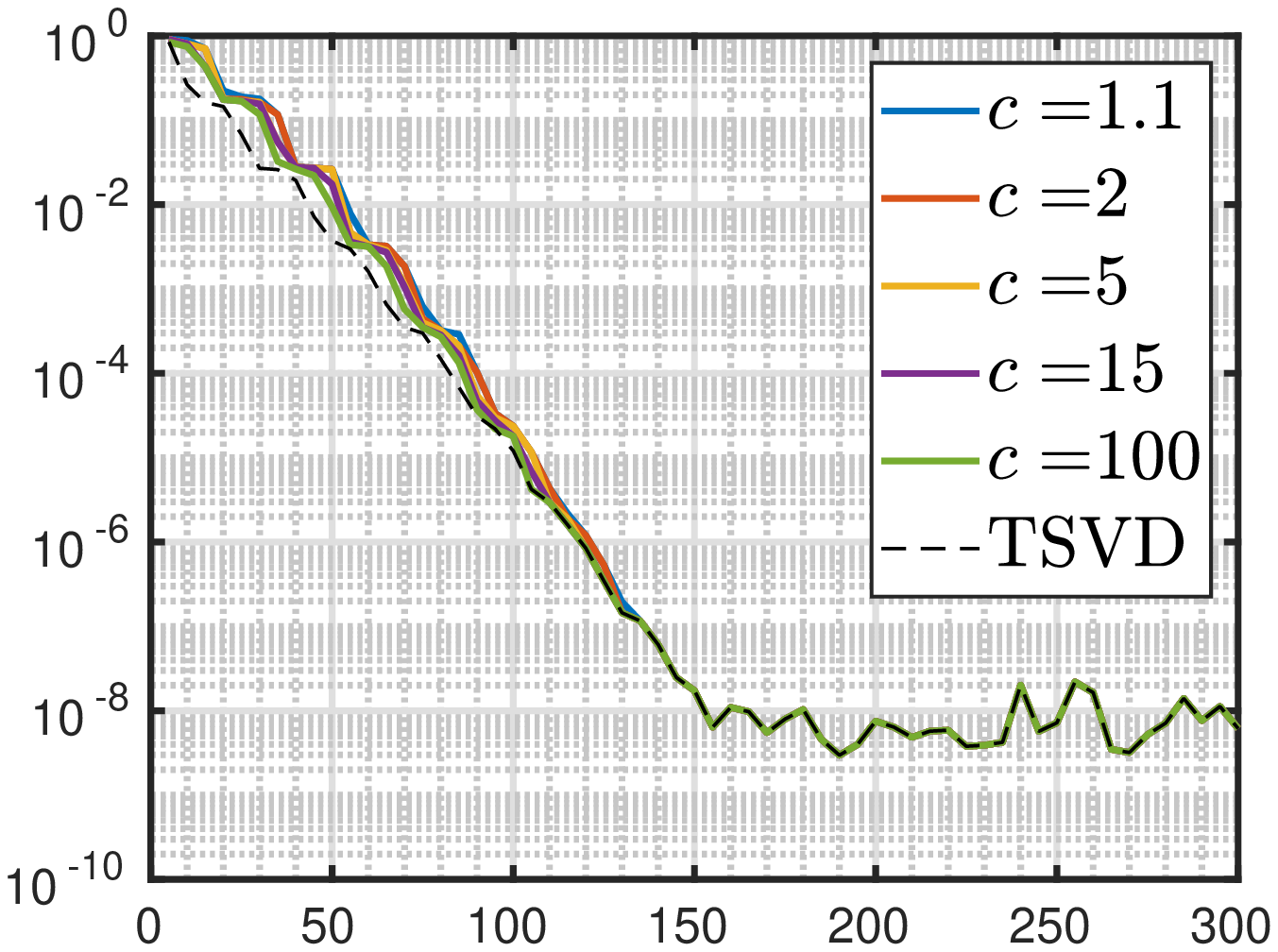} &
\includegraphics[width=\four]{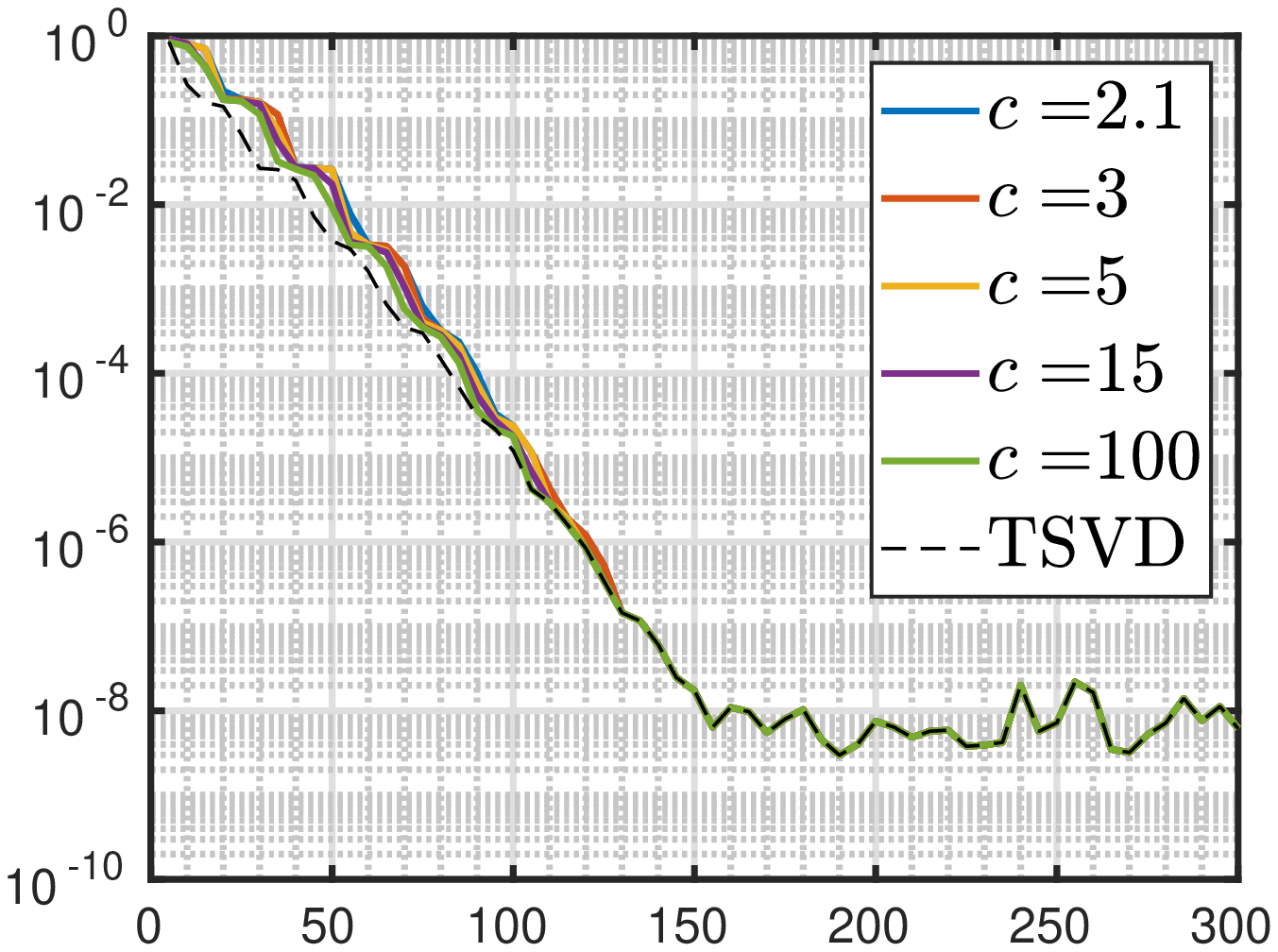} &
\includegraphics[width=\four]{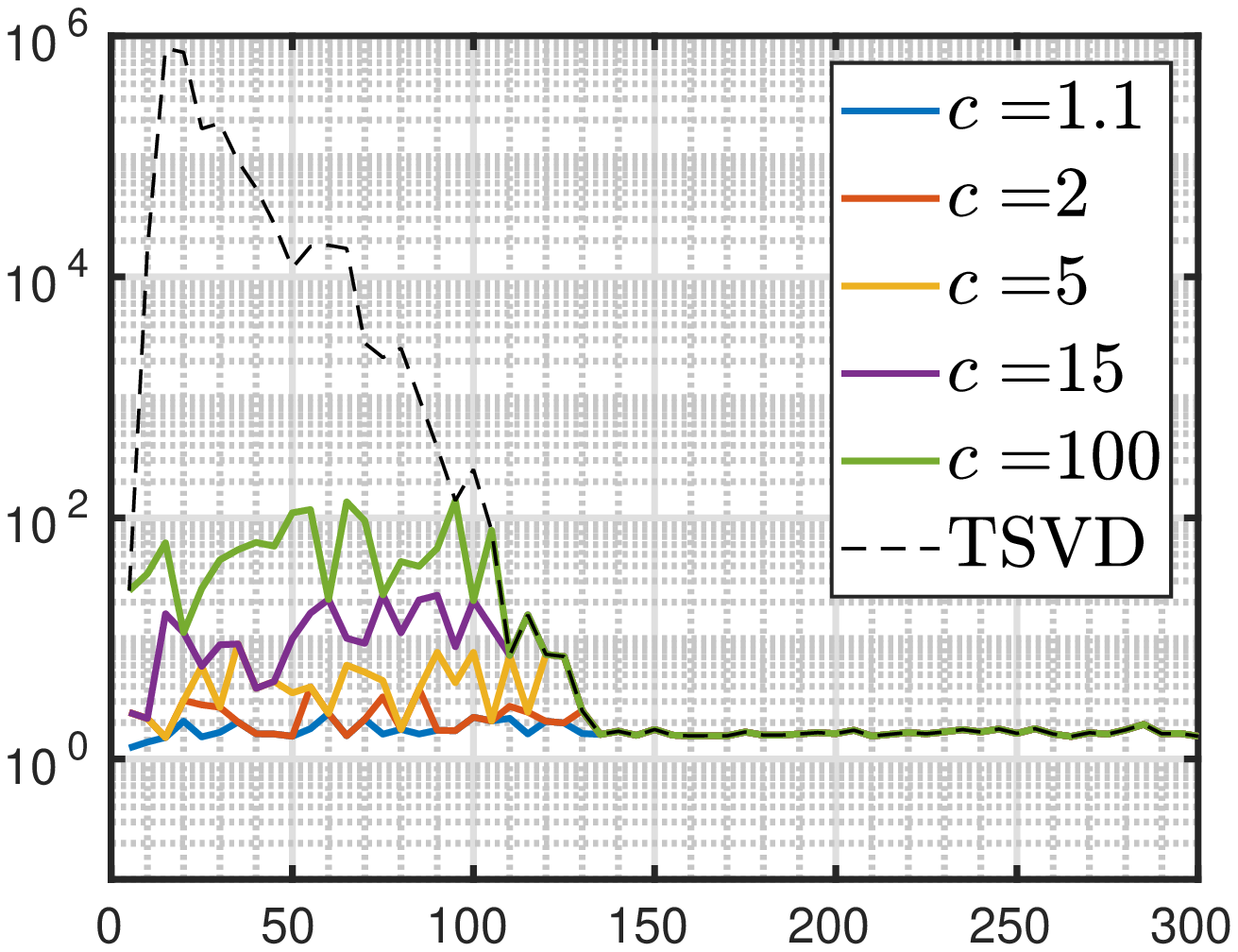} &
\includegraphics[width=\four]{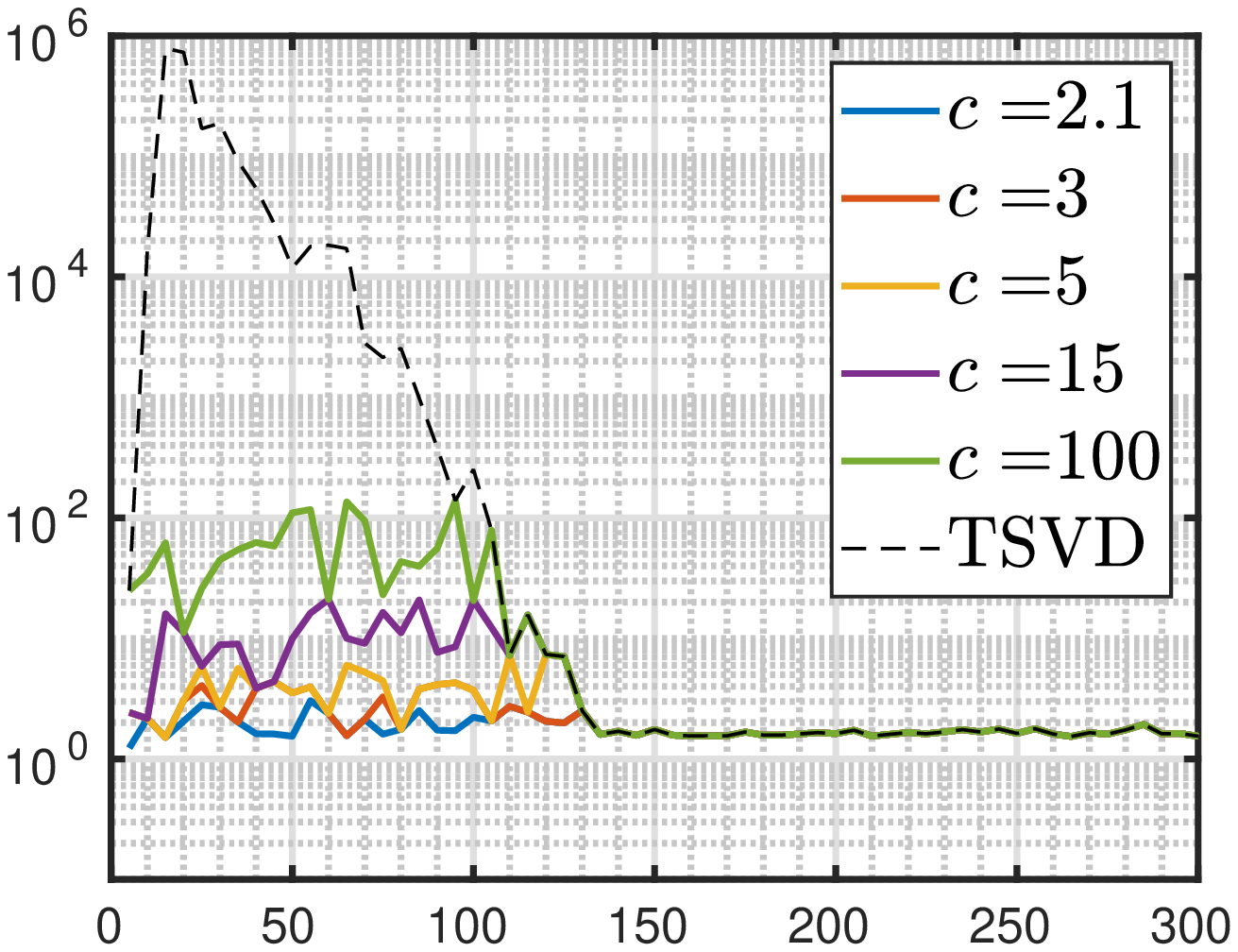}
\end{tabular}
}
\end{center}
\vspace*{-5mm}\caption{
Comparison of ASVD1 and ASVD2 for various $c$. From left to right: $\rL^2$-norm error versus $N$ for ASVD1, $\rL^2$-norm error versus $N$ for ASVD2, norm of the coefficient vectors for ASVD1, norm of the coefficient vectors for ASVD2. We also show the TSVD approximation for comparison. The functions used are $f_{1}(t)=\frac{1}{1+75t{^2}}$ (top), $f_{2}(t)=\frac{1}{0.57-t}$ (middle) and $f_{3}(t)=e^{\sin{\left(20t+0.5\right)}}\sqrt{1+t}\cos{\left(10t\right)}$ (bottom). The threshold $\epsilon$ is set to be $10^{-15}$ in all cases.
}
\label{fig:IGFNA1ASVD}
\end{figure}

In Fig.\ \ref{fig:IGFNA1ASVD} we compare ASVD1 and ASVD2 with TSVD for this frame, using various choices of $c$ and several different functions. As is evident, both ASVD methods successfully maintain small-norm coefficients, while yielding similar approximation errors to the TSVD approximation. Typically, the approximation error is slightly worse when $c$ is smaller, with the decline in performance depending on the function. The first fact is unsurprising: smaller $c$ means that coefficient vectors $\bm{z}$ with large norms but for which $\nmu{f - \cT_N \bm{z}}$ is small are not attainable by the ASVD methods, whereas they may be obtained by TSVD. The dependence on the function is also unsurprising in light of \R{LegFrameSob}. The function $f_2$ has large derivatives, and therefore the coefficient vectors of \R{LegFrameSob} have norms growing rapidly with the Sobolev order $k$. When $c$ is small, they may not be realizable by the ASVD approximations. Conversely, the derivatives of $f_1$ grow more slowly, meaning such coefficient vectors are potentially realizable.

\begin{figure}[t]
\begin{center}
{\small
\begin{tabular}{ccc}
\includegraphics[width=\two]{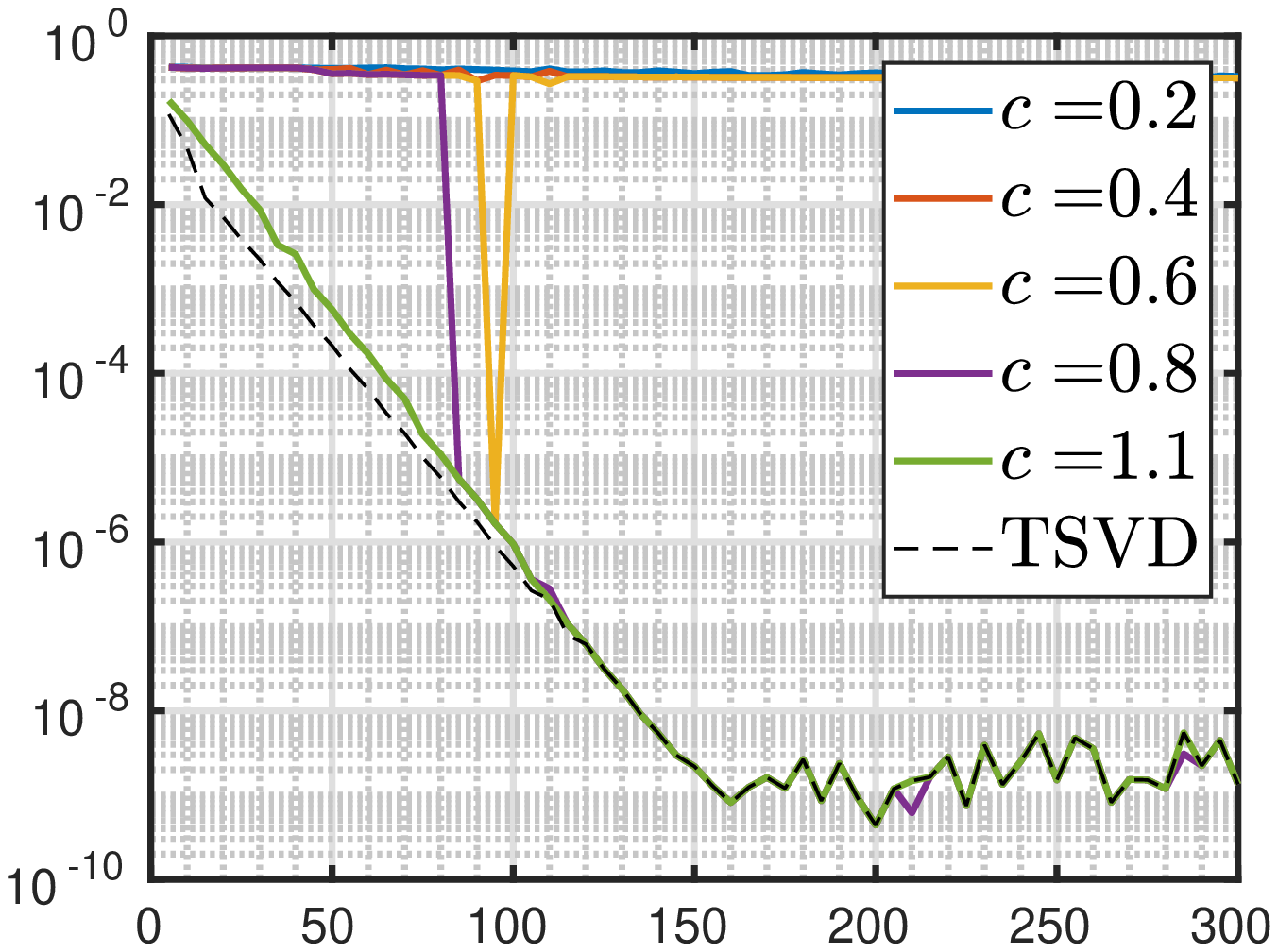}&
\includegraphics[width=\two]{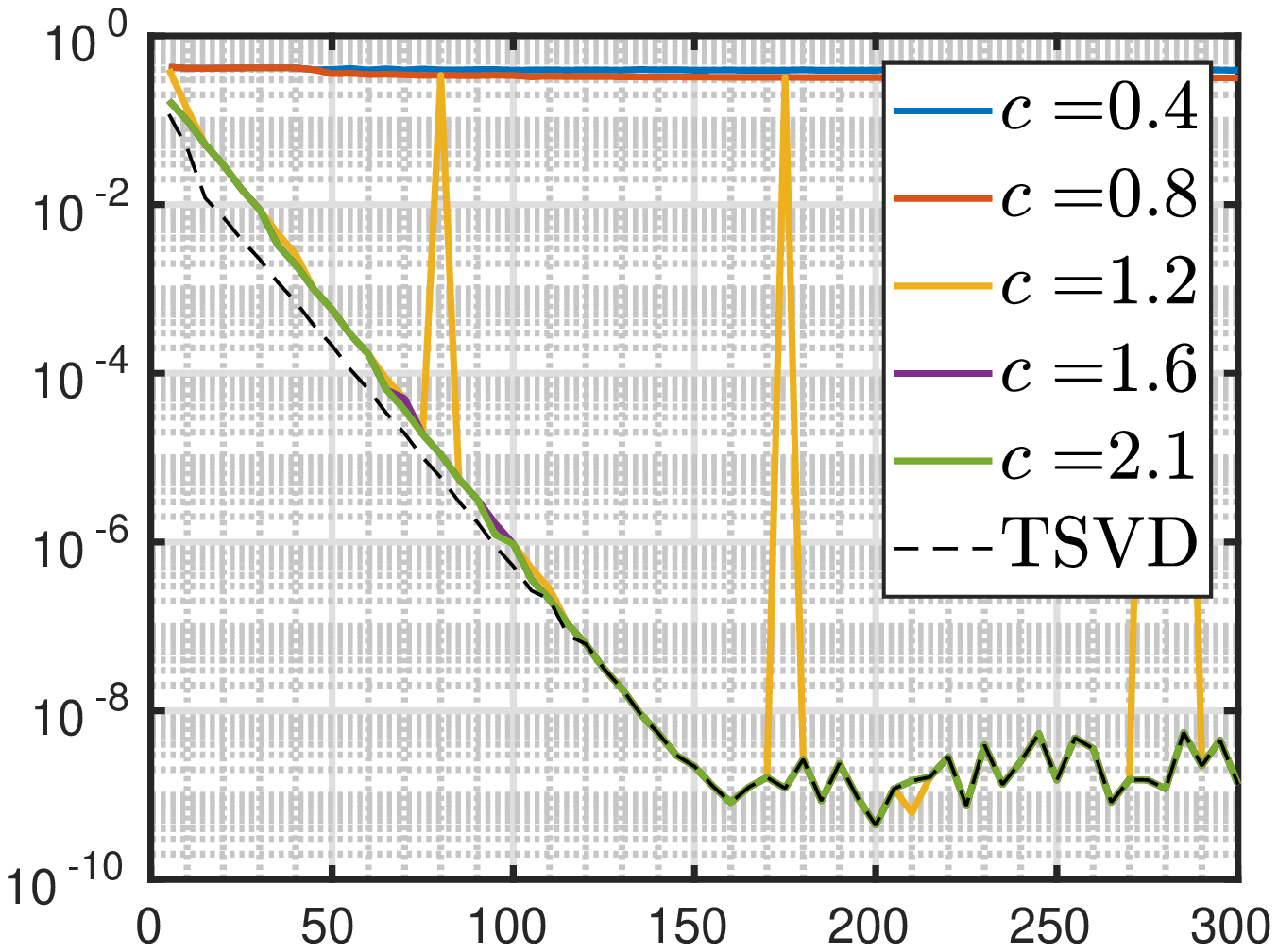}
\end{tabular}
}
\end{center}
\vspace*{-5mm}\caption{
Comparison of ASVD1 (left) and ASVD2 (right) for various $c$. The plots show the ${\rL}^2$-norm errors versus $N$. The function used was $f(t)=\frac{1}{1+75t{^2}}$. The threshold $\epsilon$ is set to be $10^{-15}$ in all cases.
}
\label{fig:IGFNA1ASVDC}
\end{figure}

In Fig.\ \ref{fig:IGFNA1ASVDC} we examine the influence of $c$ further. Notice that when $c$ is too small -- specifically, less than one for ASVD1 and less than two for ASVD2 -- the error behaves badly, either converging very slowly, or jumping from small to large. In \S \ref{section : theory}, we will support these conclusions with theory: specifically, we show that $c > 1/A$ and $c>2/A$ respectively are sufficient for convergence, where $A$ is the lower frame bound (recall that $A = 1$ for the Legendre polynomial frame). See Theorem \ref{thm:limsup}.

\begin{figure}[t]
\begin{center}
{\small
\begin{tabular}{ccc}
\includegraphics[width=\two]{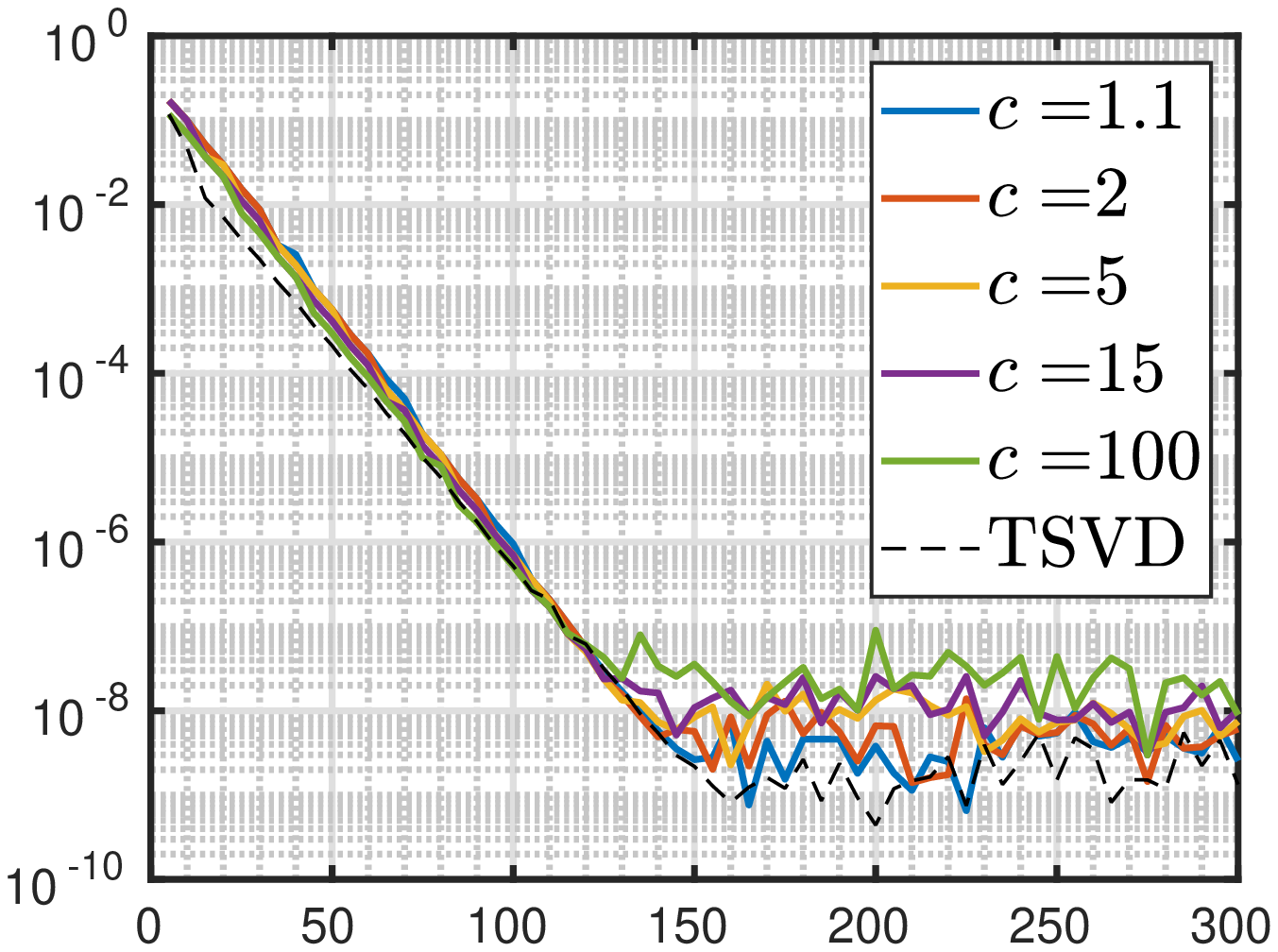}&
\includegraphics[width=\two]{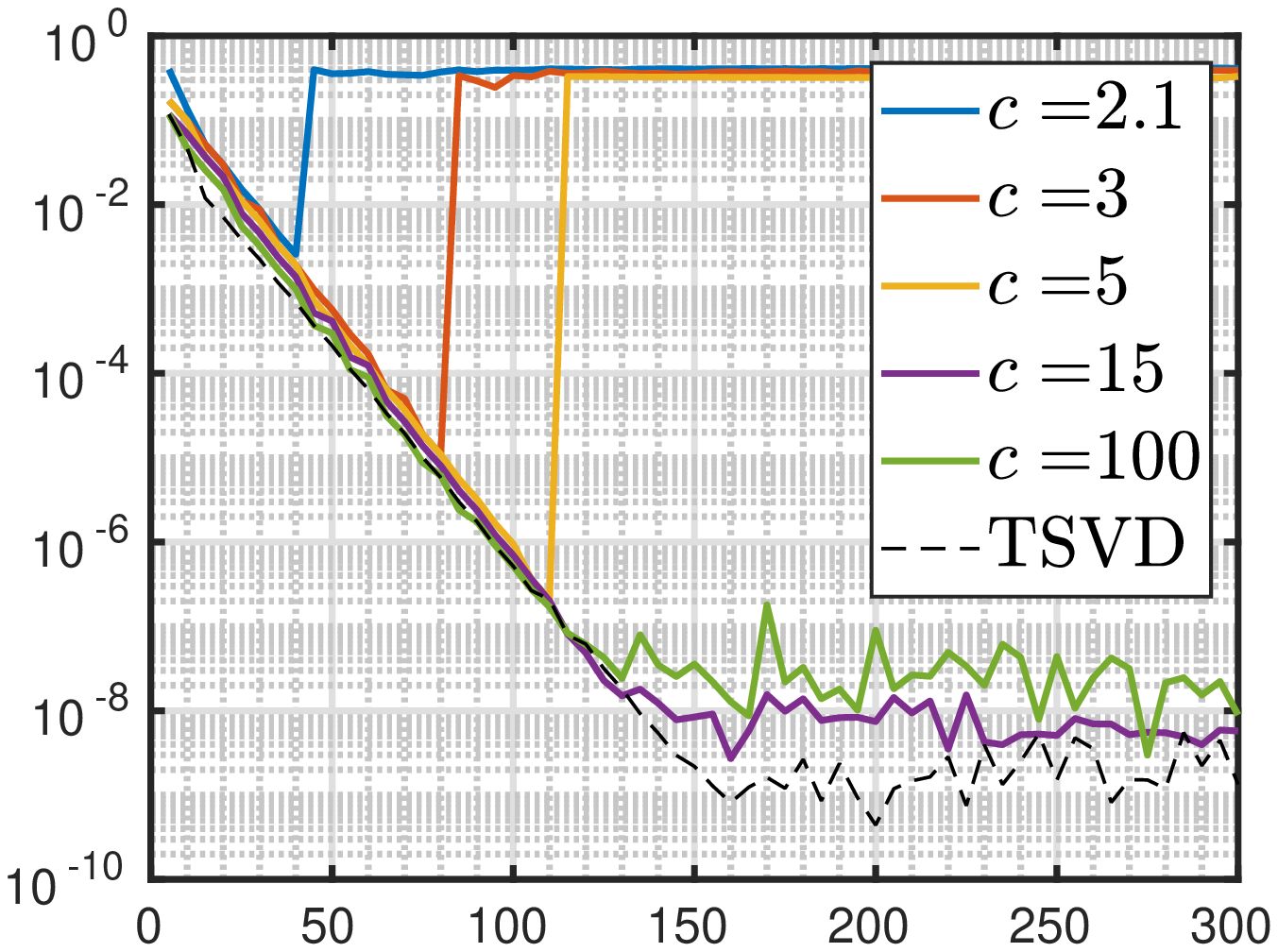}
\\  
\includegraphics[width=\two]{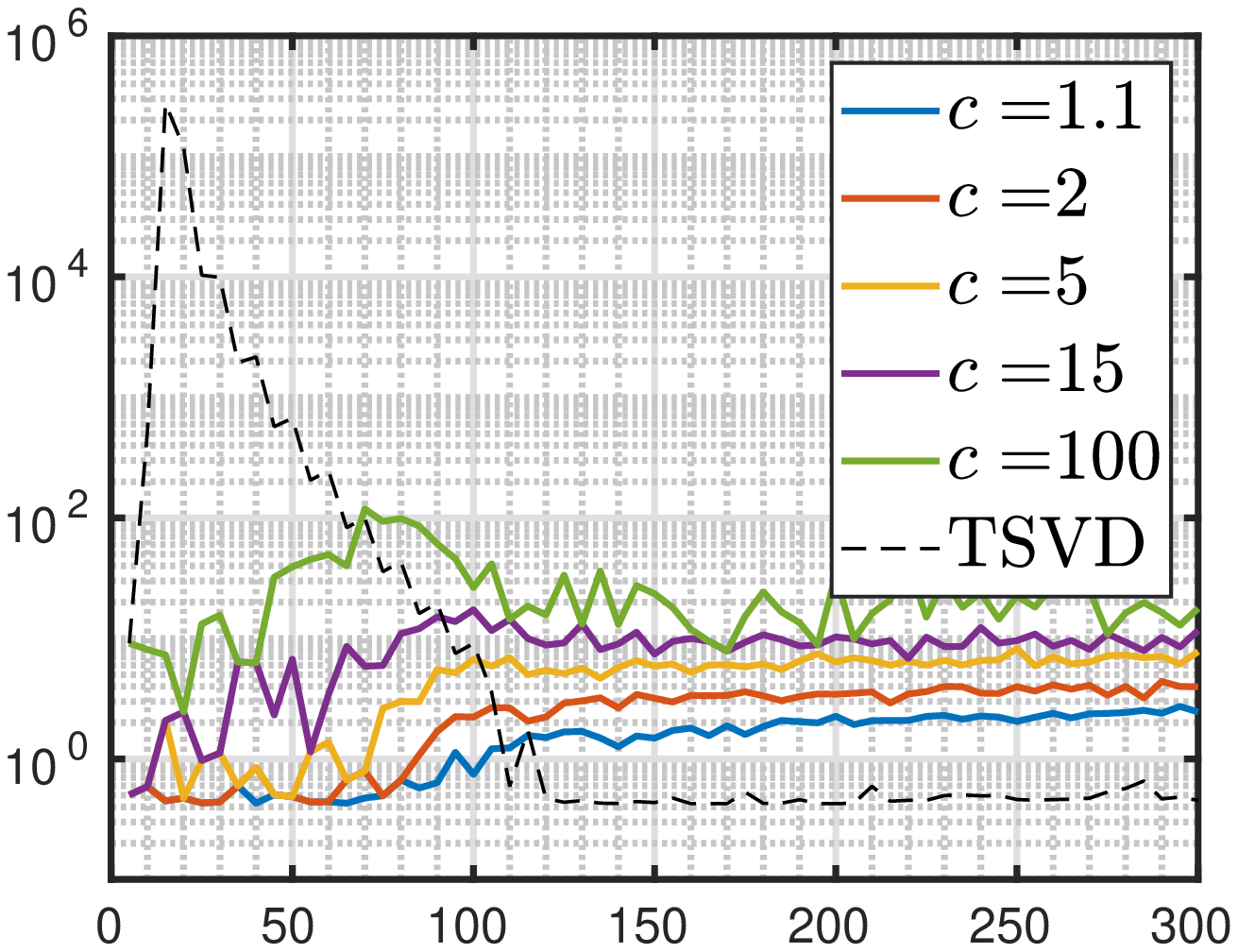}&
\includegraphics[width=\two]{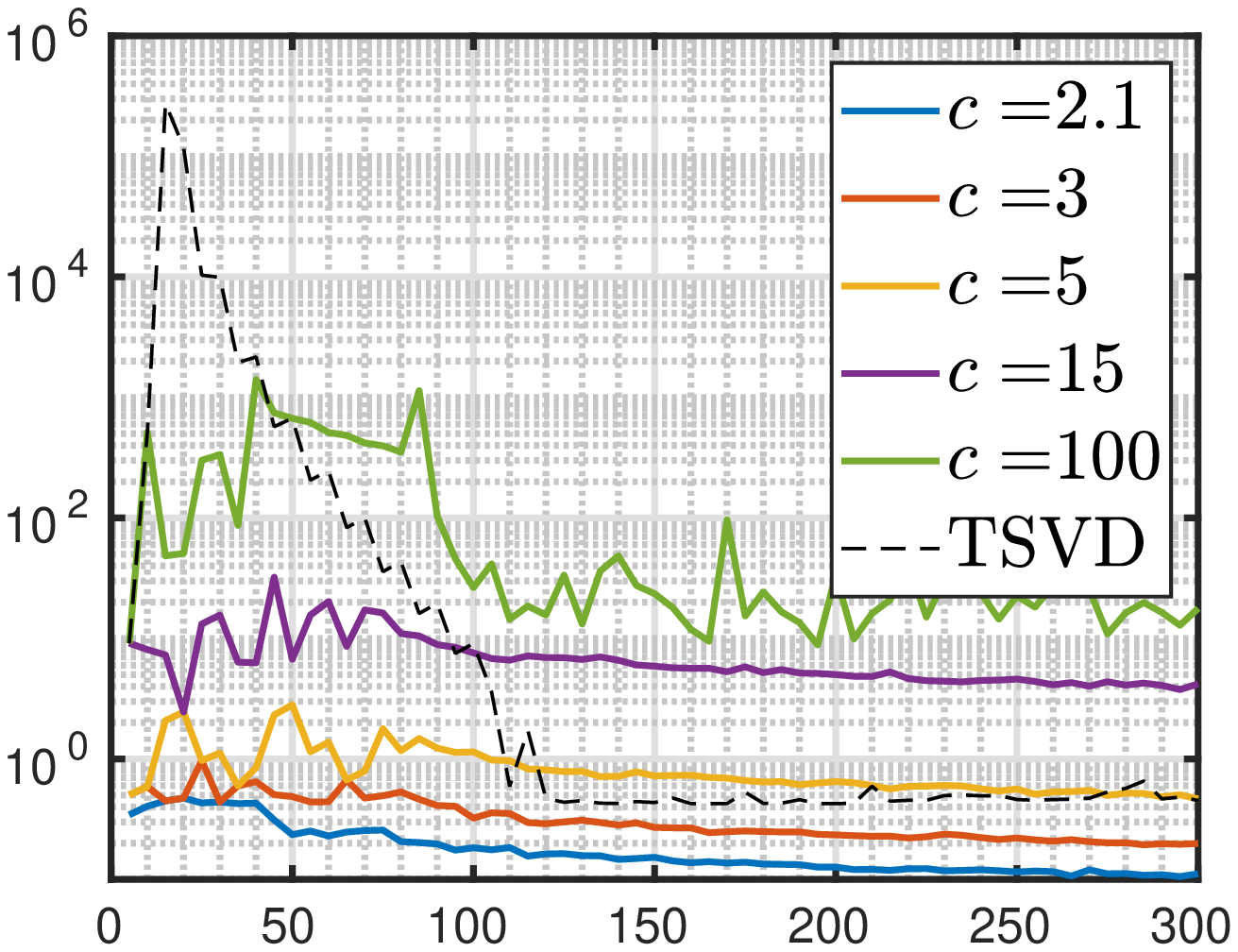}
\end{tabular}
}
\end{center}
\vspace*{-5mm}\caption{
Comparison of ASVD1 (left) and ASVD2 (right) for various $c$. Top row: ${\rL}^2$-norm errors versus $N$. Bottom row: the norms of the coefficient vectors versus $N$. The function used was $f(t)=\frac{1}{1+75t{^2}}$. The threshold $\epsilon = 0$ is used for ASVD1 and ASVD2, and for TSVD the value $\epsilon = 10^{-15}$ is used.
}
\label{fig:IGFNA1ASVDFail}
\end{figure}

Finally, in Fig.\ \ref{fig:IGFNA1ASVDFail} we examine why a truncation parameter $\epsilon > 0$ is generally a good idea in ASVD1 and ASVD2. This figure shows the behaviour of each method when no truncation occurs, i.e.\ $\epsilon = 0$. Notice that for ASVD1 the method still produces an acceptable error. Yet the coefficients are no longer $\ord{1}$ in the limit $N \rightarrow \infty$, unlike those of the TSVD approximation. Indeed, the limiting accuracy and coefficient norm both appear to scale with $c$. This is undesirable; we would certainly prefer the ASVD methods to retain the desirable limiting behaviour of the TSVD approximation. For ASVD2 the behaviour is even worse. For smaller $c$, the error jumps to $\ord{1}$, indicating lack of (or at best very slow) convergence of the approximation. Incorporating a threshold appears vital for rapid convergence in this case. Our theoretical results in \S \ref{section : theory} support this conclusion.

\subsection{Approximating weakly singular functions}\label{ss:numexamp2}

We now consider approximating functions of the form
\be{ \label{eq: Functions with singularities}
f(t) = w(t)g(t) + h(t), \quad t \in [0, 1],
}
where $g$ and $h$ are smooth functions, and $w \in \rL^{2}(0, 1)$ is a known function that may be singular.  We shall consider the case $w(t) = \log(t)$, i.e.\ $f(t)$ has a logarithmic singularity at $t = 0$.
This problem was considered in \cite{BADHFramesPart2}, and the approach used therein is employs the frame
\be{
\label{frameweuse}
\Phi = \{ \varphi_n \}^{\infty}_{n=0} \cup \{ \psi_j \}^{K-1}_{j=0},
}
where $\varphi_n$ is the orthonormal Legendre polynomial basis on $(0,1)$, $\psi_j(t) = \log(t) \varphi_{j}(t)$ and $K \in \bbN$. {This forms a frame for any finite $K$. It is an example of a general construction of a frame, in which an orthonormal basis of a Hilbert space $\rH$ is augmented with a finite number of arbitrary elements from $\rH$.}
The rationale for using \R{frameweuse} for functions of the form \R{eq: Functions with singularities} that the $\psi_j$ capture the logarithmic singularity up to some finite order depending on $K$ and the polynomials capture the remaining smooth part.  To see this, let 
\bes{
\Phi_N =\{ \psi_j \}^{K-1}_{j=0} \cup \{ \varphi_n \}^{N-K-1}_{n=0},\qquad N > K,
}
Similar to the previous example, we can show that exists an approximation to $f$ from $\Phi_N$ with bounded coefficients that converges at an algebraic rate depending on $K$. Let $p = \sum^{K-1}_{k=0} z_k \varphi_k$ be the $K^{\rth}$ order Taylor polynomial of $g$ around $t = 0$, so that $p^{(j)}(0) = g^{(j)}(0)$ for $j = 0,\ldots,K-1$. Then let
\bes{
q = \cQ_{N-K} (h + w(g-p)) = \sum^{N-K-1}_{n=0} z_{K+n} \varphi_n,\qquad z_{K+n} = \ip{h + w ( g - p)}{\varphi_n},
}
be the orthogonal projection of $h + w(g-p)$ onto $\spn \{ \varphi_{n} \}^{N-K-1}_{n=0}$. Clearly,
\bes{
|z_j| \leq C_{K} \max_{0 \leq l < K} | g^{(l)}(0) | \leq C_K \nm{g}_{H^{K}}, \qquad j = 0,\ldots,K-1,
}
for some constant $C_{K} > 0$ depending on $K$. This implies that the coefficients  $\bm{z} = (z_n)^{N-1}_{n=0}$ are uniformly bounded in $N$. Indeed, we have
\bes{
\sqrt{\sum^{N-K-1}_{n=0} |z_{K+n}|^2} \leq \nm{h + w(g-p)} \leq \nm{h} + \nm{w(g-p)} \leq \nm{h} + C_{K} \nm{g}_{H^{K}} ,
}
where $C_{K}$ is a possibly different constant, and therefore
\bes{
\nm{\bm{z}} \leq  \nm{h} + C_{K} \nm{g}_{H^{K}}.
}
Moreover, the approximation $\cT_{N} \bm{z}$ converges algebraically fast to $f$. To see why, observe that, by construction
\bes{
\nm{ f - \cT_{N} \bm{z}} = \nm{w (g-p) + h - \cQ_{N-K}(w(g-p) + h)} .
}
Since $h$ is smooth by assumption and $g-p$ vanishes along with its first $K-1$ derivatives at $t = 0$, the function $w (g-p) + h$ has roughly $K$ orders of smoothness. Hence $\nm{ f - \cT_{N} \bm{z}}$ decays at an algebraic rate in $N$ depending on $K$. 

\begin{figure}
\begin{center}
{
\begin{tabular}{cccc} 
\includegraphics[width=\four]{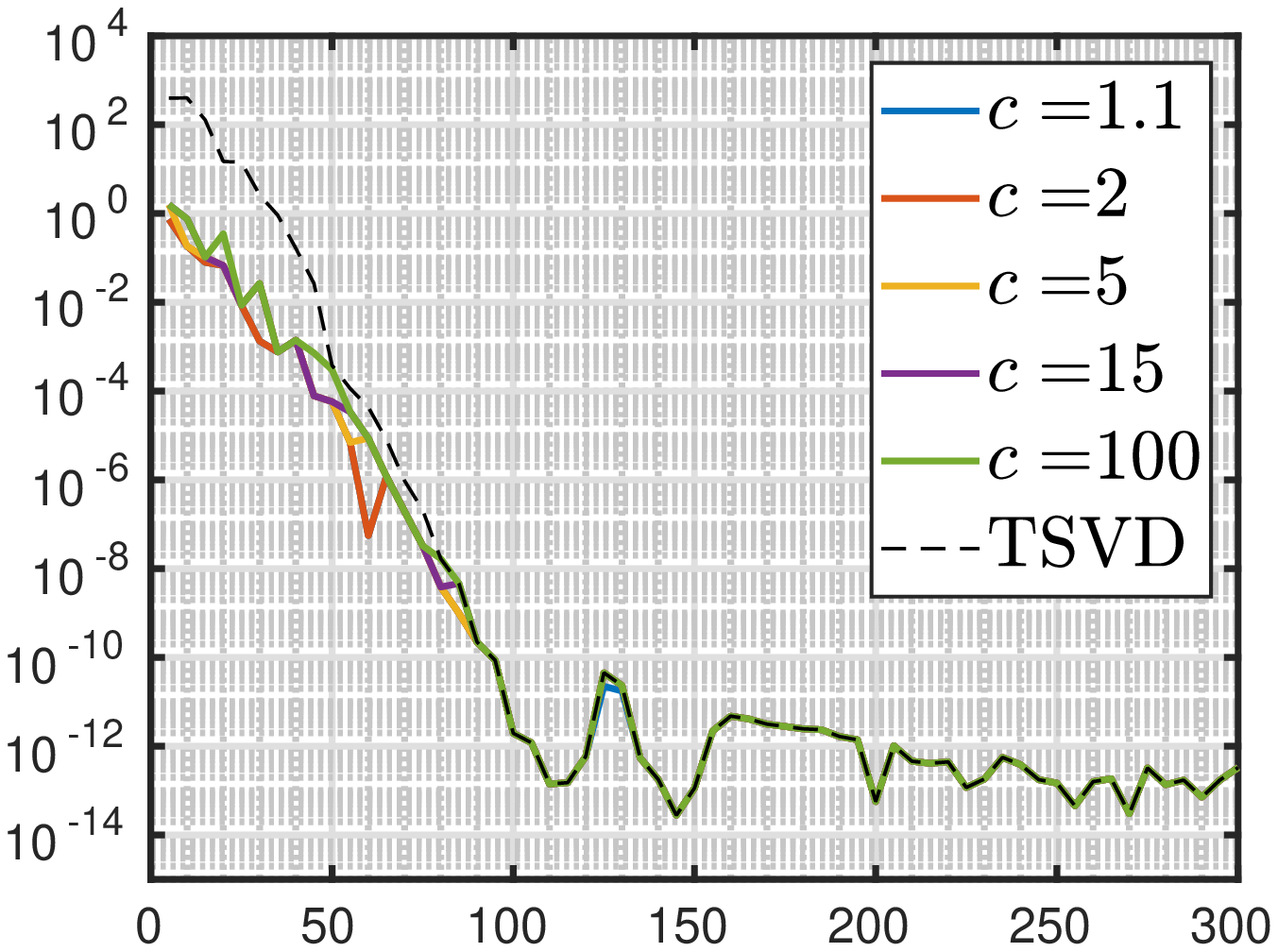}&
\includegraphics[width=\four]{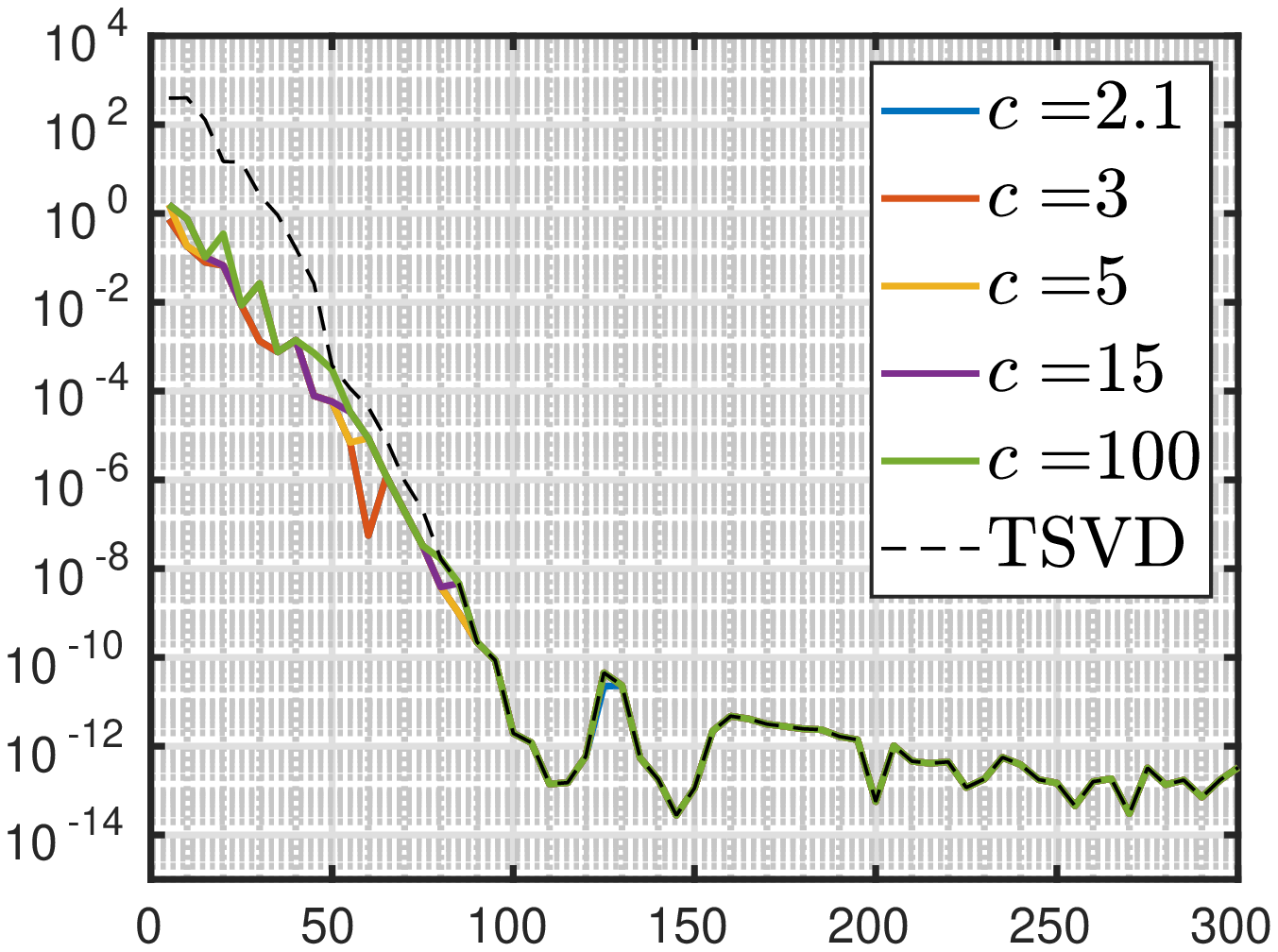}&
\includegraphics[width=\four]{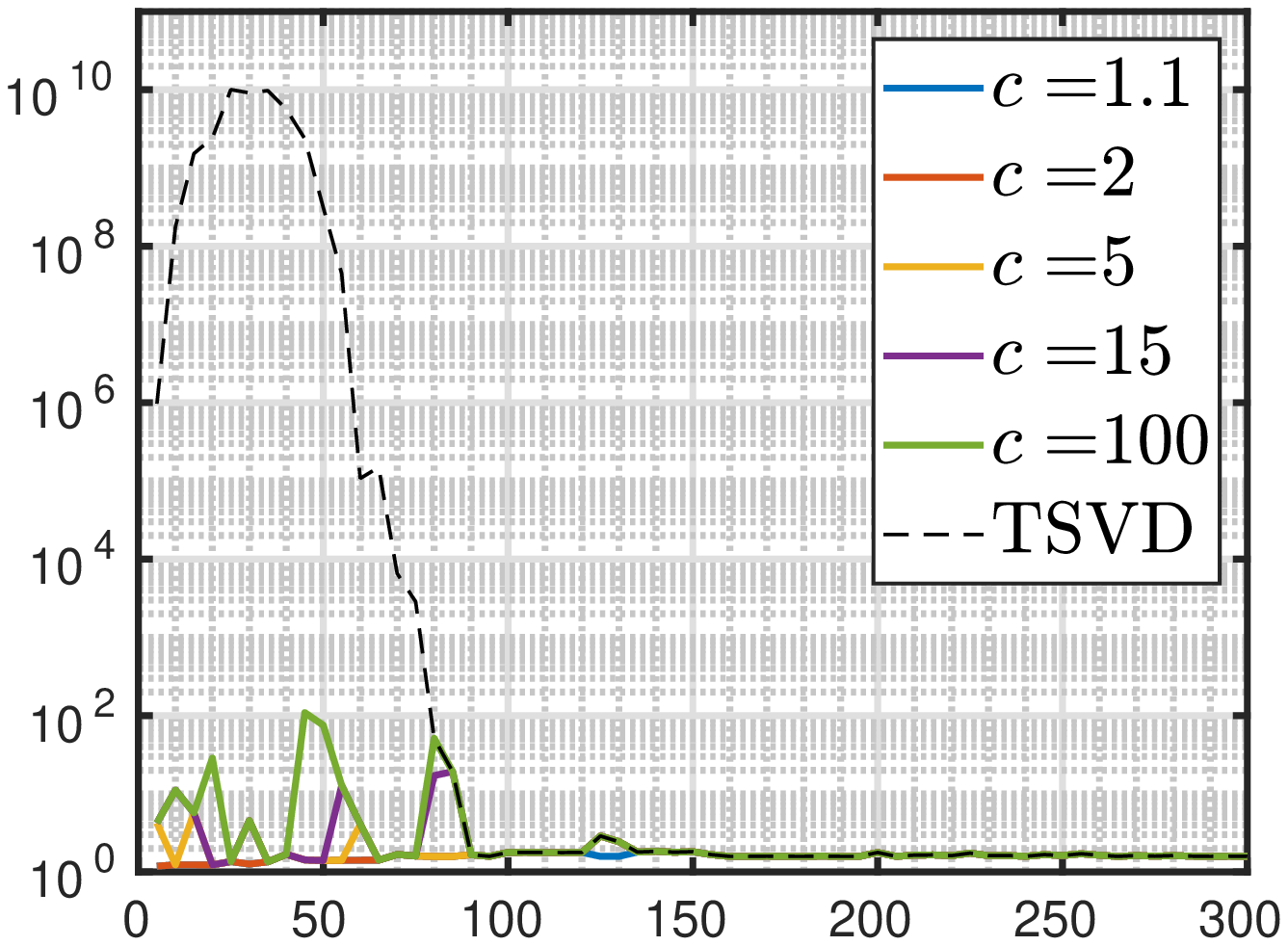}&
\includegraphics[width=\four]{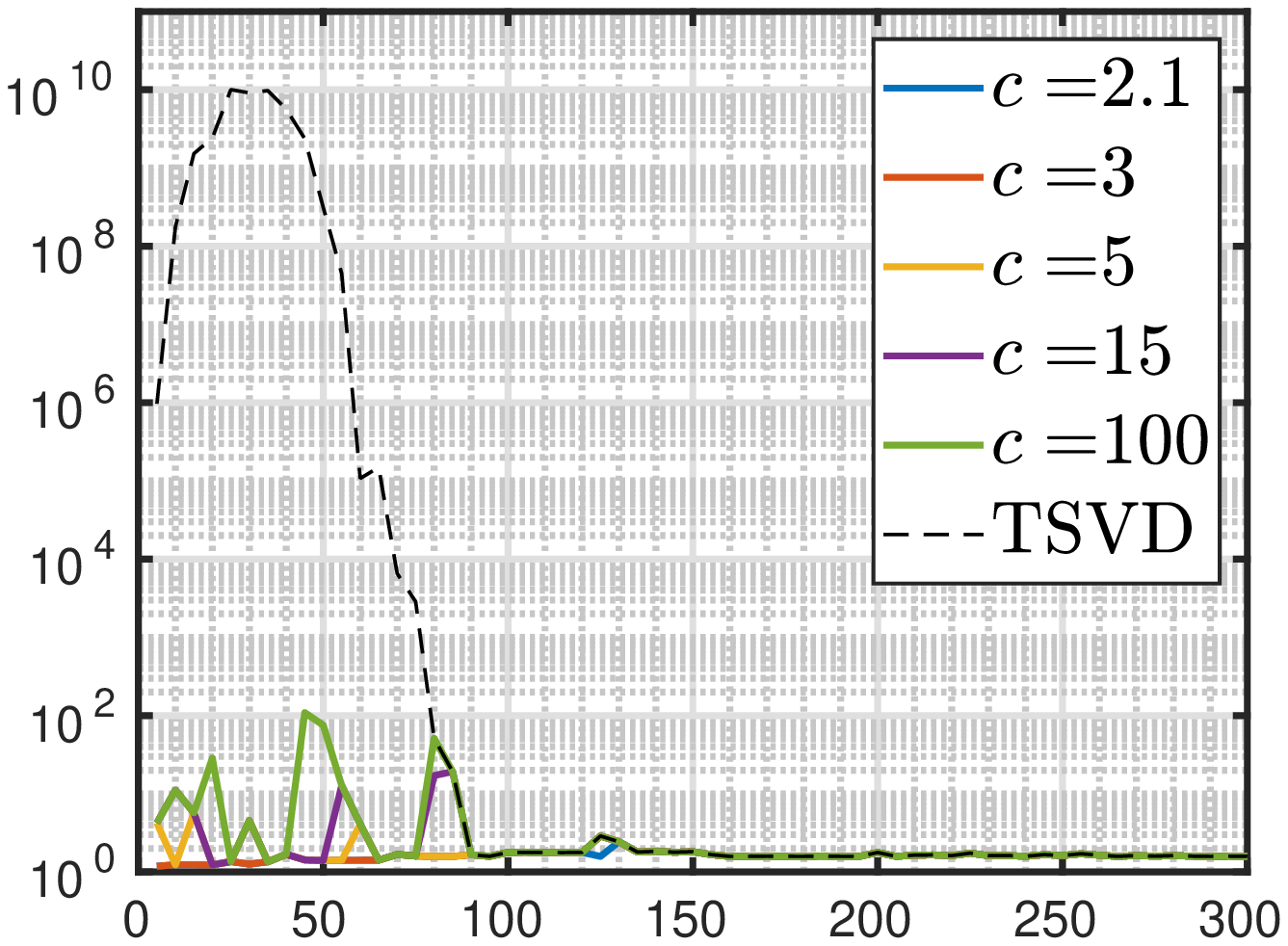}
\\
\includegraphics[width=\four]{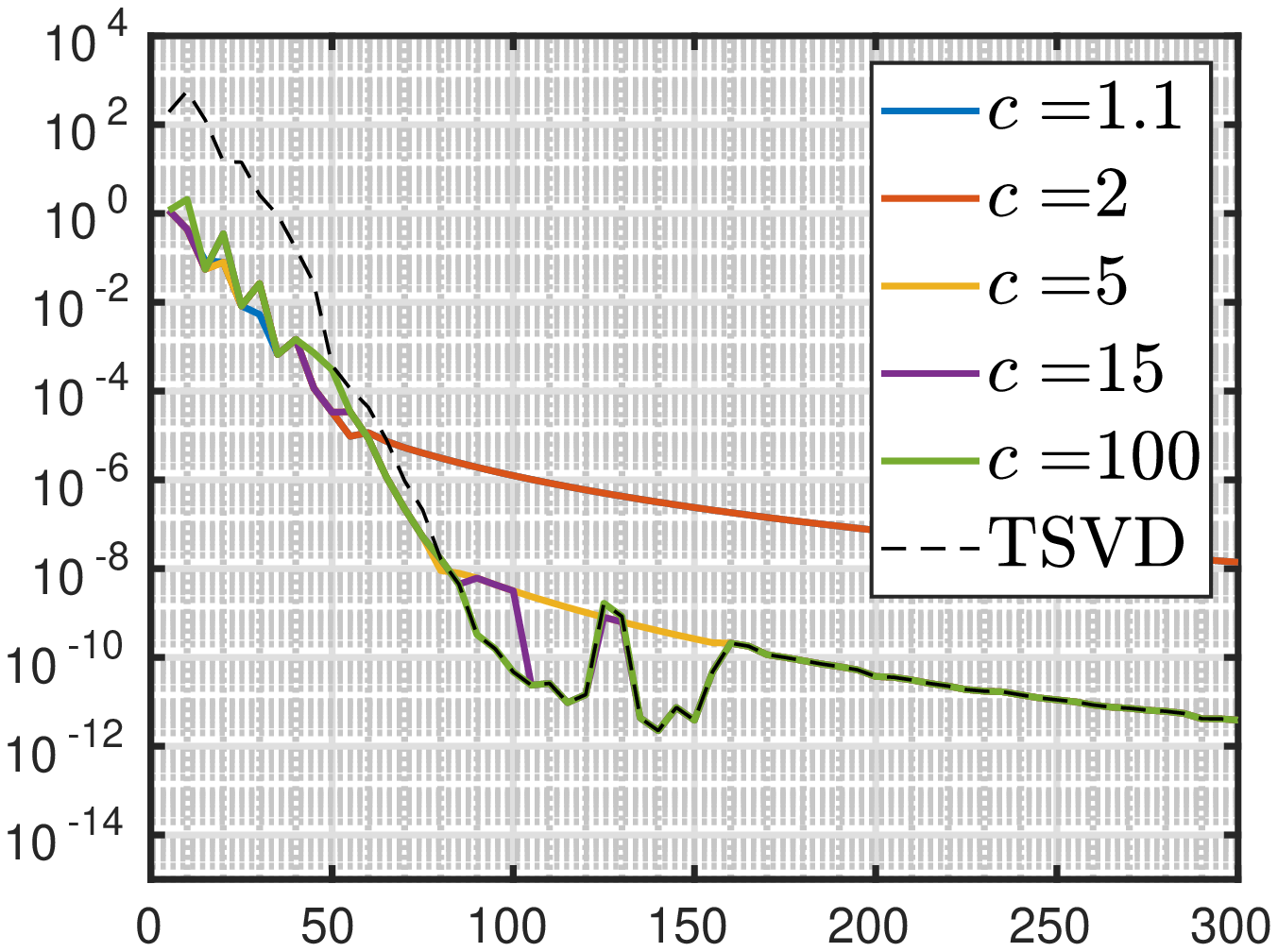}&
\includegraphics[width=\four]{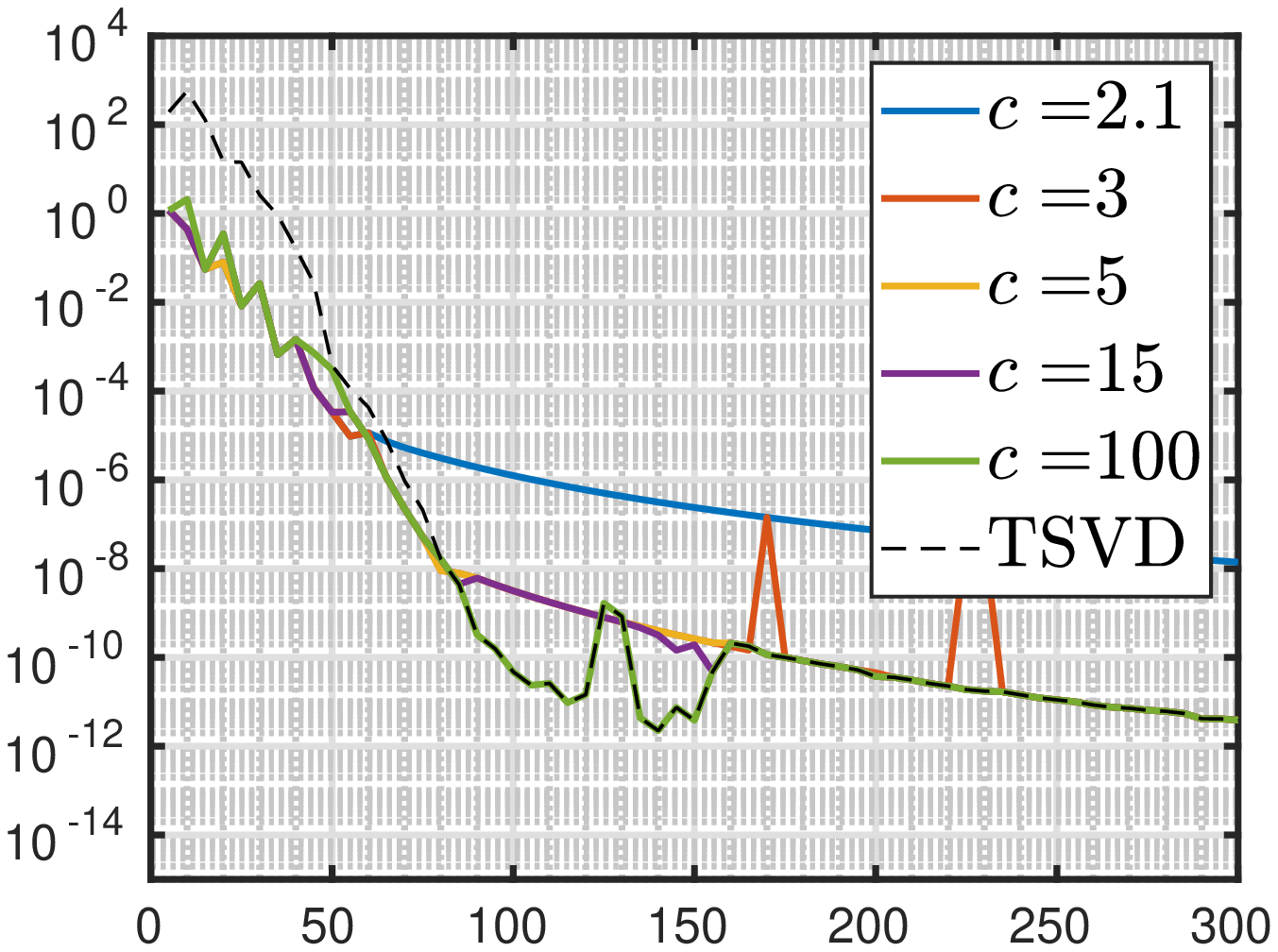}&
\includegraphics[width=\four]{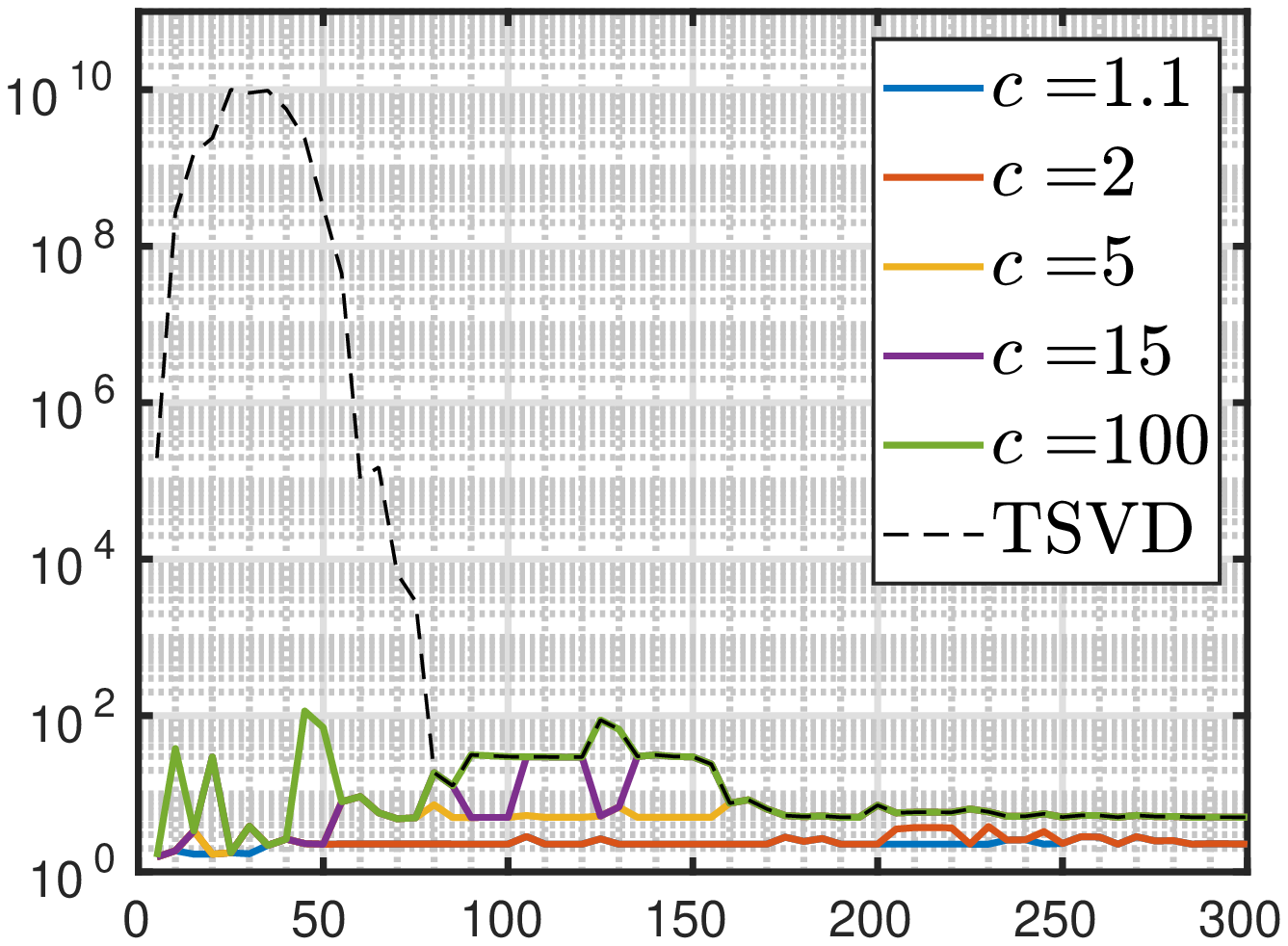}&
\includegraphics[width=\four]{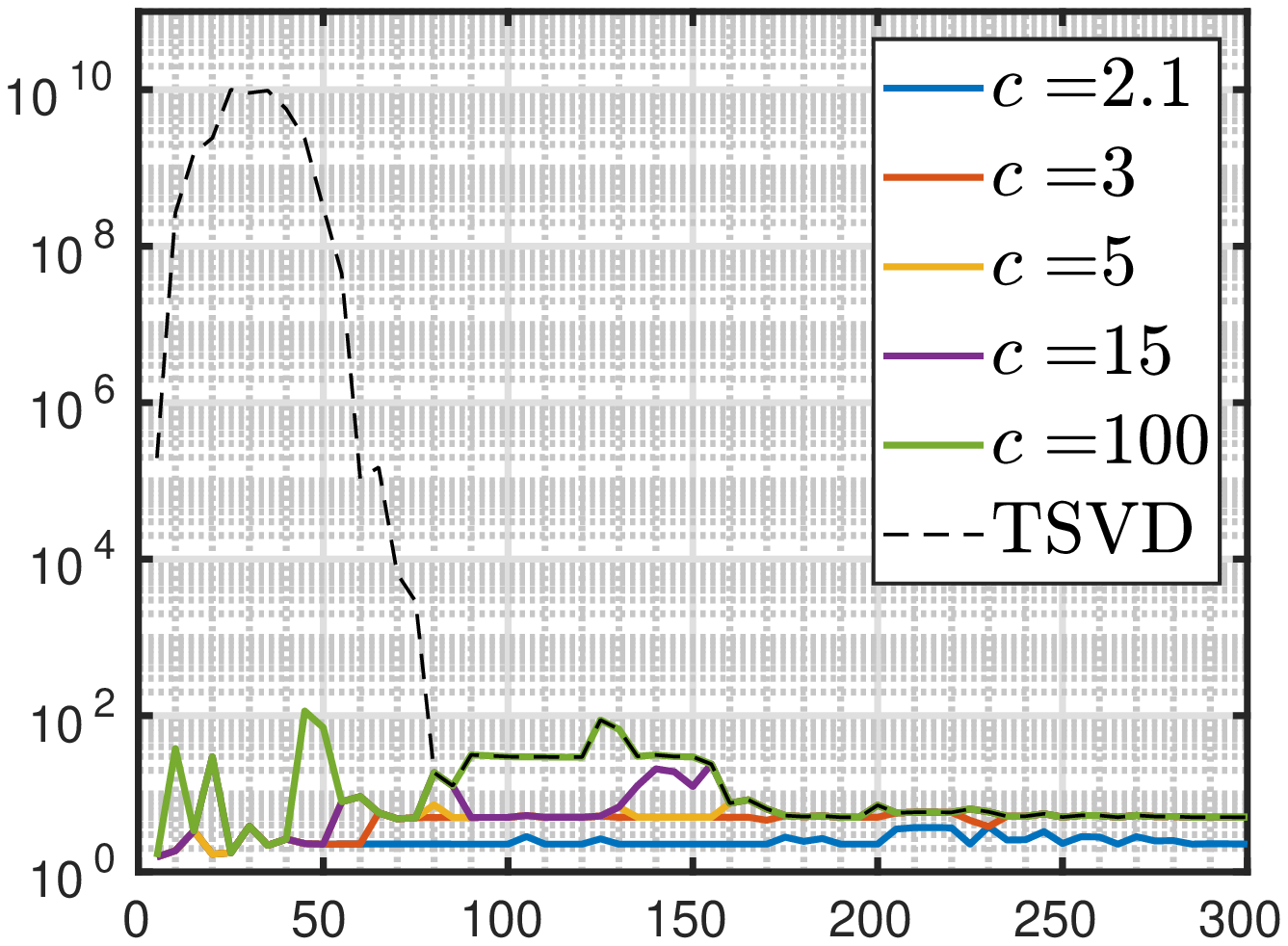}
\\
\includegraphics[width=\four]{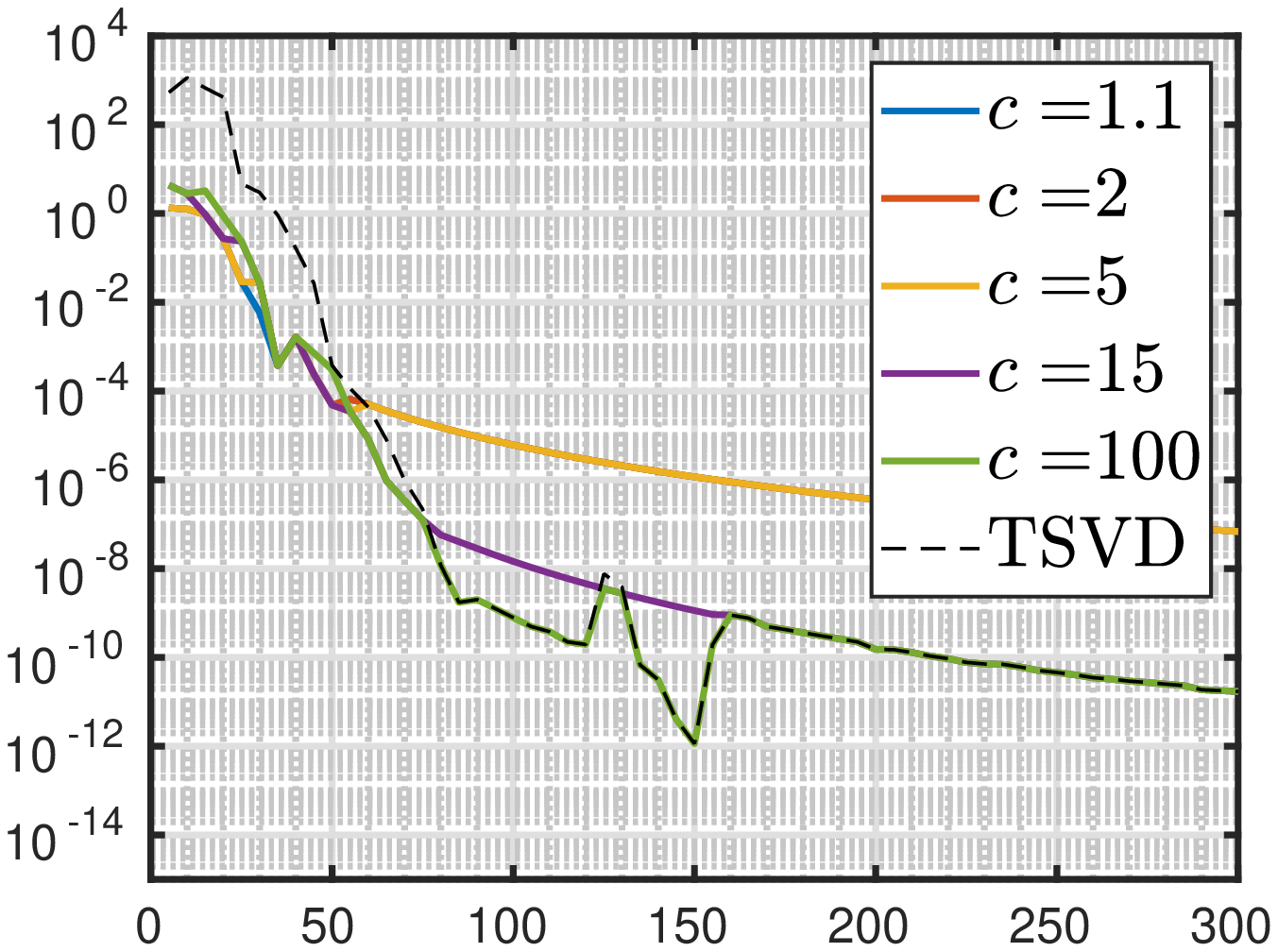} &
\includegraphics[width=\four]{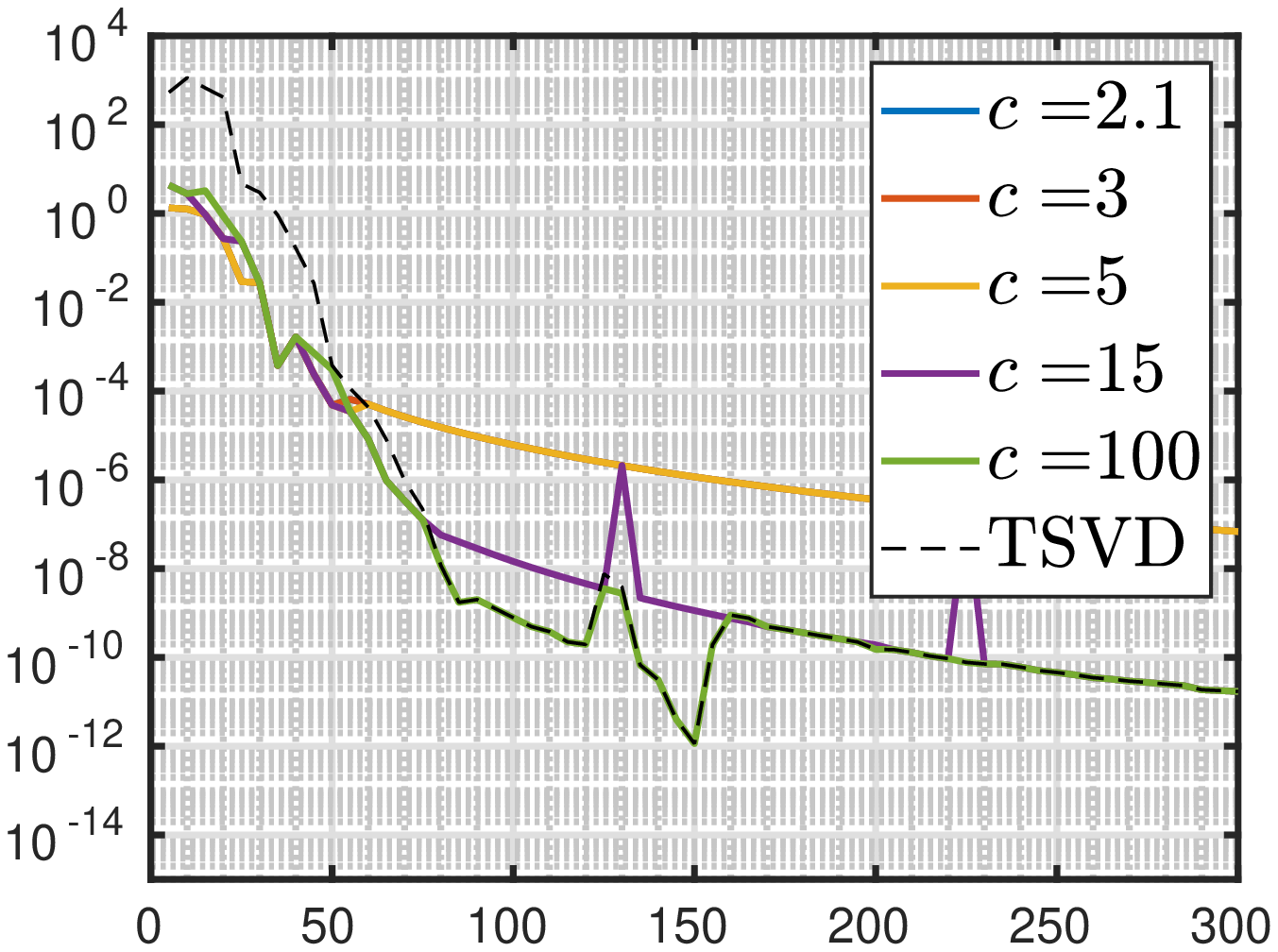} &
\includegraphics[width=\four]{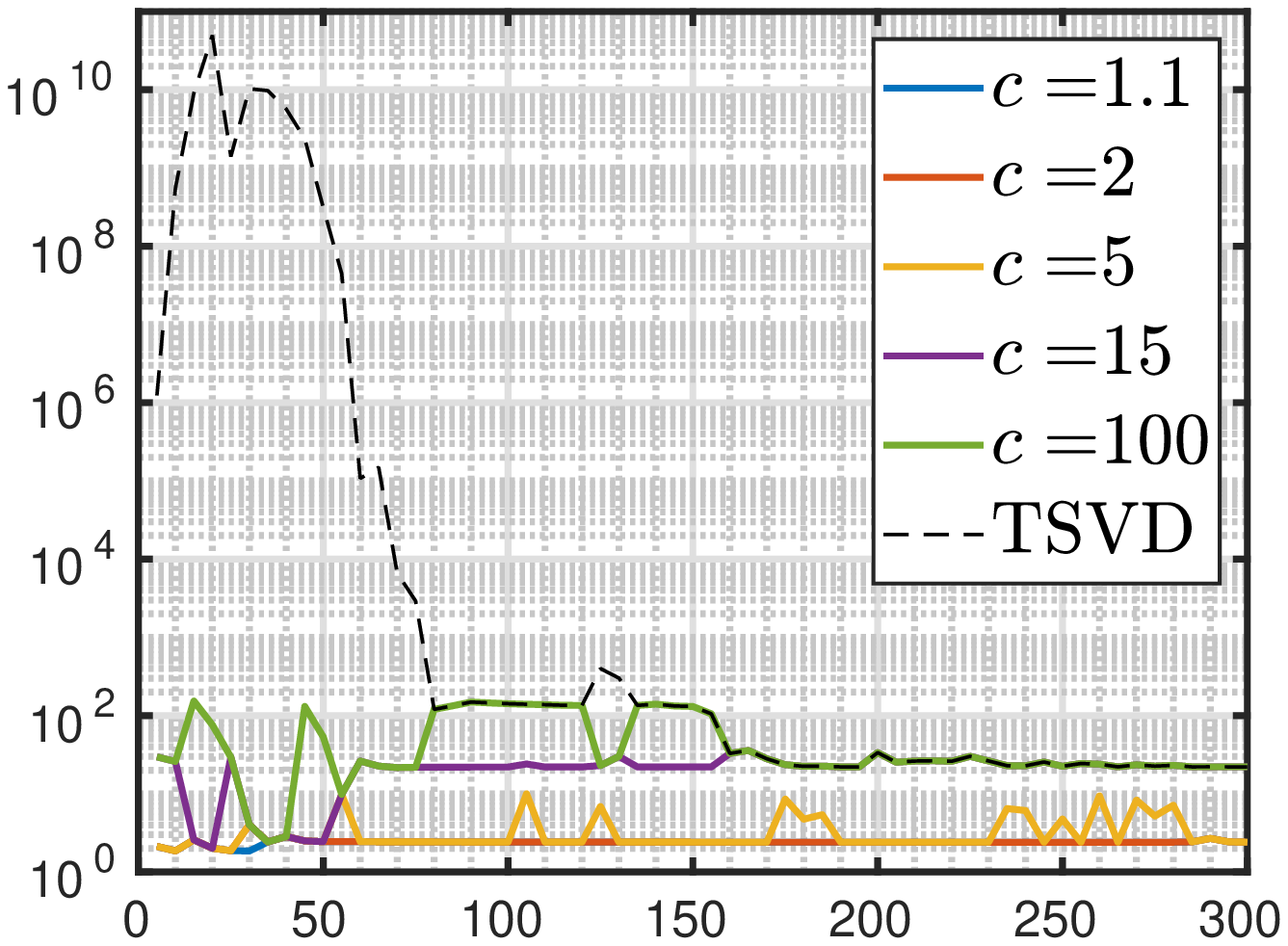} &
\includegraphics[width=\four]{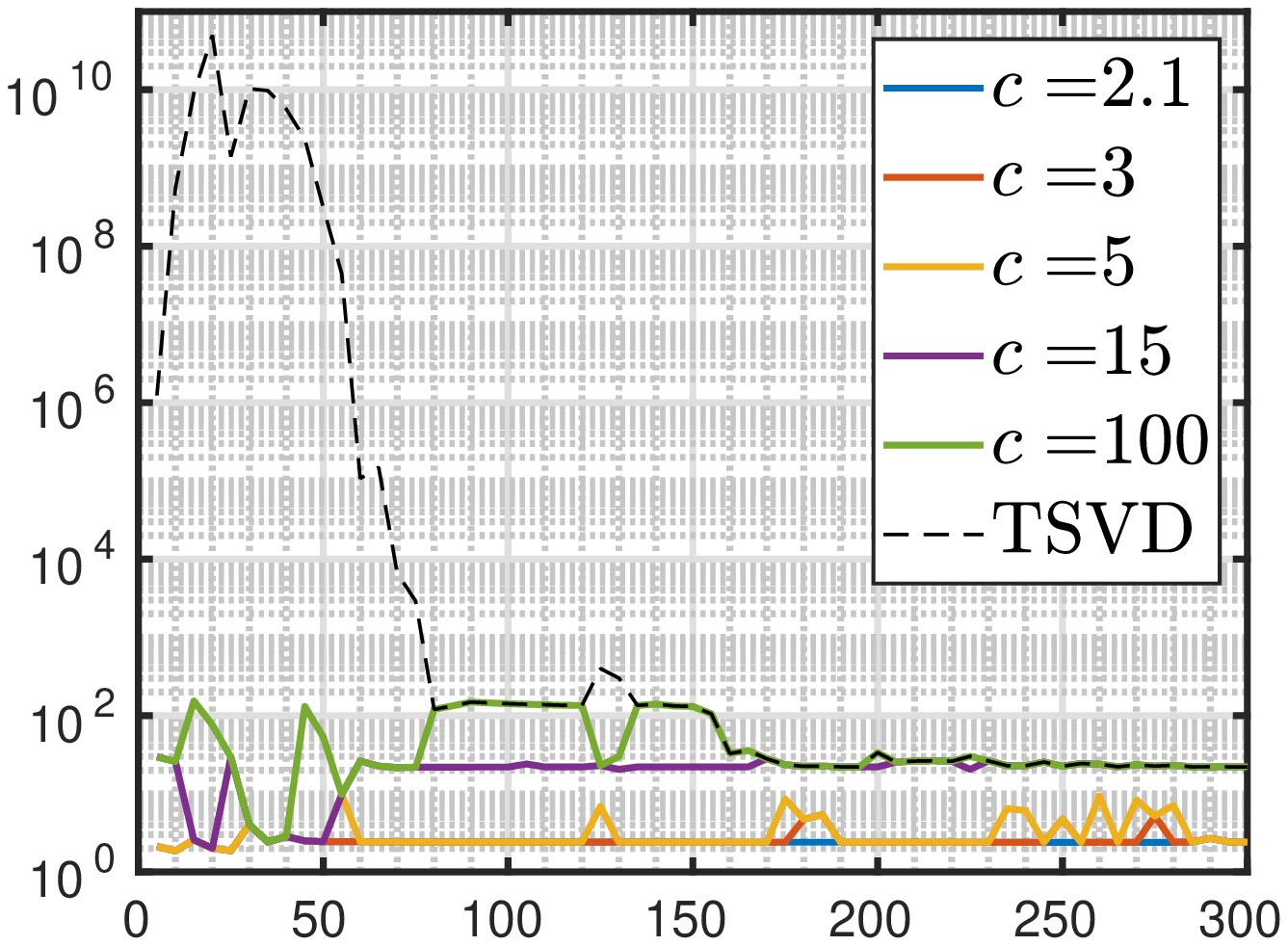}
\end{tabular}
}
\end{center}
\vspace*{-5mm}\caption{
Comparison of ASVD1 and ASVD2 for various values of the parameter $c$. The methods are applied to the system \R{LSsing} with data \R{Datasing} and $M/N = 2$. Left to right: $L^2$-norm error ASVD1, $L^2$-norm error for ASVD2, coefficient norm for ASVD1, coefficient norm for ASVD2. The function considered is $f(t) = \E^{\sin(15 t + 0.5)} + \log(t) \cos(\alpha  t)$ for $\alpha = 1$ (top), $\alpha = 20$ (middle) and $\alpha = 40$ (bottom). The threshold $\epsilon$ is set to $10^{-15}$ in all cases.
}
\label{fig:SingFNA1ASVD1}
\end{figure}

Having provided the rationale for using this frame, we now present numerical results. Rather than solving \R{eq: Best approximation}, we consider the setup of \cite{BADHFramesPart2} (see Remark \ref{r:FNA2setup}) and compute an approximation in $\rH_N$ using the data
\be{
\label{Datasing}
\bm{y} = \{ f(t_{m,M}) \}^{M}_{m=1},\qquad t_{m,M} = \frac{\cos(\frac{2m-1}{2M} \pi ) + 1}{2},\qquad m = 1,\ldots,M.
}
Note that the points $t_{m,M}$ are Chebyshev nodes in $(0,1)$. To do this, we simply modify the TSVD, ASVD1 and ASVD2 methods to approximately solve the least-squares problem
\be{
\label{LSsing}
\bm{G}_{M,N} \bm{z} \approx \bm{y},\qquad \bm{G}_{M,N} = \left \{ \phi_n(t_{m,M}) \right \}^{M,N-1}_{m=1,n=0} \in \bbC^{M \times N},
}
In Fig.\ \ref{fig:SingFNA1ASVD1} we compare these methods for this problem with the oversampling $M/N = 2$.

This experiment reveals several interesting phenomena. First, when $c$ is small, the convergence of ASVD1 and ASVD2 suffers, in particular when $\alpha \gg 1$ is larger. This is as expected. For the function considered, the $K^{\rth}$ derivative of $g(t) = \log(t) \cos(\alpha t)$ scales like $\alpha^{K}$. When $\alpha$ is large, the coefficients constructed above might not be attainable by ASVD1 or ASVD2. The slower rates of convergence shown in the figure support this argument. Second, for small $N$, ASVD1 and ASVD2 actually achieve slightly smaller errors than TSVD. Informally, this might be explained as follows: by constraining the singular values, we prohibit larger-norm solutions, which asymptotically might give the optimal rate of convergence, but are associated with larger constants in the error bounds.

\section{Theoretical analysis}
\label{section : theory}
We now present our main theoretical results. As with TSVD, we aim to establish estimates of the form \R{TSVDerr}--\R{TSVDcoefflimit}. Proofs of the results in this section are found in \S \ref{s:proofs}.

\subsection{Main results}

We commence with the error bounds:

\thm{
\label{t:errbds}
Let $f \in \rH$ and $\cP_{\Lambda} f$ be its frame approximation corresponding to TSVD, ASVD1 or ASVD2. Then, for TSVD,
\be{
\label{TSVD_err}
\Vert f-\mathcal{P}_{\Lambda}f\Vert\leq\Vert f-\mathcal{T}_{N}\bm{z}\Vert+\sqrt{\epsilon}\Vert\bm{z}\Vert,\qquad\forall\bm{z}\in\mathbb{C}^{N},
}
for ASVD1,
\be{
\label{ASVD1_err}
\Vert f-\mathcal{P}_{\Lambda}f\Vert\leq\Vert f-\mathcal{T}_{N}\bm{z}\Vert+\max\left\{\sqrt{\epsilon}, \frac{\Vert f-\mathcal{T}_{N}\bm{z}\Vert}{c\nm{\bm{y}}-\Vert\bm{z}\Vert}\right\}\Vert\bm{z}\Vert,\qquad\forall\bm{z}\in\mathbb{C}^{N}, \ \nm{\bm{z}} < c \nm{\bm{y}},
}
and for ASVD2,
\be{
\label{ASVD2_err}
\begin{split}
\Vert f-\mathcal{P}_{\Lambda}f\Vert\leq  \Vert f-\mathcal{T}_{N}\bm{z}\Vert+&\max\left\{\sqrt{\epsilon},\frac{\Vert f-\mathcal{T}_{N}\bm{z}\Vert}{c\nm{\bm{y}}-\left(\frac{1}{\sqrt{\epsilon}}\right)\Vert f-\mathcal{T}_{N}\bm{z}\Vert-2\Vert\bm{z}\Vert}\right\}\Vert\bm{z}\Vert,
\\
& \forall\bm{z}\in\mathbb{C}^{N},\ \left(1/\sqrt{\epsilon}\right) \Vert f-\mathcal{T}_{N}\bm{z}\Vert + 2 \Vert\bm{z}\Vert < c \nm{\bm{y}}.
\end{split}
}
}
This result explains the effect of constraining $\nmu{\bm{x}^{\Lambda}}$.  While \R{TSVD_err} allows for arbitrary coefficients $\bm{z}$, in \R{ASVD1_err} and \R{ASVD2_err} the error is controlled only by those $\bm{z}$ that satisfy a certain bound. It is notable that the bound on $\bm{z}$ for ASVD2 is much stricter than it is for ASVD1. 
We discuss further consequences of Theorem \ref{t:errbds} in a moment. First, we consider the coefficient norms:

\thm{
\label{t:coeffbds}
Let $f \in \rH$ and $\bm{x}^{\Lambda}$ be the coefficients of the frame approximation corresponding to TSVD, ASVD1 or ASVD2. Then, for TSVD,
\be{
\label{TSVD_coeff}
\nm{\bm{x}^{\Lambda}} \leq\Vert f-\mathcal{T}_{N}\bm{z}\Vert/\sqrt{\epsilon}+\Vert\bm{z}\Vert,\qquad\forall\bm{z}\in\mathbb{C}^{N},
}
for ASVD1,
\be{
\label{ASVD1_coeff}
\nm{\bm{x}^{\Lambda}} \leq \min \left \{ c \sqrt{B} \sqrt{N} \nm{f},\Vert f-\mathcal{T}_{N}\bm{z}\Vert/\sqrt{\epsilon}+\Vert\bm{z}\Vert \right \},\qquad\forall\bm{z}\in\mathbb{C}^{N},
}
and for ASVD2,
\be{
\label{ASVD2_coeff}
\nm{\bm{x}^{\Lambda}} \leq \min \left \{ c \sqrt{B}\nm{f},\Vert f-\mathcal{T}_{N}\bm{z}\Vert/\sqrt{\epsilon}+\Vert\bm{z}\Vert \right \},\qquad\forall\bm{z}\in\mathbb{C}^{N}.
}
}

This agrees with the numerical results. The coefficients $\bm{x}^{\Lambda}$ of ASVD1 or ASVD2 remain bounded (or nearly bounded) in the pre-asymptotic regime, and then match the coefficient bound for TSVD asymptotically as $N \rightarrow \infty$.

\subsection{Error behaviour}

We now give further insight into Theorem \ref{t:errbds}. To do so, we define the following:

\defn{
Let $a> 0$ and $0 < \delta < 1$.  Then $f \in \rH$ has the $\left(a,\delta\right)$-stable approximation property if there exists a $\bm{z} \in \bbC^N$ such that $\nm{f - \cT_N \bm{z}} \leq \delta \nm{f}$ and $\nm{\bm{z}} \leq a \nm{f}$.
}

\prop{
Let $a, \delta > 0$ with $a \geq 1/\sqrt{A}$ and $f \in \rH$.  Then there exists an $N_0 = N_0(a,\delta,f)$ such that $f$ has  the $\left(a,\delta\right)$-stable approximation property for all $N \geq N_0$.
}
\prf{
Recall that $f = \sum_{n \in I} a_n \phi_n$, where $a_n$ are the frame coefficients, and this series converges in $\rH$. Hence, there is an $N_0$ such that if $N \geq N_0$ and $\bm{a}_N = (a_n)_{n \in I_N}$ then $\nm{f - \cT_N \bm{a}_N}  \leq \delta \nm{f}$.
Moreover,  $\Vert \bm{a}_{N} \Vert^{2} \leq \nm{\bm{a}}^2 \leq 1/A \nm{f}^2$ by \R{dualframebd}. This gives the result.
}

\thm{
\label{thm:stableapproxprop}
Let $f\in \rH$ have $\left(a,\delta\right)$-stable approximation property.  Then for the TSVD,
\be{
\nm{f - \cP_{\Lambda} f } \leq (\delta + a \sqrt{\epsilon} ) \Vert f \Vert , 
}
and for ASVD1,
\be{
\label{ASVD1adelta}
\nm{f - \cP_{\Lambda} f } \leq \left ( \delta +  \max \left \{a \sqrt{\epsilon} , \frac{a^2}{c(1-\delta) - a^2} \delta \right \} \right ) \nm{f},\qquad \forall c > \frac{a^2}{1 - \delta}.
}
Further, write $\delta = b \sqrt{\epsilon}$ for $b>0$. Then, for ASVD2,
\be{
\nm{f - \cP_{\Lambda} f } \leq \left (b + \max \left \{ a , \frac{a^2 b}{c \left(1- b\sqrt{\epsilon} \right) - ab - 2a^2} \right \} \right ) \sqrt{\epsilon} \Vert f \Vert,\qquad \forall c > \frac{a \left( 2 a + b \right)}{1 - b\sqrt{\epsilon}}.
}
}

Notice that if we replace the maximum by a sum, then the ASVD1 error \R{ASVD1adelta} can be replaced by the simpler expression
\bes{
\nm{f - \cP_{\Lambda} f } \leq \left ( \frac{c(1-\delta)}{c(1-\delta)-a^2} \delta + a \sqrt{\epsilon} \right )\nm{f},\qquad \forall c > \frac{a^2}{1 - \delta}.
}
Hence this theorem shows that whenever there is sequence of coefficients $\bm{z}$ approximating $f$ to a relative error of $\delta$, and which do not grow too large (depending on $c$), then TSVD and ASVD1 are guaranteed the same error bound up to a constant depending on $c$, $\delta$ and $a$.
Larger $c$ allows means this holds for bigger $a$, thus, as expected, making it easier in some sense to achieve an error of size $\delta$. Notice that for ASVD2 this result only applies when $\delta = \ord{\sqrt{\epsilon}}$. This is due to the condition on $\bm{z}$ in Theorem \ref{t:errbds}, and suggests that ensuring fast error decay with ASVD2 may be more challenging, especially for smaller $c$.

It is informative to examine this result for the Legendre polynomial frame of \S \ref{ss:mainexamp} . Recall that \R{LegFrameSob} asserts the existence of bounded coefficients yielding algebraic rates of convergence. Setting $\delta = D_{k,f} N^{-k}$, where $D_{k,f} = C_k \nm{f}_{H^k} / \nm{f}_{L^2}$ and $a = D_{k,f}$, it follows that 
\bes{
\nm{ f - \cP_{\Lambda} f} \leq  \left ( N^{-k} + \sqrt{\epsilon} \right ) D_{k,f} \nm{f}  ,
}
for TSVD. For ASVD1, if $N \geq (2 D_{k,f})^{1/k}$ (the factor of $2$ is arbitrary), we have
\bes{
\nm{ f - \cP_{\Lambda} f} \leq \left (  \frac{c}{c - 2 D^2_{k,f}}  N^{-k} + \sqrt{\epsilon} \right ) D_{k,f} \nm{f},\qquad c > 2 D^2_{k,f},
}
This agrees with the numerical examples. For sufficiently large $c$, we expect little or no deterioration in the rate of error decay. But for small $c$, in particular, when $\nm{f}_{H^k}$ grows rapidly with $k$, we may see slower decay due to unattainability in ASVD1 of the coefficient vectors giving the faster algebraic rates.

\subsection{Limiting behaviour}

To conclude this section, we consider the limiting behaviour of each method:

\thm{ 
\label{thm:limsup}
For either TSVD, ASVD1 with $c > 1/A$ or ASVD2 with  $c > 2/A$ the approximation satisfies
\begin{enumerate}[label=(\roman*)]
\item \label{i}$\limsup_{N \rightarrow \infty} \nm{\bm{x}_{\Lambda}} \leq \nm{\bm{a}},$
\item \label{ii}$\limsup_{N \rightarrow \infty} \nm{\bm{a} - \bm{a}_{\Lambda}} \leq \sqrt{\epsilon/A}\nm{\bm{a}},$ 
\item \label{iii}$\limsup_{N \rightarrow \infty} \nm{f - \cP_{\Lambda} f } \leq \sqrt{\epsilon} \nm{\bm{a}},$
\end{enumerate}
where $\bm{a} = \{ \langle f,\mathcal{S}^{-1}\phi_{n} \rangle \}_{n \in I}$ are the frame coefficients of $f$ and $\bm{a}_{\Lambda} \in \ell^2(I)$ is the extension of $\bm{x}_{\Lambda}$ by zero.
}

Recall that $\nm{\bm{a}} \leq \nm{f} / A$ by \R{dualframebd}.
This result agrees with our numerical examples.  For all three methods, the approximation eventually reaches within $\ord{\sqrt{\epsilon}}$ of $f$, and in the case of ASVD1 and ASVD2 this occurs only when $c$ is larger than a constant depending on the frame bound. Recall that $A = 1$ for the Legendre polynomial frame. The frame introduced in \S \ref{ss:numexamp2} also has lower frame bound $A = 1$, since it contains an orthonormal basis.

\section{Proofs of the main results}\label{s:proofs}

To commence, we require the following two lemmas:

\lem{
\label{err_lem}
For any fixed $\Lambda\subseteq I_{N}$ the orthogonal projection $\mathcal{P}_{\Lambda}$ satisfies
\be{
\label{err_bound_lemma}
\Vert f-\mathcal{P}_{\Lambda}f\Vert\leq\inf\{\Vert f-\mathcal{T}_{N}\bm{z}\Vert+\max_{n\in I_N \backslash \Lambda}\{\sqrt{\sigma_{n}}\}\Vert\bm{z}\Vert:\bm{z}\in\mathbb{C}^{N}\}.
}
}

\prf{
Let $\bm{z}\in\mathbb{C}^{N}$. Since $\mathcal{P}_{\Lambda}f$ is an orthogonal projection onto $\rH_{\Lambda}$, we have 

\bes{
\Vert f-\mathcal{P}_{\Lambda}f\Vert\leq\Vert f-\mathcal{P}_{\Lambda}\mathcal{T}_{N}\bm{z}\Vert\leq\Vert f-\mathcal{T}_{N}\bm{z}\Vert+\Vert\mathcal{T}_{N}\bm{z}-\mathcal{P}_{\Lambda}\mathcal{T}_{N}\bm{z}\Vert.
}
Note that $\mathcal{T}_{N}\bm{z}=\mathcal{P}_{N}\mathcal{T}_{N}\bm{z}$ since $\mathcal{T}_{N}\bm{z}\in \rH_{N}$. Hence (\ref{orthog_proj}) and the orthogonality of the $\xi_{n}$'s gives

\bes{
\Vert\mathcal{T}_{N}\bm{z}-\mathcal{P}_{\Lambda}\mathcal{T}_{N}\bm{z}\Vert^{2}\leq\left\Vert\sum_{n\in I_N \backslash \Lambda}\frac{\langle\mathcal{T}_{N}\bm{z},\xi_{n}\rangle}{\sigma_{n}}\xi_{n}\right\Vert^{2}=\sum_{n\in I_N \backslash \Lambda}\frac{\vert\langle\mathcal{T}_{N}\bm{z},\xi_{n}\rangle\vert^{2}}{\sigma_{n}}.
} 
Observe that $\langle\mathcal{T}_{N}\bm{z},\xi_{n}\rangle=\langle\mathcal{T}_{N}\bm{z},\mathcal{T}_{N}\bm{v}_{n}\rangle=\langle\bm{z},\bm{G}_{N}\bm{v}_{n}\rangle=\sigma_{n}\langle\bm{z},\bm{v}_{n}\rangle$ and therefore

\bes{
\Vert\mathcal{T}_{N}\bm{z}-\mathcal{P}_{\Lambda}\mathcal{T}_{N}\bm{z}\Vert^{2}=\sum_{n\in I_N \backslash \Lambda}\sigma_{n}\vert\langle\bm{z},\bm{v}_{n}\rangle\vert^{2}\leq\max_{n\in I_N \backslash \Lambda}\{\sigma_{n}\}\sum_{n\in  I_{N}}\vert\langle\bm{z},\bm{v}_{n}\rangle\vert^{2}=\max_{n\in I_N \backslash \Lambda}\{\sigma_{n}\}\Vert\bm{z}\Vert^{2},
}
where in the last step we use the fact that the vectors $\{\bm{v}_{n}\}_{n\in I_N}$ are orthonormal. 
}
\lem{
The coefficients $\bm{x}_{\Lambda}$ of the orthogonal projection $\mathcal{P}_{\Lambda}$ satisfy

\be{
\label{coef_bound_lemma}
\Vert\bm{x}_{\Lambda}\Vert\leq\inf{\left\{ \left(1/\min_{n\in\Lambda}{\left\{\sqrt{\sigma_n}\right\}}\right)
\Vert f-\mathcal{T}_{N}\bm{z}\Vert+\Vert\bm{z}\Vert:\bm{z}\in\mathbb{C}^{N}\right\}}.
}
Moreover, if $\bm{a}_{\Lambda}\in\ell^{2}(I)$ is the extension of $\bm{x}_{\Lambda}$ by zero, then 
\be{
\label{coef_conv_lemma}
\Vert \bm{a} - \bm{a}_{\Lambda} \Vert \leq \left( 1 + \sqrt{B}/\min_{n\in\Lambda}{\left\{ \sqrt{\sigma_{n}} \right\}} \right) \sqrt{\sum_{n\in I\setminus I_{N}} \vert a_{n} \vert^{2}} + \left( \max_{n\in I_N \backslash \Lambda}{\left\{ \sqrt{\sigma_{n}} \right\}}/\sqrt{A} \right) \Vert \bm{a} \Vert
}
where $\bm{a}=\left\{\langle f,\mathcal{S}^{-1}\phi_{n}\rangle\right\}_{n\in I}$ are the frame coefficients of $f$.
}

\prf{
For the first part, we use (\ref{coef_formula}) to write
\bes{
\bm{x}_{\Lambda}=\sum_{n\in \Lambda}\frac{\langle f,\xi_{n}\rangle}{\sigma_{n}}\bm{v}_{n}=\sum_{n\in \Lambda}\frac{\langle f-\mathcal{T}_{N}\bm{z},\xi_{n}\rangle}{\sigma_{n}}\bm{v}_{n}+\sum_{n\in \Lambda}\frac{\langle \mathcal{T}_{N}\bm{z},\xi_{n}\rangle}{\sigma_{n}}\bm{v}_{n}.
}
Consider the first term on the right-hand side. By (\ref{prolate_orthog}) and (\ref{orthog_proj}) we have 
\bes{
\left\Vert\sum_{n\in \Lambda}\frac{\langle f-\mathcal{T}_{N}\bm{z},\xi_{n}\rangle}{\sigma_{n}}\bm{v}_{n}\right\Vert^{2}=\sum_{n\in \Lambda}\frac{\vert\langle f-\mathcal{T}_{N}\bm{z},\xi_{n}\rangle\vert^{2}}{\sigma_{n}^{2}}\leq\left(1/\min_{n\in\Lambda}\left\{{\sqrt{\sigma_n}}\right\}\right)\Vert \mathcal{P}_{\Lambda}\left(f-\mathcal{T}_{N}\bm{z}\right)\Vert,
}
and hence 
\bes{
\left\Vert\sum_{n\in \Lambda}\frac{\langle f-\mathcal{T}_{N}\bm{z},\xi_{n}\rangle}{\sigma_{n}}\bm{v}_{n}\right\Vert^{2}\leq\left(1/\min_{n\in\Lambda}\left\{{\sqrt{\sigma_n}}\right\}\right)\Vert f-\mathcal{T}_{N}\bm{z}\Vert.
}
For the second term, we notice that $\langle\mathcal{T}_{N}\bm{z},\xi_{n}\rangle=\sigma_{n}\langle\bm{z},\bm{v}_{n}\rangle$, and therefore
\bes{
\left\Vert\sum_{n\in \Lambda}\frac{\langle \mathcal{T}_{N}\bm{z},\xi_{n}\rangle}
{\sigma_{n}}\bm{v}_{n}\right\Vert^{2}=\sum_{n\in \Lambda}\left\vert\langle\bm{z},\bm{v}_{n}\rangle\right\vert^{2}\leq\Vert\bm{z}\Vert^{2}.
}
Combining these two bounds now gives the first result. For the second result, we first let $\bm{a}_{N}\in\mathbb{C}^{N}$ be the vector with $n^{\text{th}}$ entry $a_{n}=\langle f,\mathcal{S}^{-1}\phi_{n}\rangle$ for $n\in I_{N}$. Then 
\bes{
\Vert\bm{a}-\bm{a}_{\Lambda}\Vert\leq\sqrt{\sum_{n\in I\setminus I_{N} }\vert a_{n}\vert^{2}}+\Vert\bm{a}_{N}-\bm{x}_{\Lambda}\Vert.
}
Hence it suffices to estimate $\Vert\bm{a}_{N}-\bm{x}_{\Lambda}\Vert$. For this, we note that  
$f=\mathcal{S}\mathcal{S}^{-1}f=\mathcal{S}_{N}\mathcal{S}^{-1}f+\left(\mathcal{S}-\mathcal{S}_{N}\right)\mathcal{S}^{-1}f$.
Since $\mathcal{S}_{N}$ is self-adjoint and $\mathcal{S}_{N}\xi{n}=\mathcal{T}_{N}\mathcal{T}_{N}^{*}\mathcal{T}_{N}\bm{v}_{n}=\sigma_{n}\mathcal{T}_{N}\bm{v}_{n}=\sigma_{n}\xi_{n}$ we have
\bes{
\langle f,\xi_{n}\rangle=\langle\mathcal{S}_{N}\mathcal{S}^{-1}f,\xi_{n}\rangle+\langle\left(\mathcal{S}-\mathcal{S}_{N}\right)\mathcal{S}^{-1}f,\xi_{n}\rangle=\sigma_{n}\langle\mathcal{S}^{-1}f,\xi_{n}\rangle+\langle\left(\mathcal{S}-\mathcal{S}_{N}\right)\mathcal{S}^{-1}f,\xi_{n}\rangle.
}
Therefore
\be{
\label{splited_coef}
\bm{x}_{\Lambda}=\sum_{n\in \Lambda}\frac{\langle f,\xi_{n}\rangle}{\sigma_{n}}\bm{v}_{n}=\sum_{n\in \Lambda}\langle\mathcal{S}^{-1}f,\xi_{n}\rangle\bm{v}_{n}+\sum_{n\in \Lambda}\frac{1}{\sigma_{n}}\langle\left(\mathcal{S}-\mathcal{S}_{N}\right)\mathcal{S}^{-1}f,\xi_{n}\rangle\bm{v}_{n}.
}	
Conversely, since $\bm{a}_{N}=\mathcal{T}_{N}^{*}\mathcal{S}^{-1}f$ we have $\langle\bm{a}_{N},\bm{v}_{n}\rangle=\langle\mathcal{S}^{-1}f,\mathcal{T}_{N}\bm{v}_{n}\rangle=\langle\mathcal{S}^{-1}f,\xi_{n}\rangle.$ Hence
\be{
\label{truncated_frame_coef}
\bm{a}_{N}=\sum_{n\in_{I_N}}\langle\bm{a}_{N},\bm{v}_{n}\rangle\bm{v}_{n}=\sum_{n\in {I_N}}\langle\mathcal{S}^{-1}f,\xi_{n}\rangle\bm{v}_{n}.
}
Combining (\ref{splited_coef}) and (\ref{truncated_frame_coef}) now gives
\be{
\label{coef_resid}
\Vert\bm{a}_{N}-\bm{x}_{\Lambda}\Vert\leq\left\Vert\sum_{n\in I_N \backslash \Lambda}\langle\mathcal{S}^{-1}f,\xi_{n}\rangle\bm{v}_{n}\right\Vert+\left\Vert\sum_{n\in\Lambda}\frac{1}{\sigma_{n}}\langle\left(\mathcal{S}-\mathcal{S}_{N}\right)\mathcal{S}^{-1}f,\xi_{n}\rangle\bm{v}_{n}\right\Vert
}
Consider the first term. By orthogonality
\be{
\label{first_term_cof_conv}
\begin{split}
\left\Vert\sum_{n\in I_N \backslash \Lambda}\langle\mathcal{S}^{-1}f,\xi_{n}\rangle\bm{v}_{n}\right\Vert \leq  \max_{n\in I_N \backslash \Lambda}\left\{\sigma_{n}\right\}\sum_{n\in I_N \backslash \Lambda}\frac{1}{\sigma_{n}}\vert\langle\mathcal{S}^{-1}f,\xi_{n}\rangle\vert^{2}
\leq &
\max_{n\in I_N \backslash \Lambda}\left\{\sigma_{n}\right\}\Vert\mathcal{S}^{-1}f\Vert^{2}
\\
\leq & \left(\max_{n\in I_N \backslash \Lambda}{\left\{\sigma_{n}\right\}}/A\right)\Vert\bm{a}\Vert^{2}.
\end{split}
}
Now consider the second term:
\eas{
\left\Vert\sum_{n\in\Lambda}\frac{1}{\sigma_{n}}\langle\left(\mathcal{S}-\mathcal{S}_{N}\right)\mathcal{S}^{-1}f,\xi_{n}\rangle\bm{v}_{n}\right\Vert^{2} & = \sum_{n\in\Lambda}\frac{1}{\sigma_{n}^{2}}\left\vert\langle\left(\mathcal{S}-\mathcal{S}_{N}\right)\mathcal{S}^{-1}f,\xi_{n}\rangle\right\vert^{2}
\\
& \leq  \left(1/\min_{n\in\Lambda}{\left\{ \sigma_{n}\right\}}\right)\sum_{n\in\Lambda}\frac{1}{\sigma_{n}}\left\vert\langle\left(\mathcal{S}-\mathcal{S}_{N}\right)\mathcal{S}^{-1}f,\xi_{n}\rangle\right\vert^{2}
\\
& \leq \left(1/\min_{n\in\Lambda}{\left\{ \sigma_{n}\right\}}\right)\left\Vert\left(\mathcal{S}-\mathcal{S}_{N}\right)\mathcal{S}^{-1}f\right\Vert .
}
Observe that 
\bes{
\left\Vert\left(\mathcal{S}-\mathcal{S}_{N}\right)\mathcal{S}^{-1}f\right\Vert=\left\Vert\sum_{n\in I\setminus I_{N}}a_{n}\phi_{n}\right\Vert=\sup_{\substack{g\in \rH \\ g\neq 0}}\left\{\frac{\left\vert\sum_{n>N}a_{n}\overline{\langle g,\phi_{n}\rangle}\right\vert}{\Vert g\Vert}\right\}\leq\sqrt{B}\sqrt{\sum_{n\in I\setminus I_{N}}\vert a_{n}\vert^{2}},
}
and therefore
\bes{
\left\Vert\sum_{n\in\Lambda}\frac{1}{\sigma_{n}}\langle\left(\mathcal{S}-\mathcal{S}_{N}\right)\mathcal{S}^{-1}f,\xi_{n}\rangle\bm{v}_{n}\right\Vert^{2}\leq\left(B/\min_{n\in\Lambda}{\left\{\sigma_{n}\right\}}\right)\sum_{n\in I\setminus I_{N}}\vert a_{n}\vert^{2}.
}
Substituting this and (\ref{first_term_cof_conv}) into (\ref{coef_resid}) gives the result.
}

\prf{[Proof of Theorem \ref{t:errbds}]
Utilizing Lemma \ref{err_lem}, it suffices to estimate the largest singular value being discarded by each method. For the TSVD, note that one simply keeps every singular value greater than $\epsilon$. Consequently, the largest singular value being discarded is strictly smaller than or equal to $\epsilon$ and the result directly follows by substituting $\epsilon$ into (\ref{err_bound_lemma}). 

For ASVD1, observe that if $m \in I_{N} \setminus \Lambda$ and $\sigma_{m} > \epsilon$, then it must be the case that $\vert \langle \bm{y},\bm{v}_{m} \rangle \vert / \sigma_{m} > c \nm{\bm{y}}$,
and therefore
\be{
\label{ASVD1_sing_intermeidate}
\sigma_{m} < \frac{\vert \langle \bm{y},\bm{v}_{m} \rangle \vert}{c \nm{\bm{y}}}.
}
Consider the numerator. Recall that $\langle\bm{y},\bm{v}_{n}\rangle=\langle f,\xi_{n}\rangle$ for every $n\in I_{N}${. T}herefore
\be{
\label{ASVD1_triang}
{\vert \langle \bm{y},\bm{v}_{m} \rangle \vert} = \vert \langle f, \xi_{m} \rangle \vert \leq \vert \langle  f- \mathcal{T}_{N}\bm{z}, \xi_{m} \rangle \vert + \vert \langle \mathcal{T}_{N}\bm{z}, \xi_{m} \rangle \vert, \qquad \forall \bm{z} \in \mathbb{C}^{N}.
}
Consider the first term on the right-hand side. Since $\Vert\xi_{n}\Vert=\sqrt{\sigma_{n}}$ by (\ref{prolate_orthog}) we have
\be{
\label{ASVD1_first-term}
\vert \langle f - \mathcal{T}_{N}\bm{z}, \xi_{m} \rangle \vert \leq \sqrt{\sigma_{m}} \Vert f - \mathcal{T}_{N}\bm{z} \Vert, \qquad \forall \bm{z} \in \mathbb{C}^{N}.
}
For the second term on the right-hand side, observe that
\be{
\label{ASVD1_second-term}
\langle \mathcal{T}_{N}\bm{z}, \xi_{{m}} \rangle = \langle \bm{z}, \mathcal{T}_{N}^{*}\xi_{{m}} \rangle = \langle \bm{z}, \mathcal{T}_{N}^{*}\mathcal{T}_{N}\bm{v}_{{m}} \rangle = \langle \bm{z}, \bm{G}_{N} \bm{v}_{{m}} \rangle =  \sigma_{{m}} \langle \bm{z},  \bm{v}_{{m}} \rangle, \qquad \forall \bm{z} \in \mathbb{C}^{N}.
} 
Combining (\ref{ASVD1_first-term}) and (\ref{ASVD1_second-term}) with (\ref{ASVD1_triang}) now gives
\be{
\label{split-sing_val}
\vert \langle \bm{y}, \bm{v}_{m} \rangle \vert \leq \sqrt{\sigma_{m}} \Vert f - \mathcal{T}_{N}\bm{z} \Vert + \sigma_{m} \Vert \bm{z} \Vert, \qquad \forall \bm{z} \in \mathbb{C}^{N}.
}
Therefore
\bes{
\sigma_{m} \leq \frac{\sqrt{\sigma_{m}} \Vert f - \mathcal{T}_{N}\bm{z} \Vert + \sigma_{m} \Vert \bm{z} \Vert}{c \nm{\bm{y}}}, \qquad \forall \bm{z} \in \mathbb{C}^{N}.
}
Rearranging terms and simplifying gives
\be{
\label{ASVD1_sing_final}
\sqrt{\sigma_{m}} \leq \frac{\Vert f- \mathcal{T}_{N}\bm{z} \Vert}{c \nm{\bm{y}} - \Vert \bm{z} \Vert}, \qquad \forall \bm{z} \in \mathbb{C}^{N},\ \nm{\bm{z}} < c \nm{\bm{y}}.
}
Substituting (\ref{ASVD1_sing_final}) into (\ref{err_lem}) gives the result for the ASVD1. 

For the ASVD2, observe that if $m \in I_{N} \setminus \Lambda$ and $\sigma_{m} > \epsilon$, then it must be the case that
\bes{
\Vert \bm{x}_{\Lambda} \Vert^{2} + \frac{\vert \langle \bm{y},\bm{v}_{m} \rangle \vert^{2}}{{\sigma_{m}}^{2}} > c^{2} \nm{\bm{y}}^{2},
}
and therefore
\bes{
\Vert \bm{x}_{\Lambda} \Vert + \frac{\vert \langle \bm{y},\bm{v}_{m} \rangle \vert}{\sigma_{m}} > c \nm{\bm{y}}.
}
It directly follows from Lemma \ref{coef_bound_lemma} that
\be{
\label{ASVD2_lemma2}
\Vert \bm{x}_{\Lambda} \Vert \leq \left(1/\sqrt{\epsilon}\right) \Vert f-\mathcal{T}_{N}\bm{z} \Vert
 + \Vert \bm{z} \Vert, \qquad \forall \bm{z} \in \mathbb{C}^{N}.
}
Combining this with (\ref{split-sing_val}) gives
\bes{
\left(1/\sqrt{\epsilon}\right) \Vert f-\mathcal{T}_{N}\bm{z} \Vert
 + \Vert \bm{z} \Vert + \frac{\sqrt{\sigma_m} \nm{f-  \cT_N \bm{z}} + \sigma_m \nm{\bm{z}}}{\sigma_m} > c \nm{\bm{y}}.
}
Simplifying and rearranging yields
\ea{
\label{ASVD2_sing_final} \nonumber
\sqrt{\sigma_{m}} < & \frac{\Vert f - \mathcal{T}_{N} \bm{z}\Vert}{c \nm{\bm{y}} - \left(1/\sqrt{\epsilon}\right) \Vert f-\mathcal{T}_{N}\bm{z}\Vert - 2 \Vert\bm{z}\Vert},
\\ 
& \forall \bm{z} \in \mathbb{C}^{N},\ \left(1/\sqrt{\epsilon}\right) \Vert f-\mathcal{T}_{N}\bm{z}\Vert + 2 \Vert\bm{z}\Vert < c \nm{\bm{y}}.
}
Substituting this into (\ref{err_lem}) finishes the proof. 
}

\prf{[Proof of Theorem \ref{t:coeffbds}]
This follows immediately from Lemma \ref{coef_bound_lemma} and, for ASVD1 and ASVD2, the discussion immediately after Definitions \ref{d:ASVD1} and \ref{d:ASVD2} respectively.
}

In order to prove the next result, we first observe the following. Let $f \in \rH$ and $\bm{y} = \cT^*_N f$. Then, for any $\bm{z} \in \bbC^N \backslash \{0 \}$,
\bes{
\nm{\bm{y}} \geq \frac{\left\vert \langle \bm{y},\bm{z} \rangle \right\vert}{\nm{\bm{z}}} = \frac{\left\vert \langle f,\mathcal{T}_{N}\bm{z} \rangle \right\vert}{\nm{\bm{z}}}= \frac{\left\vert \nm{f}^{2} + \langle f,f - \mathcal{T}_{N}\bm{z} \rangle \right\vert}{\nm{\bm{z}}} \geq \frac{\nm{f}^2 - | \ip{f}{f - \cT_{N} \bm{z}} |}{\nm{\bm{z}}}.
}
We deduce the following inequality:
\be{
\label{(f,y)-NormConv}
\nm{\bm{y}} \geq \frac{\nm{f}^2 - \nm{f - \cT_{N} \bm{z}} \nm{f} }{\nm{\bm{z}}},\qquad \forall \bm{z} \in \bbC^N \backslash \{ 0 \},f \in \rH,\ \mbox{where}\ \bm{y} = \cT^*_N f.
}
%
%

\prf{
[Proof of Theorem \ref{thm:stableapproxprop}]
Since $f$ has $\left(a,\delta\right)$-stable approximation property, by definition, there exists a vector $\bm{z}\in\mathbb{C}^{N}$ such that
\be{\label{pr_stable}
\nm{f - \cT_N \bm{z}} \leq \delta \nm{f},\quad \nm{\bm{z}} \leq a \nm{f}.
}
For the TSVD we apply this and (\ref{TSVD_err}) to obtain 
\bes{
\nm{f-\mathcal{P}_{\Lambda}f} \leq  \nm{f-\mathcal{T}_{N} \bm{z}} + \sqrt{\epsilon} \nm{\bm{z}} \leq \left ( \delta + a \sqrt{\epsilon} \right ) \nm{f},
}
as required.
For the ASVD1, we first observe that
\bes{
\nm{\bm{y}} \geq \frac{\nm{f}^2 - \delta \nm{f}^2}{a \nm{f}} = \frac{1-\delta}{a} \nm{f},
}
by \R{(f,y)-NormConv}.  Since $c > \frac{a^2}{1-\delta}$, we have $\nm{\bm{z}} \leq a \nm{f} < c \nm{\bm{y}}$.
Hence, using (\ref{ASVD1_err}) we obtain
\bes{
\nm{f-\mathcal{P}_{\Lambda}f} \leq \delta \nm{f} + \max \left \{ \sqrt{\epsilon} , \frac{\delta}{c \nm{\bm{y}} - \nm{\bm{z}} } \right \} \nm{\bm{z}}.
}
We now use the bounds $\nm{\bm{z}} \leq a \nm{f}$ and $c \nm{\bm{y}}  > c (1-\delta) / a  \nm{f}$ to get the result.

Now let $\delta = b \sqrt{\epsilon}$ and consider ASVD2. Since $c > \frac{a(2a+b)}{1-b \sqrt{\epsilon}}$, we have
\bes{
\left( \frac{1}{\sqrt{\epsilon}}\right) \nm{f - \mathcal{T}_{N}\bm{z}} + 2 \nm{\bm{z}} <  \left( b + 2 a \right) \nm{f} \leq \frac{c \left( 1 - b \sqrt{\epsilon}\right)}{a} \nm{f} \leq c \nm{\bm{y}},
}
where, once again, the last inequality directly follows from \R{(f,y)-NormConv}. Thus, similarly to the first two parts, using (\ref{ASVD2_err}) one can write
\bes{
\nm{f-\mathcal{P}_{\Lambda}f} \leq b \sqrt{\epsilon} \nm{f} + \max \left\{ \sqrt{\epsilon}, \frac{b \sqrt{\epsilon} \nm{f}}{c \nm{\bm{y}} - \left( \frac{1}{\sqrt{\epsilon}}\right) \nm{f - \mathcal{T}_{N}\bm{z}} - 2 \nm{\bm{z}}} \right\} \nm{\bm{z}}.
}
By the bounds derived above, we get
\bes{
\nm{f-\mathcal{P}_{\Lambda}f} \leq b \sqrt{\epsilon} \nm{f} + \max \left\{ \sqrt{\epsilon}, \frac{b \sqrt{\epsilon} \nm{f}}{\frac{c\left( 1 - b \sqrt{\epsilon}\right)}{a} \nm{f} - b \nm{f} - 2a \nm{f}} \right \} a \nm{f}.
}
Simplifying and rearranging which finishes the proof.
}

\prf{[Proof of Theorem \ref{thm:limsup}]
For part \ref{i}, we use (\ref{coef_bound_lemma}) to write
\bes{
\nm{\bm{x}_{\Lambda}} \leq \left( 1/\min_{n \in \Lambda} \{ \sqrt{\sigma_{n}} \} \right) \nm{f - \mathcal{T}_{N}\bm{a}_{N}} + \nm{\bm{a}_{N}}.
}
We note that $\min_{n \in \Lambda}\{ \sigma_{n} \} {>} \epsilon$ in all three methods. Therefore
$
\nm{\bm{x}_{\Lambda}} \leq \frac{1}{\sqrt{\epsilon}} \nm{f - \mathcal{T}_{N}\bm{a}_{N}} + \nm{\bm{a}_{N}}.
$
Notice that $\nm{f - \mathcal{T}_{N}\bm{a}_{N}} \rightarrow 0$ and $\nm{\bm{a}_{N}} \rightarrow \nm{\bm{a}}$ as $N \rightarrow \infty$. Hence the result follows.

For part \ref{ii}, we use (\ref{coef_conv_lemma}) to write  
\bes{
\nm{\bm{a} - \bm{a}_{\Lambda}} \leq \left( 1 + \sqrt{B/\epsilon}\right) \sqrt{\sum_{n\in I\setminus I_{N}} \vert a_{n} \vert^{2}} + \left( \max_{n\in I_N \backslash \Lambda}{\left\{ \sqrt{\sigma_{n}} \right\}}/\sqrt{A} \right) \nm{\bm{a}}.
}
First, notice that the first term of the right-hand side vanishes as $N$ grows large, since $\sum_{n\in I\setminus I_{N}}\vert a_{n} \vert^{2} \rightarrow 0$ as $N \rightarrow \infty$. Now, consider the second term of the right-hand side. For the TSVD approximation, observe that $\max_{n\in I_N \backslash \Lambda}{\left\{ \sqrt{\sigma_{n}} \right\}} = \sqrt{\epsilon}$. Hence the result follows immediately. For either of the other two approximations, notice that $\nm{\bm{a}_N} \leq \nm{\bm{a}} \leq 1/\sqrt{A} \nm{f}$ and $\lim_{N \rightarrow \infty} \nm{\bm{y}} \geq \sqrt{A} \nm{f}$. For the ASVD1 approximation, since $c > 1/A$, we deduce that $c\nm{\bm{y}} > \nm{\bm{a}_{N}}$ for all sufficiently large values of $N$. Similarly, for the ASVD2 approximation, we deduce that $c\nm{\bm{y}} > \frac{1}{\sqrt{\epsilon}} \nm{f - \mathcal{T}_{N}\bm{a}_{N}} + 2\nm{\bm{a}_{N}}$ for all sufficiently large values of $N$, since $c > 2/A$ therein. Therefore, using (\ref{ASVD1_sing_final}) and {(\ref{ASVD2_sing_final})} we write
\bes{
\max_{n\in I_N \backslash \Lambda}{\left\{ \sqrt{\sigma_{n}} \right\}} \leq   \max{\left\{ \sqrt{\epsilon}, \frac{\nm{f - \mathcal{T}_{N}\bm{a}_{N}}}{c\nm{\bm{y}} - \nm{\bm{a}_{N}}}\right\}},
}
and
\bes{
\max_{n\in I_N \backslash \Lambda}{\left\{ \sqrt{\sigma_{n}} \right\}} \leq \max{\left\{ \sqrt{\epsilon}, \frac{\nm{f - \mathcal{T}_{N}\bm{a}_{N}}}{c\nm{\bm{y}} - \left(\frac{1}{\sqrt{\epsilon}}\right) \nm{f - \mathcal{T}_{N}\bm{a}_{N}} - 2\nm{\bm{a}_{N}}}\right\}},
}
for all sufficiently large values of $N$, forASVD1 and the ASVD2 respectively. Again, since $\nm{f - \mathcal{T}_{N}\bm{a}_{N}} \rightarrow 0$ and $\nm{\bm{a}_{N}} \rightarrow \nm{\bm{a}}$ as $N \rightarrow \infty$, we deduce that
\bes{
\limsup_{N \rightarrow \infty} \nm{\bm{a} - \bm{a}_{\Lambda}} \leq \sqrt{\epsilon/A}\nm{\bm{a}}.
}
For part \ref{iii}, we use (\ref{err_bound_lemma}) to write
\bes{
\nm{f - \mathcal{P}_{\Lambda}f} \leq \nm{f - \mathcal{T}_{N}\bm{a}_{N}} + \max_{n\in I_N \backslash \Lambda}{\left\{ \sqrt{\sigma_{n}} \right\}} \nm{\bm{a}_{N}}, \quad \forall N \in I_{N},
}
and with exact the same argument as in proof of part \ref{ii} we deduce that
\bes{
\limsup_{N \rightarrow\infty} \nm{f - \mathcal{P}_{\Lambda}f} \leq \sqrt{\epsilon} \nm{\bm{a}}, 
}
for either of the TSVD, ASVD1, or ASVD2 approximation.
}

\subsection*{Acknowledgements}
The question of frame approximation with bounded coefficients was first raised during a talk by the first author at the Oberwolfach conference on ``Multiscale and High-Dimensional Problems''.  The authors would like to thank Ingrid Daubechies for raising this question. They would also like to thank Daan Huybrechs for helpful comments and suggestions. This work was supported by NSERC  through grant 611675, as well as through the PIMS CRG on ``High-dimensional Data Analysis''.

\bibliographystyle{abbrv}
\small
\bibliography{FramesBCrefs}

\begin{thebibliography}{10}

\bibitem{PEHighDim}
B.~Adcock and D.~Huybrechs.
\newblock Approximating smooth, multivariate functions on irregular domains.
\newblock {\em arXiv:1802.00602}, 2018.

\bibitem{BADHFramesPart2}
B.~Adcock and D.~Huybrechs.
\newblock Frames and numerical approximation {II}: generalized sampling.
\newblock {\em arXiv:1802.01950}, 2018.

\bibitem{BADHframespart}
B.~Adcock and D.~Huybrechs.
\newblock Frames and numerical approximation.
\newblock {\em SIAM Rev.}, 61(3):443--473, 2019.

\bibitem{boffi2015immersed}
D.~Boffi, N.~Cavallini, and L.~Gastaldi.
\newblock The finite element immersed boundary method with distributed
  {L}agrange multiplier.
\newblock {\em SIAM J. Numer. Anal.}, 53(6):2584--2604, 2015.

\bibitem{boyd2005fourier}
J.~Boyd.
\newblock Fourier embedded domain methods: extending a function defined on an
  irregular region to a rectangle so that the extension is spatially periodic
  and ${C}^{\infty}$.
\newblock {\em Appl. Math. Comput.}, 161(2):591--597, 2005.

\bibitem{christensen2003introduction}
O.~Christensen.
\newblock {\em An Introduction to Frames and {R}iesz Bases}.
\newblock Applied and Numerical Harmonic Analysis. Birkh{\"a}user, 2nd edition,
  2016.

\bibitem{kasolis2015fictitious}
F.~Kasolis, E.~Wadbro, and M.~Berggren.
\newblock Analysis of fictitious domain approximations of hard scatterers.
\newblock {\em SIAM J. Numer. Anal.}, 2015(5):2347--2362, 2015.

\bibitem{lindner2006}
M.~Lindner.
\newblock {\em Infinite Matrices and their Finite Sections}.
\newblock Frontiers in Mathematics. Birkh\"auser Verlag, Basel, 2006.

\bibitem{lui2009embedding}
S.~H. Lui.
\newblock Spectral domain embedding for elliptic {PDE}s in complex domains.
\newblock {\em J. Comput. Appl. Math.}, 225(2):541--557, 2009.

\bibitem{LyonFast}
M.~Lyon.
\newblock A fast algorithm for {F}ourier continuation.
\newblock {\em {SIAM} J. Sci. Comput.}, 33(6):3241--3260, 2012.

\bibitem{matthysen2015fastfe}
R.~Matthysen and D.~Huybrechs.
\newblock Fast algorithms for the computation of {F}ourier extensions of
  arbitrary length.
\newblock {\em SIAM J. Sci. Comput.}, 38(2):A899--A922, 2016.

\bibitem{matthysen2017fastfe2d}
R.~Matthysen and D.~Huybrechs.
\newblock {Function approximation on arbitrary domains using Fourier extension
  frames}.
\newblock {\em SIAM J. Numer. Anal.}, 56(3):1360--1385, 2018.

\bibitem{pasquettiFourEmbed}
R.~Pasquetti and M.~Elghaoui.
\newblock A spectral embedding method applied to the advection--diffusion
  equation.
\newblock {\em J. Comput. Phys.}, 125:464--476, 1996.

\bibitem{shirokoff2015volumepenalty}
D.~Shirokoff and J.-C. Nave.
\newblock A sharp-interface active penalty method for the incompressible
  {N}avier--{S}tokes equations.
\newblock {\em J. Sci. Comput.}, 62(1):53--77, 2015.

\end{thebibliography}

\end{document}